\documentclass[12pt,a4paper]{amsart}
\usepackage{mathtools}
\usepackage{xcolor}
\usepackage{amsmath,amstext,amssymb,amsthm,amsfonts}
\usepackage[all]{xy} \xyoption{poly}

\usepackage{graphicx}

\usepackage{hyperref}
\hypersetup{colorlinks,linkcolor={blue},citecolor={blue},urlcolor={red}}

\usepackage{framed}

\usepackage{setspace}

\usepackage{tikz}
\usepackage{pgffor}
\usepackage[framemethod=tikz]{mdframed}

\newcommand{\nc}{\newcommand}
\nc{\one}{\mbox{\bf 1}}
\nc{\invtensor}{\underset{\leftarrow}{\otimes}}
\nc{\const}{\operatorname{const}}

\nc{\ad}{\operatorname{ad}}

\nc{\tr}{\operatorname{tr}}
\nc{\dex}{\operatorname{dex}}
\nc{\pari}{\operatorname{dex}}
\nc{\tp}{\operatorname{top}}
\nc{\rank}{\operatorname{rank}}
\nc{\corank}{\operatorname{corank}}
\nc{\codim}{\operatorname{codim}}
\nc{\sdim}{\operatorname{sdim}}
\nc{\mult}{\operatorname{mult}}
\nc{\ds}{\operatorname{ds}}
\nc{\tail}{\operatorname{tail}}
\nc{\howl}{\operatorname{howl}}
\nc{\spn}{\operatorname{span}}
\nc{\Sym}{\operatorname{Sym}}
\nc{\sym}{\operatorname{sym}}
\nc{\id}{\operatorname{id}}
\nc{\Id}{\operatorname{Id}}
\nc{\Ree}{\operatorname{Re}}
\nc{\hi}{\operatorname{hi}}
\nc{\htt}{\operatorname{ht}}
\nc{\at}{\operatorname{at}}
\nc{\str}{\operatorname{str}}
\nc{\Iso}{\operatorname{Iso}}
\nc{\Ker}{\operatorname{Ker}}
\nc{\rker}{\operatorname{rKer}}
\nc{\im}{\operatorname{Im}}
\nc{\osp}{\mathfrak{osp}}
\nc{\sgn}{\operatorname{sgn}}
\nc{\F}{\operatorname{F}}
\nc{\Mod}{\operatorname{Mod}}
\nc{\DS}{\operatorname{DS}}
\nc{\Soc}{\operatorname{Soc}}
\nc{\Inj}{\operatorname{Inj}}
\nc{\Hom}{\operatorname{Hom}}
\nc{\End}{\operatorname{End}}
\nc{\supp}{\operatorname{supp}}
\nc{\Card}{\operatorname{Card}}
\nc{\Ann}{\operatorname{Ann}}
\nc{\Ind}{\operatorname{Ind}}
\nc{\Coind}{\operatorname{Coind}}
\nc{\wt}{\operatorname{wt}}
\nc{\ch}{\operatorname{ch}}
\nc{\sch}{\operatorname{sch}}
\nc{\mdim}{\operatorname{mdim}}
\nc{\Stab}{\operatorname{Stab}}
\nc{\Irr}{\operatorname{Irr}}
\nc{\Spec}{\operatorname{Spec}}
\nc{\Res}{\operatorname{Res}}
\nc{\Aut}{\operatorname{Aut}}
\nc{\Ext}{\operatorname{Ext}}
\nc{\Prec}{\operatorname{Prec}}
\nc{\Fract}{\operatorname{Fract}}
\nc{\gr}{\operatorname{gr}}
\nc{\deff}{\operatorname{def}}
\nc{\core}{\operatorname{core}}
\nc{\HC}{\operatorname{HC}}
\nc{\dpth}{\operatorname{dpth}}
\nc{\sw}{\operatorname{sw}}
\nc{\red}{\operatorname{red}}

\nc{\pos}{\operatorname{pos}}
\nc{\atyp}{\operatorname{atyp}}
\nc{\wdchi}{\widetilde{\chi}}
\nc{\wdH}{\widetilde{H}}
\nc{\wdN}{\widetilde{N}}
\nc{\wdM}{\widetilde{M}}
\nc{\wdO}{\widetilde{O}}
\nc{\wdR}{\widetilde{R}}

\nc{\wdV}{\widetilde{V}}

\nc{\wdC}{\widetilde{C}}

\nc{\Obj}{\operatorname{Obj}}
\nc{\Dglie}{\operatorname{{\mathcal D}glie}}
\nc{\Fin}{\operatorname{{\mathcal F}in}}
\nc{\pr}{\operatorname{pr}}
\nc{\Adm}{\operatorname{\mathcal{A}dm}}
\nc{\Sg}{{\cS(\fg)}}
\nc{\Shg}{{\cS(\fhg)}}
\nc{\Ug}{{\cU(\fg)}}
\nc{\Uhg}{{\cU(\fhg)}}
\nc{\Sh}{{\cS(\fh)}}
\nc{\Uh}{{\cU(\fh)}}
\nc{\Uhh}{{\cU(\fhh)}}
\nc{\Zg}{{{\mathcal{Z}}(\fg)}}

\nc{\Vir}{{\mathcal{V}ir}}
\nc{\NS}{{\mathcal{N}S}}

\nc{\tZg}{{\widetilde{\mathcal Z}({\mathfrak g})}}
\nc{\Zk}{{\mathcal Z}({\mathfrak k})}

\newcommand{\shr}{shr}

\nc{\Up}{{\mathcal U}({\mathfrak p})}
\nc{\Ah}{{\mathcal A}({\mathfrak h})}
\nc{\Ag}{{\mathcal A}({\mathfrak g})}
\nc{\Ap}{{\mathcal A}({\mathfrak p})}
\nc{\Zp}{{\mathcal Z}({\mathfrak p})}
\nc{\cR}{\mathcal R}
\nc{\cS}{\mathcal S}
\nc{\cT}{\mathcal{T}}
\nc{\CC}{\mathcal C}
\nc{\cA}{\mathcal A}
\nc{\cU}{\mathcal U}
\nc{\cZ}{\mathcal Z}
\nc{\cM}{\mathcal M}
\nc{\cL}{\mathcal L}
\nc{\cF}{\mathcal F}
\nc{\fg}{\mathfrak g}
\nc{\cB}{\mathcal{B}}

\nc{\fo}{\mathfrak o}

\nc{\CO}{\mathcal O}
\nc{\CR}{\mathcal R}

\nc{\cK}{\mathcal{K}}
\nc{\cW}{\mathcal{W}}
\nc{\bM}{\mathbf{M}}
\nc{\bL}{\mathbf{L}}
\nc{\bN}{\mathbf{N}}

\nc{\zq}{\mathpzc q}

\nc{\fl}{\mathfrak l}
\nc{\fn}{\mathfrak n}
\nc{\fm}{\mathfrak m}
\nc{\fp}{\mathfrak p}
\nc{\fh}{\mathfrak h}
\nc{\ft}{\mathfrak t}
\nc{\fk}{\mathfrak k}
\nc{\fb}{\mathfrak b}
\nc{\fs}{\mathfrak s}
\nc{\fB}{\mathfrak B}

\nc{\vareps}{\varepsilon}
\nc{\varesp}{\varepsilon}
\nc{\veps}{\varepsilon}

\nc{\fsl}{\mathfrak{sl}}
\nc{\fgl}{\mathfrak{gl}}
\nc{\fso}{\mathfrak{so}}
\nc{\fosp}{\mathfrak{osp}}
\nc{\fsp}{\mathfrak{sp}}
\nc{\fq}{\mathfrak q}
\nc{\fsq}{\mathfrak{sq}}
\nc{\fpsl}{\mathfrak{psl}}

\nc{\fhg}{\hat{\fg}}
\nc{\fhn}{\hat{\fn}}
\nc{\fhh}{\hat{\fh}}
\nc{\fhb}{\hat{\fb}}
\nc{\hrho}{\hat{\rho}}

\nc{\hsl}{\hat{\fsl}}
\nc{\fpo}{\mathfrak{po}}
\nc{\dirlim}{\underset{\rightarrow}{\lim}\,}
\nc{\nen}{\newenvironment}
\nc{\ol}{\overline}
\nc{\ul}{\underline}
\nc{\ra}{\rightarrow}
\nc{\lra}{\longrightarrow}
\nc{\Lra}{\Longrightarrow}
\nc{\bo}{\bar{1}}
\nc{\Lla}{\Longleftarrow}

\nc{\Llra}{\Longleftrightarrow}

\nc{\thla}{\twoheadleftarrow}

\nc{\lang}{(}
\nc{\rang}{)}

\nc{\hra}{\hookrightarrow}

\nc{\iso}{\overset{\sim}{\lra}}

\nc{\ssubset}{\underset{\not=}{\subset}}

\nc{\vac}{|0\rangle}

\nc{\simka}{{\ \scriptscriptstyle _{{\sim}}^\text{\tiny{k}}\ }}

\newcommand{\pla}{\ulcorner{\!\lambda}}

\makeatletter
\newcommand*{\owedge}{%
  \mathbin{%
    \mathpalette\@owedge{}%
  }%
}
\newcommand*{\@owedge}[2]{%
  \sbox0{$#1\oplus\m@th$}%
  \dimen2=.5\dimexpr\wd0-\ht0-\dp0\relax % side bearing
  \dimen@=\dimexpr\ht0+\dp0\relax
  \def\lw{.04}% line width as factor for height of \oplus
  \def\radius{.5-\lw/2}%
  \kern\dimen2 % side bearing
  \tikz[
    line width=\lw\dimen@,
    line join=round,
    x=\dimen@,
    y=\dimen@,
    baseline=\dimexpr-.5\dimen@+\dp0\relax,
  ]
  \draw
    (0,0) circle[radius=\radius]
    (225:\radius) -- (0,.5-\lw) -- (-45:\radius)
  ;%
  \kern\dimen2 % side bearing  
}
\makeatother
\nc{\Thm}[1]{Theorem~\ref{#1}}
\nc{\Prop}[1]{Proposition~\ref{#1}}
\nc{\Lem}[1]{Lemma~\ref{#1}}
\nc{\Cor}[1]{Corollary~\ref{#1}}
\nc{\Conj}[1]{Conjecture~\ref{#1}}
\nc{\Claim}[1]{Claim~\ref{#1}}
\nc{\Defn}[1]{Definition~\ref{#1}}
\nc{\Exa}[1]{Example~\ref{#1}}
\nc{\Rem}[1]{Remark~\ref{#1}}
\nc{\Note}[1]{Note~\ref{#1}}
\nc{\Quest}[1]{Question~\ref{#1}}
\nc{\Hyp}[1]{Hypoth\`ese~\ref{#1}}
\nen{thm}[1]{\label{#1}{\bf Theorem.\ } \em}{}
\nen{prop}[1]{\label{#1}{\bf Proposition.\ } \em}{}
\nen{lem}[1]{\label{#1}{\bf Lemma.\ } \em}{}
\nen{cor}[1]{\label{#1}{\bf Corollary.\ } \em}{}
\nen{conj}[1]{\label{#1}{\bf Conjecture.\ } \em}{}

\nen{claim}[1]{\label{#1}{\bf Claim.\ } \em}{}

\nen{defn}[1]{\label{#1}{\bf Definition.\ } }{}
\nen{exa}[1]{\label{#1}{\bf Example.\ } }{}
\nen{rem}[1]{\label{#1}{\em Remark.\ } }{}
\nen{note}[1]{\label{#1}{\em Note.\ } }{}
\nen{exer}[1]{\label{#1}{\em Exercise.\ } }{}
\nen{sket}[1]{\label{#1}{\em Sketch of proof.\ } }{}
\nen{quest}[1]{\label{#1}{\bf Question.\ } \em}{}

\nen{hyp}[1]{\label{#1}{\bf Hypoth\`ese.\ } \em}{}
\begin{document}
\setcounter{section}{0}
\setcounter{tocdepth}{1}

\title[Purity]{Semisimplicity of the $DS$ functor for the orthosymplectic Lie superalgebra}
\author{{\rm M. Gorelik, Th. Heidersdorf}}

\address{T. H.: Mathematisches Institut Universit\"at Bonn, Germany}
\email{heidersdorf.thorsten@gmail.com} 
\address{M. G..: Weizmann Institute of Science, Rehovot, Israel}
\email{maria.gorelik@gmail.com}

\date{}

\begin{abstract} We prove that the Duflo-Serganova functor $DS_x$ attached to an odd nilpotent element $x$ of $\mathfrak{osp}(m|2n)$ is semisimple, i.e. sends a semisimple representation $M$ of $\mathfrak{osp}(m|2n)$ to a semisimple representation of $\mathfrak{osp}(m-2k|2n-2k)$ where $k$ is the rank of $x$. We prove a closed formula for $DS_x(L(\lambda))$ in terms of the arc diagram attached to $\lambda$. 
\end{abstract}

\subjclass[2010]{17B10, 17B20, 17B55, 18D10.}

\medskip

\keywords{Representations of supergroups, Purity, Duflo-Serganova functor, Superdimensions, Orthosymplectic Lie superalgebra, Orthosymplectic supergroup, Symmetric monoidal functor}

\maketitle

%\renewcommand{\baselinestretch}{0.3}\normalsize
%\tableofcontents
%\renewcommand{\baselinestretch}{1.0}\normalsize

%\begin{spacing}{0.05}
%\tableofcontents
%\listoffigures
%\listoftables
%\end{spacing}

%\tableofcontents

\section{Introduction}

\subsection{The Duflo-Serganova functor $DS$} For a finite dimensional complex Lie superalgebra $\mathfrak{g}$ and an odd nilpotent element $x$ Duflo and Serganova defined a functor $DS_x: Rep(\mathfrak{g}) \to Rep(\mathfrak{g}_x)$ where $\mathfrak{g}_x := ker \ ad(x) / im \ ad(x)$. This functor is given by taking the cohomology of the complex $\xymatrix{ V \ar[r]^{\rho(x)} & V \ar[r]^{\rho(x)} & V} $ for $(V,\rho) \in Rep(\mathfrak{g})$. For $\mathfrak{gl}(m|n)$ we have $\mathfrak{g}_x \simeq \mathfrak{gl}(m-k|n-k)$ and for $\mathfrak{osp}(m|2n)$ we have $\mathfrak{g}_x \simeq \mathfrak{osp}(m-2k|2n-2k)$ where $k$ is the so-called rank of $x$. The $DS$-functor is a symmetric monoidal functor which allows to reduce questions about superdimensions or tensor products to lower rank. It is however very complicated to compute $DS(V)$ explicitly. In \cite{HW} the authors derived a closed formula  for $DS_x(L(\lambda))$ where $L(\lambda))$ is a finite dimensional irreducible representation of $\mathfrak{gl}(m|n)$ and arbitrary $x$. In particular $DS_x(L(\lambda))$ is always semisimple. %The $L(\lambda_i)$ are obtained from $L(\lambda)$ by the removal of certain \emph{maximal arcs or cups} in the cup diagram of $\lambda$. %This formula is first proven for $x$ of rank $1$, and the general case is deduced by noting via some spectral sequence argument that $DS_x(L(\lambda))$ is the same as $DS_1 \circ \ldots \circ DS_1(L(\lambda))$. This result has a number of important applications. In particular it gives a formula for the sudperdimension of an irreducible representation and also allows to reduce certain questions about tensor products to lower rank \cite{H-semisimple}\cite{HW-tannaka}.

\subsection{Serganova's conjecture} More generally Serganova conjectured \cite{S-ICM} that $DS_x$ should be semisimple for any basic classical Lie superalgebra (i.e. $DS_x(L)$ is semisimple for irreducible $L$). Since the conjecture is trivial in the exceptional cases (in these cases $\mathfrak{g}_x$ is a reductive Lie algebra) this leaves the $\mathfrak{osp}(m|2n)$-case. We prove Serganova's conjecture in this article and give  a closed formula for $DS_x(L(\lambda))$ for any $x$ and any $\lambda$. %We remark that $DS_x(L(\lambda))$ is not semisimple in the periplectic case \cite{IEAS}.

\subsection{Main steps and results} 

%\subsubsection{The $DS$ functor} 
We work in the full subcategory $\tilde{\cF}(\fg)$ of algebraic representations of $\mathfrak{osp}(m|2n)$. The category $\tilde{\cF}(\fg)$  is canonically
isomorphic to the category of $SOSp(m|2n)$-modules; this
category decomposes into a direct sum:
 $$\tilde{\cF}(\fg) = {\cF}(\fg) \oplus \Pi({\cF}(\fg)),$$
where $\cF(\fg)$ is the full subcategory of modules $M$, where the parity is induced
from the parity on $\Lambda_{m|n}$. We view $DS_x$ as a functor $\tilde{\cF}(\fg) \to \tilde{\cF}(\fg_x)$. For fixed $x$ of rank $r$ we also write $DS_r$. 

\subsubsection{The cases $t=0,1,2$} Recall that $\osp(M|N)$ consists
of two series: the $B$-series (when $M$ is odd) and 
the $D$-series (when $M$ is even). In the $\fosp$-case each block of atypicality $k$ is equivalent to the principal block of $\fosp(2k+t|2k)$ where $t= 1$ for the B-series and $t= 0,2$ for the D-series.  The equivalences are described in \cite{GS} (we recall
details in~\ref{notres}). We consider these three cases separately and set $\fg:=\osp(2m+t|2n)$, where $t=0,1,2$ as above. This division
is compatible with $\DS$-functors: if $N$ lies in a block of type $t$,
then $\DS_x(N)$  lies in a block of type $t$.
In Section~\ref{Fgp} we reduce the computations of multiplicities $[DS_x(L(\lambda)):L_{\mathfrak{g}_x}(\nu)]$ to
the case of principal blocks.

\subsubsection{Reduction to the principal block} Any irreducible module can be moved via a series of translation functors to a \emph{stable} irreducible module (see section \ref{sect2}). The subcategory $\cF^{g}_{st}$ of a block $\cF^g(\fg)$ consisting
 of stable modules (a module {\em  stable} if all its simple subquotients are stable) is equivalent via a functor $Res$ to the principal block of $\fosp(2k+t|2k)$ for
$t=0,1,2$. Since $DS$ commutes with both $Res$ and translation functors (see section \ref{DSResx}), we can compute it on irreducible modules in the principal block.

\subsubsection{Recursive formulae}
We then induct on the degree of atypicality. For atypicality 1 (i.e. the principal block of $\mathfrak{osp}(2|2)$, $\mathfrak{osp}(3|2)$ and $\mathfrak{osp}(4|2)$) $DS(L(\lambda))$ can be computed easily 
(see for example \cite{Gdex}). In order to treat the general case we establish recursive formulae for the multiplicities of irreducible constituents in $DS_x(L(\lambda))$.

For a category $\CC$ denote by $\Irr(\CC)$ the set of isomorphism
classes of simple modules in $\CC$. Let $V_{st}$ be the natural representation.
For blocks $\tilde{\cF}^{g_1}(\fg), \tilde{\cF}^{g_0}(\fg)$ we denote
by $T^{g_0}_{g_1}$ the translation functor
%$$T^{g'}_g:\ \cF^{g}(\fg)\to\cF^{g'}(\fg),\ \ \ \ \
$$T^{g_0}_{g_1}:\ \tilde{\cF}^{g_1}(\fg)\to\tilde{\cF}^{g_0}(\fg)$$
which maps $N$ to the projection of $N\otimes V_{st}$ to the subcategory
$\tilde{\cF}^{g_0}(\fg)$. Since blocks of atypicality $n$ for $\fg$ and blocks of atypicality $n-r$ for $\fg_x$ ($rk(x)  = r)$ correspond to each other (via so-called core diagrams), we may look at 
$$T^{g_0}_{g_1}: \tilde{\cF}(\fg)^{g_1}\to \tilde{\cF}(\fg)^{g_0},\ \ \
T^{g_0}_{g_1}: \tilde{\cF}(\fg_x)^{g_1}\to \tilde{\cF}(\fg_x)^{g_0}.$$

For  $N\in\cF(\fg)^{g_1}$ and
 $L'\in \Irr(\fg_x)^{g_0}$ we have
\begin{align*}\begin{array}{ll}
[\DS_r(T^{g_0}_{g_1}(N)):L']&=[T^{g_0}_{g_1}(DS_r(N)):L']\\
&=\displaystyle\sum_{L_1\in\Irr(\fg_x)^{g_1}}
[DS_r(N):L_1][T^{g_0}_{g_1}(L_1):L'].\end{array}\end{align*}

This formula allows us in Section~\ref{sect4}  to successively reduce the computation of the multiplicity $[DS_x(L(\lambda)):L_{\mathfrak{g}_x}(\nu)]$ to the case where $\fg_x$ is $0$ or $\mathbb{C}$.

\subsubsection{Multiplicities for the case $\rank x=1$}

For such $x$ the condition $\fg_x=0,\mathbb{C}$ gives 
$\fg=\osp(m|2)$ for $m=2,3,4$; for these cases $\DS_x(L)$ can be easily computed.
This gives the  multiplicity $[DS_x(L(\lambda)):L_{\mathfrak{g}_x}(\nu)]$
for   $\rank x=1$.  As in the $\mathfrak{gl}(m|n)$ and $\fp_n$ cases, treated in \cite{HW},
\cite{EAS} respectively, we give the final answer in terms of \emph{arc diagrams}.

For a weight $\lambda$ denote by $\howl(\lambda)$ the corresponding weight in the principal block; to $\howl(\lambda)$ we attach an arc diagram $Arc(\howl(\lambda))$, see \ref{arcs}.

\medskip
\textbf{Theorem A} (see Theorem \ref{thmDS1osp} for a more precise statement). \emph{Let $rk(x) = 1$. Then $L(\nu)$ is a subquotient of $DS(L(\lambda))$ if and only if
$Arc(\howl(\nu))$ is obtained from $Arc(\howl(\lambda))$ by removing
a maximal arc and, in addition,
in the $\osp(2m+1|2n)$-case,
if $\nu$ has a sign, then the signs of $\lambda$ and $\nu$ are equal. The multiplicity of $L(\nu)$ is either 1 or 2 (depending on the shape of the Arc diagram).}
\medskip

\subsubsection{Semisimplicity for the case $\rank x=1$}
In~\cite{Gdex} (see also~\Cor{corext}) it is shown that
 the extension graph of $\cF(\fg)$ is bipartite; the bipartition
 is given by  a sign function $\dex: \Irr(\tilde{\cF}(\fg)) \to \{\pm 1\}$ such that $Ext^1(L_1,L_2)= 0 \ \text{ if } \dex(L_1) = \dex(L_2)$. 
Theorem A implies $[\DS_x(L(\lambda)):L_{\mathfrak{g}_x}(\nu)]=0$ if
 $\dex(L(\lambda))\not=\dex(L(\nu))$ (and $\rank x=1$); this shows the
semisimplity of $\DS_x(L(\lambda))$.

\subsection{The general case}
We denote with $\cF_+(\fg)$ the full Serre subcategory of $\tilde{\cF}(\fg)$ generated by the irreducible objects with $\dex(L) = 1$. By definition it is a semisimple category. We say that a module $M$ is \emph{pure} if for any subquotient $L$ of $M$, $\Pi(L)$ is not a subquotient of $M$.

Serganova's semisimplicity conjecture follows from the following theorem.

\medskip
\textbf{Theorem B} (Purity and semisimplicity).
\emph{For each $x$ one has  $\DS_x\bigl(\cF_+(\fg)\bigr)=\cF_+(\fg_x)$. In particular $DS_x(L(\lambda))$ is pure.}

 {\em For each  $L\in\Irr(\tilde{\cF}(\fg))$ the module $\DS_x(L)$ only depends on the rank of $x$. Moreover,}
 $\DS_{r+1}(L)\cong \DS_{1}(\DS_{r}(L))$, {\em where }$\DS_r$ {\em stands for }$\DS_x$ {\em with }$\rank x=r$.
\medskip

By above, the statements follow from Theorem A for $\rank x=1$;
the general case easily follows by induction on
$\rank x$.

For $\fgl(m|n)$ purity and semisimplicity were established
in \cite{HW}; for the exceptional cases purity was checked in \cite{Gaugusta} (and semisimplicity is trivial since $\fg_x$ is reductive). For $\fp_n$
the semisimplicity does not hold ($\DS_1$ maps the standard representation
of $\fp_2$ to  the standard representation
of $\fp_{1}$ which is not semisimple).

Theorems A,B allow us to compute the superdimension of an irreducible $\mathfrak{osp}(m|2n)$-module (since $\DS$ is symmetric monoidal it preserves the superdimensions); see \cite{HW} and \cite{EAS} for analogous results in the $\mathfrak{gl}(m|n)$ and $\mathfrak{p}(n)$-case. Our main theorem also allows to reduce certain questions about tensor products to lower rank similarly to \cite{HW-tannaka}, \cite{Heidersdorf-semisimple}.

\subsection{Acknowledgments} The authors are grateful to  V.~Hinich,  V.~Serganova and C.~Stroppel for numerous 
helpful discussions.

\subsection{Index of definitions and notation} \label{sec:app-index}
Throughout the paper the ground field is $\mathbb{C}$; 
$\mathbb{N}$ stands 
 for the set of non-negative integers. We denote by $\Pi$ the parity change functor. 
In Sections \ref{sectPrelim}---\ref{Sectlast}
$\fg=\osp(M|2n)$ with $M,n\geq 0$ and $M$ either $2m$ or $2m+1$. 
One has $\osp(1|0)=\osp(0|0)=0$ and $\osp(2|0)=\mathbb{C}$.

\begin{center}
\begin{tabular}{ll}
$\Omega(N), \Sigma, \sigma, \tilde{\cF}(\fg), \cF(\fg), \Irr, \Lambda^+_{m|n}, 
L(\lambda),\ \chi_{\lambda}\ \ $  & Section~\ref{sectPrelim}\\

weight diagram, the diagram $f_-f_+$ & \ref{weightdia}\\

core symbols, core diagram & \ref{corewt}\\

$\cF^g(\fg)$ & \ref{catFgg} \\

$t$, $\ell$, $\Lambda^{(t)}_{m+\ell|n}$, core-free, $\Theta_k^{(t)}\ \ $ & \ref{reft}\\

$\fg_r, \Sigma_r, S_r$  & \ref{Sigmar}\\

$\tail$ & \ref{tailsub}\\

$\howl, \tau$ & \ref{howl}\\

$\dex$ & \ref{mapdex}\\

stable diagrams & \ref{sec:stable}\\

translation functor $T^{g'}_g$ & \ref{Transfunc}\\

$\DS_x$ & \ref{DSconst}\\

$\supp(x)$, $x_s$, $\DS_s$ & \ref{G0} \\

$\fg_x$ & \ref{subsecgx} \\

graded multiplicity & \ref{gradedmult}\\

arc diagram & \ref{arcs}

\end{tabular}
\end{center}

\section{Notation}\label{sectPrelim}

\subsection{Root lattice}\label{notat}
Our notation and a choice of triangular decomposition follow \cite{GS},\cite{GSBGG}.

We denote by $\Delta$ the set of roots of $\fg$ and by $\Delta_{0}$
(resp., $\Delta_{1}$) the set of even (resp., odd ) roots.
In this paper all modules are weight modules (i.e., $\fh$ acts diagonally) with
finite-dimensional weight spaces; for such a module $N$ we set
$$\Omega(N):=\{\nu\in\fh^*|\ N_{\nu}\not=0\}.$$

The root system $\Delta$ lies in the lattice $\Lambda_{m|n}\subset \fh^*$ spanned by $\{\vareps_i\}_{i=1}^m\cup\{\delta_i\}_{i=1}^n$. We denote by $\Lambda$ the lattice
spanned by $\{\vareps_i\}_{i=1}^{\infty}\cup\{\delta_i\}_{i=1}^{\infty}$
and view $\Lambda_{m|n}$ as a subset of $\Lambda$. We define the parity homomorphism
$p:\Lambda\to\mathbb{Z}_2$ by
$p(\vareps_i)=\ol{0}$, $p(\delta_j)=\ol{1}$ for all $i,j$.

\subsection{Triangular decomposition}\label{tri}
We fix a Cartan subalgebra $\fh\subset\fg_{\ol{0}}$.
We fix a triangular decomposition corresponding to 
a "mixed" base $\Sigma$, i.e. a base
 containing maximal possible number of odd roots.
For $\osp(2m+1|2n)$  with $m,n>0$ we take
$$\Sigma:=\left\{\begin{array}{ll}
\delta_1-\delta_2,\delta_2-\delta_3,\ldots,\delta_{n-m}-\vareps_1,\vareps_1-\delta_{n-m+1},\ldots,\vareps_{m}-\delta_{n},
\delta_n\ & \text{ for } n>m\\
\vareps_1-\vareps_2,\vareps_2-\vareps_3,\ldots, \vareps_{m-n+1}-\delta_1,\delta_1-\vareps_{m-n+2},\ldots,\vareps_m-\delta_n,\delta_n\ & \text{ for }
m\geq n\\
\end{array}\right.$$
and for  $\osp(2m|2n)$ with $m,n>0$
we take
$$\Sigma:=\left\{\begin{array}{ll}
\delta_1-\delta_2,\delta_2-\delta_3,\ldots,\delta_{n-m+1}-\vareps_1,\vareps_1-\delta_{n-m+2},\ldots,\vareps_{m-1}-\delta_{n},
\delta_n\pm\vareps_m & \text{ for } n\geq m\\
\vareps_1-\vareps_2,\vareps_2-\vareps_3,\ldots, \vareps_{m-n}-\delta_1,\delta_1-\vareps_{m-n+1},\ldots,\vareps_{m-1}-\delta_n,\delta_n\pm\vareps_m & \text{ for }
m >n\\
\end{array}\right.$$
For the remaining case  $mn=0$ all triangular decompositions are conjugate and we fix a standard base. 

Recall that  $\osp(1|0)=\osp(0|0)=0$ and $\osp(2|0)=\mathbb{C}$;
in these cases $\fh=\fg$ and $\Delta=\Sigma=\emptyset$.
Note that $\Lambda_{0|0}=0$ and $\Lambda_{1|0}=\mathbb{Z}\vareps_1$.

%We always consider $\Sigma$ as
%the ordered set with respect to the above order.
We denote by $\rho$ the Weyl vector of $\fg$.

\subsubsection{Involution $\sigma$}
The superalgebra $\osp(2m|2n)$ with  $m,n>0$ or $m>1,n=0$ admits an involutive automorphism $\sigma$
induced by the automorphism of the Dynkin diagram of $\Sigma$
(the resulting action on $\fh^*$ is given by the reflection
$r_{\vareps_m}$). The restriction of $\sigma$ to
$\mathfrak{o}_{2m}\subset\fg_0$ is the standard involution
induced by the automorphism of the Dynkin diagram of $\mathfrak{o}_{2m}$.
 For  $\osp(2|0)=\mathbb{C}$
we set $\sigma:=-\Id$. For $\osp(0|2n)=\mathfrak{sp}_{2n}$ we set
$\sigma:=\Id$. 
Note that for all cases $\sigma(\Sigma)=\Sigma$.

\subsection{Category $\cF(\fg)$}\label{tildeF}
The category $\Fin(\fg)$ of finite dimensional representations of $\mathfrak{g}$ with parity preserving morphisms is the direct sum of two categories: $\tilde{\cF}(\fg)$
with the modules whose weights lie in $\Lambda_{m|n}$ and   $\tilde{\cF}^{\perp}(\fg)$
with the modules whose weights lie in $\fh^*\setminus \Lambda_{m|n}$.
The category $\tilde{\cF}^{\perp}(\fg)$
 is semisimple and $\DS_x(\tilde{\cF}^{\perp}(\fg))=0$ for $x\not=0$;
 all simple modules  in  $\tilde{\cF}^{\perp}(\fg)$
 are typical,
their characters  are given by the Weyl-Kac character formula.

The category $\tilde{\cF}(\fg)$  is canonically
isomorphic to the category of $SOSp(m|n)$-modules; this
category decomposes into a direct sum:
 $$\tilde{\cF}(\fg) = {\cF}(\fg) \oplus \Pi({\cF}(\fg)),$$
where $\cF(\fg)$ is the full subcategory with the modules $M$, where the parity is induced
from the parity on $\Lambda_{m|n}$, i.e. $M\in\cF(\fg)$
if and only if each weight space $p(M_{\nu})=p(\nu)$ for all $\nu\in \Lambda_{m|n}$.
We will sometimes omit $\fg$ from notation if this  does not lead to ambiguity; for instance,
we may use $\cF$ instead of $\cF(\fg)$.

\subsubsection{}\label{d0d1}
For each category $\CC$ we denote by $\Irr(\CC)$
 the set of isomorphism
classes of irreducible modules in $\CC$. For each finite-dimensional 
$N\in\CC$ and $L\in\Irr(\CC)$ we consider the {\em graded multiplicity}:
we write $[N:L]=(d_0|d_1)$ if a Jordan-H\"older series
of $N$ contains $d_0$ copies of $L$ and
$d_1$ copies of $\Pi(L)$.

We denote by  $L(\lambda)$ 
a simple $\fg$-module  of the highest weight $\lambda$ 
We set
$$\Lambda^+_{m|n}:=\{\lambda\in \Lambda_{m|n}|\ \dim L(\lambda)<\infty\}.$$
For $\lambda\in\Lambda^+_{m|n}$ we fix the grading in such a way that
$L(\lambda)\in\cF(\fg)$. Then 
$$\Irr(\cF(\fg))=\{L(\lambda)|\ \lambda\in \Lambda^+_{m|n}\}.$$
We denote by $L_{\fg_x}(\nu)$ the  simple $\fg$  of the highest weight $\nu$.

We denote by $\chi_{\lambda}$ the central character
of $L(\lambda)$.

%\subsection{The algebras $\fg_i$}\label{fgi}
%For $\fgl(m|m),\osp(2m|2m),\osp(2m+1|2m)$
% we denote the base $\Sigma$ by $\Sigma_m$; for $\osp(2m+2|2m)$ with $
% we denote the base $\Sigma$ by $\Sigma_{m+1/2}$:
%$$\begin{array}{ll}
%\osp(2m|2m) &
%\Sigma_m=\{\delta_1-\vareps_1,\vareps_1-\delta_2,\ldots,\vareps_{m-1}-\delta_m,
%\delta_m\pm\vareps_m\} \ \ \ \ \ \ \ \ \ \ \rho=0;\\
%\osp(2m+2|2m) &
%\Sigma_{m+\frac{1}{2}}=\{\vareps_1-\delta_1,\delta_1-\vareps_2,\ldots,\vareps_{m}-\delta_m,
%\delta_m\pm\vareps_{m+1}\} \ \ \ \ \ \ \  \rho=0;\\
%\osp(2m+1|2m) &
%\Sigma_m=\{\vareps_1-\delta_1,\delta_1-\vareps_2,
%\ldots,\vareps_{m}-\delta_m,\delta_m\} \ \ \
%\rho=\frac{1}{2}\displaystyle\sum_{i=1}^m(\delta_i-\vareps_i).\end{array}$$
%
%If $\Sigma$ contains a base $\Sigma'_i$ isomorphic to $\Sigma_i$, we denote the corresponding Kac-Moody subalgebra of $\fg$ by $\fg_i$: for $\osp$-case $\Sigma'_i$ consists of
%last $2i$ roots of $\Sigma$; for $\fgl(m|n)$-case $\Sigma'_i$ consists
%of $2i-1$ roots, where the isotropic root is in the middle position
%(for example, $\Sigma_2=\{\vareps_{m-1}-\vareps_m,\vareps_m-\delta_1,\delta_1-\delta_2\}$). We consider $\fg_i$ is defined for the following values of $i$:
%$$\begin{array}{ll}
%i=1,\ldots,\min(m,n)\ \text{ for }\fg\not=\osp(2m|2n),\\
%i=1,3/2,2,\ldots, m$
%for $\osp(2m|2n)$ with $m\leq n$
%and $i=1,3/2,2,\ldots, n+1/2$
%for $\osp(2m|2n)$ with $m>n$.
%
\section{Weights, roots and diagrams} \label{sec:weights}
\subsection{Weight diagrams}\label{weightdia}
Many properties of a finite dimensional representation $L(\lambda)$ can be better understood by assigning a {\em weight diagram} to the weight $\lambda$ (see e.g. \cite{GS}).  Note that the conventions how to draw these weight diagrams differ. We follow essentially \cite{GS} and list some differences below. These weight diagrams were first introduced in the $\mathfrak{gl}(m|n)$-case in \cite{BS4} in their work on Khovanov algebras of type $A$. The weight diagrams introduced in \cite{ES} for representations of the orthosymplectic supergroup differ considerably from ours (see \cite{ES}, Section 6).

Take $\lambda\in\Lambda^+_{m|n}$ and write
$$\lambda+\rho=:
\displaystyle\sum_{i=1}^m a_i\vareps_i+\sum_{j=1}^n b_j\delta_j.$$

For $\osp(2m|2n)$-case all coefficients $a_i,b_j$ are integers which are non-negative except for $a_m$; for $\osp(2m+1|2n)$-case the numbers $a_i+1/2,b_j-1/2$ are
non-negative integers. Moreover,
$a_i\geq a_i', b_j\geq b'_j$ for $i<i', j<j'$ and the equalities imply
$a_i=0, b_j=0$ in $\osp(2m|2n)$-case, $a_i=-\frac{1}{2}, b_i=\frac{1}{2}$
in $\osp(2m+1|2n)$-case. The weight $\lambda$ is typical if $a_i\not=b_j$ for all $i,j$.

\subsubsection{}
We assign to $\lambda$ the (weight) diagram using the following procedure: we label the numberline $\mathbb{N}$ as follows:

for $\osp(2m|2n)$ we put  $>$  (resp., $<$) at the position with the coordinate $t$ if $|a_i|=t$ (resp., $|b_i|=t$)
for some $i$;

for $\osp(2m+1|2n)$ we put  $>$  (resp., $<$) at the position with the coordinate $t-1/2$ if $|a_i|=t$ (resp., $|b_i|=t$).

If $>,<$ occupy the same position we put the sign $\times$ ($\times^s$ stands for
$s$ signs $<$ and $s$ signs $>$; $\frac{>}{\times^s}$ stands for
$s$ signs $<$ and $s+1$ signs $>$). We put  an empty sign $\circ$ at the  non-occupied
positions with the coordinates in $\mathbb{N}$.
We call the resulting diagram an {\em unsigned weight diagram}.

For $\osp(2m|2n)$ with $m\geq 1$   we add the sign $+$ (resp., $-$)
if $a_m>0$ (resp., $a_m<0$).

For $\osp(2m+1|2n)$
we put the sign $+$ (resp., $-$) before the diagram if the zero position is occupied
by $\times^p$ for $p>0$ and $(\lambda+\rho|\vareps_i)=\frac{1}{2}$ for some $i$
(resp., $(\lambda+\rho|\vareps_i)\not=\frac{1}{2}$ for each $i$).

We call the resulting diagram a {\em  weight diagram} of $\lambda$.

Notice that for $\osp(2m|2n)$-case the action of automorphism $\sigma:\lambda\mapsto\lambda^{\sigma}$
corresponds to the change of signs of the diagrams; we denote 
this operation (the change of sign) also by $\sigma$.

\subsubsection{}
Note that for $\osp(2m+1|2n)$-case our weight diagram
is obtained from  the diagram used in \cite{GS} by the shift
by $-1/2$.

\subsubsection{Examples: $\osp(2m|2n)$}
The empty diagram corresponds to $\osp(0|0)=0$; the diagram
$>$ (respectively, $-\circ>$)
corresponds to the weight $0$ (respectively, $-\vareps_1$) for $\osp(2|0)=\mathbb{C}$.
One has $L(\emptyset)=\mathbb{C}$ and $L(>)$ is the trivial $\osp(2|0)$-module.

For $s>0$ the diagram $\times^s$  is assigned to the weight $0$  for
$\osp(2s|2s)$; the diagrams
 $\overset{\times^s}{>}>,\times^s>>$ are assigned to the 
$\osp(2s+4|2s)$-weights $0$ and $\vareps_1+\vareps_2$ respectively; the diagrams
$\pm\circ \times \circ \times$ are assigned to the $\osp(4|4)$-weights
$3(\delta_1+\vareps_1)+(\delta_2\pm\vareps_2)$.

\subsubsection{Examples: $\osp(2m+1|2n)$}
The empty diagram corresponds to $\osp(1|0)=0$; the diagram
$>$ (resp., $<$) corresponds to the weight $0$ for
$\osp(3|0)=\mathfrak{o}_3$  (resp., for
$\osp(1|2)$).

For $s>0$ the diagram $- \times^s$ is assigned to the weight $0$  for
$\osp(2s+1|2s)$; the diagram $+ \times^s$ is assigned to the  weight $\vareps_1$
for $\osp(2s+1|2s)$.

The diagram $-\times^n>>$
is assigned to the $\osp(2n+5|2n)$-weight $0$
and the diagram $+\times^n<<$
is assigned to $\osp(2n+1|2n+4)$-weight $\vareps_1$.

\subsubsection{}
The above procedure gives a one-to-one correspondence
between $\Lambda^+_{m|n}$ and the diagrams containing $k$ symbols $\times$,
$m-k$ symbols $>$ and $n-k$ symbols $<$ (where $k\leq \min(m,n)$)
with the following additional properties:

\begin{enumerate}
\item
the coordinates of the occupied position lie in $\mathbb{N}$ and
each non-zero occupied
position contains exactly one of the signs $\{>,<,\times\}$;

\item for $\osp(2m|2n)$ with $m>0$ 
the zero position contains any number of $\times$, no $<$ and at most one $>$;
the diagram has a sign if and only if  the zero position is empty;

\item for $\osp(2m+1|2n)$ the zero position contains any number of $\times$ and at most one of the symbols $>, <$;
the diagram has a sign if and only if the zero position is occupied by $\times^i$
for $i\geq 1$.
\end{enumerate}

\subsubsection{}
The atypicality of $\lambda$ is equal to the number of symbols $\times$ in the
diagram of $\lambda$.

\subsubsection{Notation}
We sometimes identify a dominant weight and its weight diagram;
for instance, $f\in\Lambda^+_{m|n}$ means that $f$ is a weight diagram
assigned to a weight in $\Lambda^+_{m|n}$. 

We always consider $t\in\{0,1,2\}$ to be fixed.
For a weight diagram $f$
we sometimes use the notation $L(f)$ for the corresponding highest weight module.
For instance, $L(\emptyset)=\mathbb{C}$ and $L(>)$ is the trivial $\osp(2|0)$-module;
$L(\times^s)$ (resp., $L(-\times^s)$) 
stands for the trivial $\osp(2s|2s)$ (resp., $\osp(2s+|2s)$-module)
and $L(+\times^s)$ stands for the standard $\osp(2s+1|2s)$-module.

\subsubsection{Remark: $OSp$-modules}\label{OSP2m2n}
Take $\fg=\osp(2m|2n)$.
By~\cite[Proposition 4.11]{ES} the simple $OSp(2m|2n)$-modules are
either of the form $L(\lambda)$ if $\lambda\in\Lambda^+_{m|n}$ is $\sigma$-invariant or $L(\lambda)\oplus L(\lambda^{\sigma})$.
Thus the simple $OSp(2m|2n)$-modules are
in one-to-one correspondence with the unsigned
$\osp(2m|2n)$-diagrams. For $\osp(2m+1|2n)$ and any $\lambda \in \Lambda_{m|n}^+$ there are two irreducible $OSp(2m+1|2n)$-modules $L(\lambda,+)$ and $L(\lambda,-)$ which restrict to $L(\lambda)$.

\subsubsection{}
We denote by $f_-f_+$
the diagram obtained by "gluing"
the diagrams $f_-$ and $f_+$ (where $f_+$ does not have sign);
for instance,
$$\begin{array}{lll}
\overset{\times}{\times}\circ\times=f_-f_+ & \text{ where } &
f_-=\overset{\times}{\times}\circ,\ \ \  f_+=\times\\
+\times^2 \circ\times=f_-f_+ & \text{ where } &
f_-=+\times^2,\ \ \  f_+=\circ\times
\end{array}
$$

\subsection{Core diagrams and central characters}
We say that a $\fg$-central character is {\em dominant} if
$\cF(\fg)$ contains modules with this central character.
By~\cite{GS}, the blocks in $\cF(\fg)$ are parametrized by the dominant central characters
and the dominant central characters can be described
in terms of typical dominant weights, see below.

\subsubsection{}\label{corewt}
We call the symbols $>,<$ the {\em core symbols}.
A {\em core diagram} is a
 weight diagram which does not contain symbols $\times$
 and does not have $-$ sign.

For a weight diagram $f$ we denote by $\core(f)$ the core diagram which is
obtained from the diagram of $f$ by replacing all symbols $\times$ by $\circ$
and by adding the sign $+$ for $\osp(2m|2n)$-case
if the zero position is empty. 

For $\lambda\in\Lambda^+_{m|n}$ we denote by
$\core(\lambda)$  the weight corresponding to $\core(f)$,
where $f$ is the diagram of $\lambda$. 

For instance,
for $\osp(4|2)$ one has 
$$\core(\times >)=\core(\pm \circ >\times)=+\circ >,\ \ 
\core(\vareps_1)=\core (\pm \vareps_2+\vareps_1+\delta_1)=\vareps_1.$$
Note that for $\osp(2m|2n)$ one has $\core(\lambda)=\core(\lambda^{\sigma})$.

\subsubsection{}\label{corechi}
For $\osp(2m+1|2n)$-case  the dominant central characters are parametrized by the core diagrams, i.e.
for $\lambda,\nu\in\Lambda^+_{m|n}$
$$\chi_{\lambda}=\chi_{\nu} \ \Longrightarrow\ \ \ \core(\lambda)=\core(\nu);$$
for $\osp(2m|2n), \osp(2m+2|2n)$ one has
$$\chi_{\lambda}\in\{\chi_{\nu},\chi_{\nu^{\sigma}}\} \ \Longrightarrow\ \ \ \core(\lambda)=\core(\nu).$$
(Note that the atypical central characters of $\osp(2m|2n)$ are $\sigma$-invariant:
$\chi_{\nu}=\chi_{\nu^{\sigma}}$ if $\nu$ is atypical.)

\subsection{Categories  $\cF^g(\fg)$}\label{catFgg}
For each core diagram $g$ we
denote by $\cF^g(\fg)$ the Serre subcategory  of $\cF(\fg)$ \footnote{  by Serre subcategory generated by a set of simple modules
we mean the full subcategory consisting of the modules of finite length whose all simple
subquotients lie in a given set.} generated by $L(\lambda)$ for $\lambda$ having diagrams $f$ with $\core(f)=g$.

By above,  the categories
$\cF^g(\osp(2m+1|2n))$ are the blocks in $\cF(\osp(2m+1|2n))$.
For $\fg:=\osp(2m|2n)$ all atypical blocks and $\sigma$-invariant
typical blocks 
are of the form $\cF^g(\fg)$; for a typical block $\cB$ with
$\cB\not=\cB^{\sigma}$
we have
$\cF^g(\fg)=\cB\oplus\cB^{\sigma}$ for a suitable diagram $g$.

Similarly we use the notation $\tilde{\cF}^g(\fg)$ for $\cF^g(\fg) \oplus \Pi \cF^g(\fg)$.

\subsubsection{}
Let $g$ be a core diagram with $m'$ symbols $>$ and $n'$ symbols $<$.

The category $\cF^g(\osp(2m+1|2n))$ is non-zero if and only if
$m-m'=n-n'\geq 0$ and that  $\cF^g(\osp(2m|2n))$ is non-zero if and only if
$m-m'=n-n'\geq 0$ and, in addition, $g$ does not have $<$ at the zero position
for $m>0$.

The  modules in $\cF^g(\fg)$ have atypicality $m-m'=n-n'$.

\subsection{Cases $t=0,1,2$}\label{reft}
Recall that $\osp(M|N)$ consists of two series: $B$ (for odd $M$) and $D$.
We will distinguish the following cases ($t=0,1,2$):

for the $B$-series (and any core diagram $g$) we put $t:=1$;

for the $D$-series  and a core diagram $g$
without $>$ at the zero position we put $t:=0$;

for the $D$-series and  a core diagram $g$
with $>$ at the zero position we put $t:=2$.

\subsubsection{}
We say that a block has type $t$ ($t=0,1,2$) if the core
diagram of the corresponding central character has type $t$.
We say that $\lambda\in\Lambda^+_{m|n}$ has type $t$ if
the diagram of $\core(\lambda)$ has type $t$. Then
 $L(\lambda)$ lies 
in a block of type $t$ if and only if $\lambda\in\Lambda^+_{m|n}$ has type $t$.

\subsubsection{}\label{ell}
For $t=0,1,2$ we take
 $\fg=\osp(2m+t|2n)$;  the  weights lattice of $\fg$ is
$\Lambda_{m+\ell|n}$, where 
$$ \ell:=\left\{\begin{array}{ll}
0 & \text{ for }t=0,1\\
1&  \text{ for }t=2.
\end{array}\right.$$

We denote by $\Lambda^{(t)}_{m+\ell|n}$ the set 
of dominant weights of type $t$ in $\Lambda^{+}_{m+\ell|n}$:

$\Lambda^{(1)}_{m+\ell|n}=\Lambda^+_{m|n}$ for the B-series;

$\Lambda^{(0)}_{m+\ell|n}$ is the set of  $\osp(2m|2n)$-dominant weights 
in $\Lambda^+_{m|n}$ with the diagrams
without symbol $>$ at the zero position;

$\Lambda^{(2)}_{m+\ell|n}$ is the set of $\osp(2m+2|2n)$-dominant weights 
in $\Lambda^+_{m+1|n}$ with  the diagrams
having $>$ at the zero position.

We call $\lambda\in\Lambda^{(t)}_{m+\ell|n}$ {\em core-free} if $\core(\lambda)=\emptyset$ for $t=0,1$ and $\core(\lambda)=>$ for $t=2$.

\subsubsection{}\label{false}
Observe that for $\lambda\in\Lambda_{m+\ell|n}^{(t)}$  we have
$$\lambda+\rho=:\sum_{i=1}^{m+\ell} a_i\vareps_i+\sum_{j=1}^n b_j\delta_j=
\sum_{i=1}^{m} a_i\vareps_i+\sum_{j=1}^n b_j\delta_j,$$
since for the case $\ell\not=0$ one has $\ell=1, a_{m+1}=0$.

For $t=0,1$ the weight  $\lambda\in\Lambda_{m+\ell|n}^{(t)}$ has atypicality $k$
if and only if $\core (\lambda)\in \Lambda^+_{m-k|n-k}$.
For $t=2$ the weight  $\lambda\in\Lambda_{m+\ell|n}^{(t)}$ has atypicality $k$
if and only if $\core (\lambda)\in \Lambda^+_{m+1-k|n-k}$; note that in this case
the core diagram of $\lambda$ has $>$ at the zero position, so
$\core(\lambda)$ lies in $\Lambda^+_{m-k|n-k}\subset \Lambda^+_{m+1-k|n-k}$.
Hence in all cases
$$\core (\lambda)\ \text{ is a typical weight in }\Lambda^+_{m-k|n-k}.$$

\subsubsection{}
For $t=0,1,2$
we denote by $\Theta_{k}^{(t)}$ the set of dominant central characters of atypicality $k$ corresponding to the $t$-case. By above,  
 the  map
$$\chi\mapsto \core (\chi)$$
gives a correspondence 
between $\Theta^{(t)}_{k}$ and the set of typical weights $\Lambda^+_{m-k|n-k}$. For  $t=1$ this is a one-to-one correspondence.
For $t=0,2$ the image  consists of the typical weights $\eta$
satisfying $(\eta|\vareps_{m-k})\not=0$; the map is injective except for the case
$t=0$ and $k=0$.

\subsection{Algebra $\fg_r$}\label{Sigmar}
For $t=0,1,2$ and  $r>0$ we set
$$\fg_r:=\osp(2r+t|2r).$$

Let $\Sigma$ be the base of simple roots for $\fg=\osp(2m+t|2n)$.
For $0<r\leq \min(m;n)$ we denote by $\Sigma_r$ a subset of $\Sigma$
which is a base of an algebra isomorphic to $\fg_r$ (such $\Sigma_r$ is unique):
$$\Sigma_r:=\left\{\begin{array}{ll}
\vareps_{m-r+1}-\delta_{n-r+1},\delta_{n-r+1}-\vareps_{m-r+2},\ldots,\vareps_m-\delta_n,\delta_n & \text{ for }t=1\\
\delta_{n-r+1}-\vareps_{m-r+1},\ldots,\vareps_{m-1}-\delta_n,\delta_n-\vareps_m,
\delta_n+\vareps_m & \text{ for }t=0\\
\vareps_{m-r+1}-\delta_{n-r+1},\delta_{n-r+1}-\vareps_{m-r+2},\ldots,\vareps_{m}-\delta_n,\delta_n\pm\vareps_{m+1} & \text{ for }t=2.
\end{array}\right.$$

For $r=0$ we set $\Sigma_r=\emptyset$ and $\rho_r=0$. 
Recall that $\osp(0|0)=\osp(1|0)=0$;
 for $t=2$ we identify  $\osp(2|0)$ with $\mathbb{C}\vareps_m^*\subset\fh$
(where $\vareps_m^*\in \fh$ is such that $\mu(h)=(\mu|\vareps_m)$
for each $\mu\in\fh^*$).

We identify $\fg_r$ with the corresponding  subalgebra of $\osp(2m+t|2n)$ (then $\Sigma_r$ is the base of of $\fg_r$). We denote by $\rho_r$ the  Weyl vector of $\fg_r$;
one has $\rho_r=\rho|_{\fg_r\cap \fh}\ $ and
$$\rho_r=0\ \text{ for }t=0,2;\ \ \ \ 
\rho_r=\frac{1}{2}\sum_{i=0}^{r-1}(\delta_{n-i}-\vareps_{m-i}) \ \ \text{ for } t=1.$$

We denote by $S_r$ the following set: $S_0=\emptyset$ and
$$S_r:=\left\{\begin{array}{ll}
\{\delta_{n-i}-\vareps_{m-i}\}_{i=0}^{r-1} & \text{ for $t=0$}\\
\{\vareps_{m-i}-\delta_{n-i}\}_{i=0}^{r-1} & \text{ for $t=1,2$.}\\
\end{array}\right.$$
Notice that $S_r$ consists of $r$ isotropic mutually orthogonal roots and $S_r\subset\Sigma_r$.

\subsection{Tail}\label{tailsub}
Take  $\lambda\in\Lambda^{(t)}_{m+\ell|n}$
and let $f$ be the diagram assigned to $\lambda$.

Let $s\leq \min(m,n)$ be the maximal number satisfying
$(\lambda|\Sigma_s)=0$, where $\Sigma_s$ as above (notice that $\Sigma_s$
depends on $t$).
We call $s$ the {\em length of the tail} of $\lambda$ and write
$$|\tail(\lambda)|=|\tail(f)|:=s.$$

For $\osp(2m|2n)$-case $s$ is equal to the  number of symbols $\times$
in the zero position.
For $\osp(2m+1|2n)$  the zero position contains $s$ (resp., $s+1$) symbols $\times$ if the diagram has sign  $-$ (resp., $+$); for instance
$$|\tail((-)\times^m)|=|\tail(0)|=m,\ \ \ \ \
|\tail((+)\times^m)|=|\tail(\vareps_1)|=m-1.$$

\subsection{Howl}\label{howl}
The block of the trivial module is called the {\em principal block} (there are two principal blocks which differ by $\Pi$).

Each block of atypicality $k$ is equivalent to the principal block of
$\osp(2k+t|2k)$.
The equivalences are described in~\cite{GS}; we give some details below.
For a dominant weight $\lambda$ we denote by $\howl(\lambda)$  the corresponding weight in the principal block (roughly speaking,
the passage from $\lambda$ to $\howl(\lambda)$ essentially amounts to removing the core symbols $<,>$ from the weight diagram, see the details  below);
note that $\howl(\lambda)$ has  core-free diagram.

\subsubsection{}
Let  $\lambda\in\Lambda^{(t)}_{m+\ell|n}$ be a weight of atypicality $k$ 
and let $f$ be the corresponding diagram.
For each $i=1,\ldots,k$
 let $s_i(f)$ be  the number of the positions to the left of $i$th symbol $\times$, which do not contain core symbols.

 \subsubsection{Case $t=1$}\label{2r12r}
In this case $\fg=\osp(2m+1|2n)$ and the
diagram  $\howl(f)\in\Lambda^{(1)}_{k|k}$ is  the diagram
  without core symbols, where
 $i$th symbol $\times$ occupies   $s_i(f)$th position and the sign of $\howl(f)$
 is such that the  tail lengths
of $\howl(f)$ and $f$ are the same. For instance,
$$\begin{array}{l}
\howl(<\times<\times)=+\times\times;\ \ \howl(\pm \times^2 \circ>\times)=\pm\times^2\circ \times;\ \ \\
\howl(\overset{\times}{>}\circ \circ >\times)=-\times\circ \times;\ \
\howl(\overset{\times}{>}\times \circ >\times)=+\times^2\circ \times
\end{array}$$
and $\howl(f)=\emptyset$ if and only if $k=0$.

 \subsubsection{Case $t=0$}\label{2r2r}
In this case $\fg=\osp(2m|2n)$ with no $>$ at the zero position.
The diagram  $\howl(f)\in\Lambda^{(0)}_{k|k}$ is the diagram
  without core symbols, where the 
 $i$-th symbol $\times$ occupies the $s_i(f)$-th position and the sign  of $\howl(f)$
 coincides with the sign of $f$ if $\howl(f)$ requires the sign.  For instance,
 $$\howl(\pm\circ>>\times)=\pm \circ \times;\ \ \ \howl(\times^2>\times)=\times^2\times;\ \ \
 \howl(\pm \circ >\circ <)=\emptyset.$$
 Notice that
  $|\tail(f)|=|\tail(\howl(f))|$,  so the only case when the diagrams
  $f$ and $\howl(f)$ do not have the same sign is when
  $\howl(f)=\emptyset$.

 \subsubsection{Case $t=2$}\label{2r2r2}
In this case $\fg=\osp(2m+2|2n)$ and the zero position of $f$
 is occupied by $\overset{\ \ \times^{p}}{>}$ ($p\geq 0$).  The diagram $\howl(f)\in\Lambda^{(2)}_{k+1|k}$ has $\overset{\ \ \times^{p}}{>}$
at the zero position; for $i=p+1,\ldots, k$ the $i$th symbol $\times$
in $\howl(f)$ occupies the position $s_i+1$. For instance,
$$\begin{array}{l}
\howl(\overset{\times^2}{>}\circ \circ >\times)=\overset{\times^2}{>}\circ \circ \times\\
\howl(>\times<\times)=>\times\times;\ \
\howl(><)=>
\end{array}
$$

 \subsubsection{}
For $\lambda\in\Lambda^{(t)}_{m+\ell|n}$ with a diagram $f$ 
let
 $\howl(\lambda)\in\Lambda^+$ be  the weight corresponding to $\howl(f)$.
Notice that for
$\lambda\in\Lambda^{(t)}_{m+\ell|n}$  one has
$\howl(f)\in\Lambda^{(t)}_{k+\ell|k}$, where $k$ is atypicality of $\lambda$.

If $\ell\not=0$, then $t=2$ and $\howl(f)$ has $>$ at the zero position,
so $\howl(f)$ lies in $\Lambda^+_{k|k}$. Hence in all cases
$\howl(f)\in\Lambda^+_{k|k}$.

 \subsubsection{}
Note that
$\howl$ preserves  the tail length:
$$ |\tail(f)|=|\tail(\howl(f))|.$$

\subsubsection{Connection between the cases $t=1$ and $t=2$}\label{t12}
Below we describe a remarkable bijection between the core-free diagrams in $\Lambda_{(n|n)}^{(1)}$ and in $\Lambda_{(n|n)}^{(2)}$.

To a diagram
$\overset{\ \ \times^{p}}{>}\circ f$ we assign the diagram $-\times^p f$
if $p>0$ and the diagram $\circ f$ if $p=0$;
to a diagram
$\overset{\ \ \times^{p}}{>}\times f$ we assign the diagram $+\times^{p+1} f$.

This assignment gives a bijection $\tau$ between the core-free
diagrams for $\osp(2n+2|2n)$ and  for $\osp(2n+1|2n)$. Notice that
$$|\tail(\tau(f))|=|\tail(f)|.$$

\subsubsection{Map $\pari$}\label{mapdex}
We introduce a map $\pari:\Lambda^{(t)}_{m+\ell|n}\to \{\pm 1\}$ by
$$\pari(\lambda):=\left\{\begin{array}{ll}
(-1)^{p(\howl(\lambda))}& \text{ for } t=0,1\\
(-1)^{p(\tau(\howl(\lambda)))} & \text{ for } t=2
\end{array}
\right.
$$
(where $p$ is the parity) and the map $\Irr\bigl(\tilde{\cF}(\fg)\bigr)\to \{\pm 1\}$ by
$$\pari(L(\lambda)):=\pari(\lambda),\ \ \ \pari(\Pi(L(\lambda))=-\pari(\lambda).$$

\subsection{Stable diagrams}\label{sec:stable}
For $\osp(2m+1|2n)$
a diagram  is called {\em stable} if all symbols $\times$ precede
all core symbols ($<,>$).
For $\osp(2m|2n)$ a diagram  is called {\em stable} if all symbols $\times$ preceed
all core symbols, except, possibly, the symbol $>$ at the zero position.

A weight $\lambda$ is called {\em stable}
if the corresponding diagram is stable.

\subsubsection{}\label{stablecore}
Take $\lambda\in\Lambda^{(t)}_{m+\ell|n}$  
and write
$$\lambda+\rho=:\sum_{i=1}^{m} a_i\vareps_i+\sum_{j=1}^n b_j\delta_j.$$
If $\lambda$ has  atypicality $k>0$, then
$$\lambda\ \text{ is stable }\ \Longleftrightarrow\
\core (\lambda)=\sum_{i=1}^{m-k} a_i\vareps_i+\sum_{j=1}^{n-k} b_j\delta_j.$$
(The same holds for $k=0$ with $t=1,2$.)

\subsubsection{}\label{howlst}
If $f$ is a stable diagram
for $\osp(2m+1|2n)$ (resp., for $\osp(2m|2n)$) the diagram
$\howl(f)$ is obtained from $f$ by replacing
all core symbols (resp., all core symbols in the non-zero position) by the empty symbols.

In other words, for $\lambda$ as above,
$\lambda$ is stable of atypicality $k$ if and only if
$$\howl(\lambda)+\rho_k=\sum_{i=1}^{k} a_{m-k+i}\vareps_i+\sum_{j=1}^k b_{n-k+j}\delta_j,$$
where $\rho_k$ is the Weyl vector of $\fg_k$.

\section{Stabilization}\label{sect2}
We call a module {\em  stable} if all its simple subquotients are of the form $L(\lambda)$ for stable weights
$\lambda$. The aim of this section is to show that any module 
  in $\cF(\fg)$ can be moved with a translation functor to a stable module.

\subsection{Translation functors}\label{Transfunc}
Let $V_{st}$ be the natural representation.
For a core diagrams $g,g'$ we denote
by $T^{g'}_g$ the translation functors
$$T^{g'}_g:\ \cF^{g}(\fg)\to\cF^{g'}(\fg),\ \ \ \ \
T^{g'}_g:\ \tilde{\cF}^{g}(\fg)\to\tilde{\cF}^{g'}(\fg)$$
which map $N$ to the projection of $N\otimes V_{st}$ to the subcategory
$\cF(\fg)^{g'}$ (resp., $\tilde{\cF}^{g'}(\fg)$).

We write $T^{g'}_g(f )= f'$ if $T^{g'}_g(L(\lambda))=L(\lambda')$ and
$f,f'$ are the weight diagrams assigned to $\lambda,\lambda'$ respectively;
similarly, we write
$T^{g'}_g(f)=f'_1\oplus f'_2$ if $T^{g'}_g(L(\lambda))=L(\lambda'_1)\oplus L(\lambda'_2)$.

\subsection{Some useful translation functors}\label{stabilization}
Let $Trans_a$ be the set of the  translation functors $T^{g'}_g$, where
$g$ is a core diagram with an occupied position $a$ and an empty position $a+1$
and  $g'$ is the core diagram obtained from $g$ by interchanging the symbols
in the positions $a,a+1$. For example, $Trans_1$ contains
$T^{\circ\circ >}_{\circ >\circ,}$, $T^{>\circ <<}_{><\circ<}$ and so on.

Each translation functor $T^{g'}_{g}\in Trans_a$ is an equivalence of categories except for the case $\osp(2m|2n)$ with $a=0$. In the $\osp(2m+1|2n)$-case (resp., in the $\osp(2m|2n)$-case)
let $Trans$ be the functors which can be written as
compositions  of functors from $Trans_a$ for all $a$ (resp., $a\not=0$).
All functors in $Trans$ are equivalences of categories.

\subsubsection{Case $a\not=0$}
In this case each $T\in Trans_a$ is an equivalence of categories acting
on simple modules by
interchanging  the symbols in positions $a,a+1$ in the corresponding diagrams
(and preserving the sign); for example,
$$T^{\circ\circ >}_{\circ >\circ}(\times > *)=\times * >, \ \ \
T^{\circ\circ >}_{\circ >\circ,}(+\circ > *)=+\circ  * >$$
where $*\in\{\times,\circ\}$.

\subsubsection{Case $\fg=\osp(2m+1|2n), a=0$}
In this case each $T\in Trans_0$ is an equivalence of categories acting
on simple modules by the following rules:
$$\begin{array}{lccl}
>\circ f\mapsto \circ>f & & & \overset{\times^i}{>}\circ f\mapsto -\times^i >f\\
 >\times f\mapsto +\times > f & & &
\overset{\times^i}{>}\times f\mapsto +\times^{i+1} >f
\end{array}$$
for each $i>0$ and similar formulae, where $>$ is changed by $<$.

\subsection{}
\begin{cor}{cor1}
(i) For $T\in Trans$ one has
$\howl(T(f))=\howl(f)$.

(ii) For any $\lambda\in\Lambda^+_{m|n},\lambda'\in\Lambda^+_{m-i|n-i}$   satisfying $$\core(\lambda)=\core(\lambda')$$
there exists $T\in Trans$ such that
$T(\lambda), T(\lambda')$ are stable.

(iii) For any module  $N\in\cF^{g}(\fg)$ there exists $T\in Trans$ such that
$T(N)$ is stable.
\end{cor}

\section{$\DS$-functor} \label{sec:DS}
In this section we recall the construction of the $\DS$-functor and 
describe the algebra $\fg_x$. We prove Corollaries \ref{corst},\ref{propDsstable}
which will be used later. We distinguish the cases $t=0,1,2$ and take
$\fg=\osp(2m+t|2n)$.

\subsection{}
\label{sectionDS}
The $\DS$-functor was introduced in \cite{DS}. We recall definitions
and some results below.
For a $\fg$-module $M$ and $g\in\fg$ we set
$$M^g:=\Ker_M g.$$

\subsection{Construction}\label{DSconst}
We set $\fg_x:=\fg^{x}/[x,\fg]$; note that $\fg^{x}$ and $\fg_x$ are Lie superalgebras.
For a $\fg$-module $M$ we set
$$\DS_x(M)=M^x/xM.$$
Observe that $M^x, xM$ are $\fg^{x}$-invariant and $[x,\fg] M^x\subset xM$,
so $\DS_x(M)$ is a $\fg^{\ad x}$-module and $\fg_x$-module.
Thus $\DS_x: M\to \DS_x(M)$ is a functor from the category of $\fg$-modules to
the category of $\fg_x$-modules.

There are canonical isomorphisms $\DS_x(\Pi(N))=\Pi(\DS_x(N))$ and
$$\DS_x(M)\otimes\DS_x(N)=\DS_x(M\otimes N).$$
If $N$ is a finite-dimensional $\fg$-module, then $\DS_x(N^*)\cong (\DS_x(N))^*$.

\subsubsection{Algebraic representations} The $DS$ functor restricts to a functor \[ DS_x: \tilde{\cF}(\fg) \to \tilde{\cF}(\fg_x).\] It does not however preserve the subcategory $\cF(\fg)$. As already noted in \cite{CH} it induces a symmetric monoidal functor between the algebraic representations of $OSp(m|2n)$ and $OSp(m-2r|2n-2r)$.

\subsubsection{DS and core diagrams}\label{DSgT}
By \cite{DS}, Sect. 7 (see also Thm. 2.1 in~\cite{Skw}),
 the $\DS$-functors preserve the core diagrams, i.e. for a core diagram $g$ one has
\begin{equation}\label{DSg}
\DS_x(\cF^g(\fg))\subset \tilde{\cF}^g(\fg_x),
\end{equation}
where $\fg_x:=\DS_x(\fg)$. Warning: $\DS_x(\cF^g(\fg))$ is in general not in $\cF^g(\fg_x)$, since $\DS_x$ does not preserve $\cF(\fg)$.

\subsubsection{DS and Translation functors} 
Since $DS$ is a symmetric monoidal functor  
$$\DS_x(N\otimes V_{st})=\DS_x(N)\otimes \DS_x(V_{st}).$$
Since $\DS_x(V_{st})$ is the natural representation of $\fg_x$,
 the translation functors
"commute with the DS-functors", i.e. the following
 diagram is commutative
$$\begin{array}{ccc}
\cF^{g_1}(\fg) & \overset{T^{g_2}_{g_1}}{\longrightarrow} & \cF^{g_2}(\fg)\\
\ & &\\
\DS_x\downarrow & & \DS_x\downarrow\\
\ & &\\
\tilde{\cF}^{g_1}(\fg_x) &  \overset{T^{g_2}_{g_1}}{\longrightarrow} & \tilde{\cF}^{g_2}(\fg_{x})
\end{array}$$

\subsection{}
\begin{cor}{corstabi}
For any $\lambda\in\Lambda^+_{m|n},\nu\in\Lambda^+_{m-s|n-s}$ with
$$\core(\lambda)=\core(\nu)$$
there exist stable weights $\lambda_{st}\in\Lambda^+_{m|n},\nu_{st}\in\Lambda^+_{m-s|n-s}$ such that
$$\begin{array}{c}
\core(\lambda_{st})=\core(\nu_{st});\\
\howl(\lambda)=\howl(\lambda_{st}),\ \ \ \ \ \howl(\nu)=\howl(\nu_{st})
\end{array}$$
and $\ \ [\DS_x(L(\lambda):L(\nu)]=[\DS_x(L(\lambda_{st}):L(\nu_{st})]$ for each $x$ of rank $s$.
\end{cor}
\begin{proof}
This follows from~\Cor{cor1} and the fact that
the translation functors commute with $\DS_x$ (see \ref{DSgT}).
\end{proof}

\subsection{$\DS$ and automorphisms}\label{inner}
Let $\phi:\fg'\to \fg$ be a homomorphism of Lie superalgebras;
for  each $\fg$-module $N$
denote by $N^{\phi}$ the $\fg'$-module
(the vector space $N$ with  the action $g'.v:=\phi(g')v$).

Each  $\phi\in\Aut(\fg)$ induces
an isomorphism $\ol{\phi}: \fg_x\iso \fg_{\phi(x)}$ and
$$\DS_{\phi(x)}(N^{\phi})=(\DS_x(N))^{\ol{\phi}}.$$

Let $a\in\fg_0$ be an ad-nilpotent element and $\phi:=e^{\ad a}$
be the corresponding inner automorphism of $\fg$.
If $a$ acts nilpotently on a $\fg$-module $N$,
then $e^a: N\iso N^{\phi}$. Therefore
$$\DS_{\phi(x)}(N)=(\DS_x(N))^{\ol{\phi}}.$$

Let $G_0=O_{2m+\ell}\times Sp_{2n}$ be the adjoint group of $\fg_{\ol{0}}$, i.e.
 the subgroup of $\Aut \fg$
generated by $e^{\ad a}$, where $a\in\fg_0$ is ad-nilpotent.
By above, if $N$ is a finite-dimensional $\fg$-module, then
\begin{equation}\label{Nfin}
\DS_{\phi(x)}(N)=(\DS_x(N))^{\ol{\phi}}
\end{equation}
for any inner automorphism $\phi\in G_0$.

\subsection{Choice of $x$}\label{G0}
We call $S\subset\Delta_{1,re}$ an {\em isotropic set} if  $S$ is a basis of an isotropic subspace in $\fh^*$. Write $g\in\fg_{\ol{1}}$ as
$$g=\sum_{\alpha\in \supp(g)}  g_{\alpha},$$
 where $g_{\alpha}\in \fg_{\alpha}\setminus\{0\}$.

By \cite{DS}, Sect. 5 the $G_0$-orbits in
$$X_{iso}:=\{x\in \fg_1|\ [x,x]=0\}$$
 are enumerated by $0,1\ldots,\min(m+\ell,n)$: the $0$-th orbit is
$\{0\}$; the $s$-th orbit contains all elements $x$ where
$\supp(x)$ is an isotropic set of cardinality $s$. We say that
$x\in X_{iso}$ has {\em rank $s$} if $x$ lies in the $s$-th orbit.

The rank of $x\in X_{iso}$ is at most $\min(m+\ell,n)$. Assume that
the rank of $x$ is greater than $m$. In this case
$\fg=\osp(2n|2n),t=2$ and $x$ has rank $n$. However blocks of type $t=2$ 
in $\cF(\osp(2n|2n))$ have atypicality at most $n-1$, so
$\DS_x$ annihilates the modules in such blocks.

Thus we can (and will) always assume that the rank of $x$
is $s$, where
\begin{equation}\label{sminmn}
0<s\leq \min(m,n).\end{equation}

%\subsubsection{}
%Consider the case 
%$\fg:=\osp(2(m+\ell)|2n)$ (i.e., $t=0,2$).
%Let us show that
%\begin{equation}\label{NsigmaN}
%\DS_x(N^{\sigma})\cong (\DS_x(N))^{\sigma}.
%\end{equation}
%
%For each $s$ as in~(\ref{sminmn}) take
% $y\in  X_{iso}$ with
%$$\supp(y)=\{\vareps_{m+\ell-i}-\delta_{n+1-i}\}_{i=1}^s.$$
%Note that  $y$ has rank $s$ and $\sigma(y)=y$.
%Moreover, $\sigma$ induces the standard involution $\sigma$ on
%$\DS_y(\fg)\cong \osp(2(m+\ell-s)|2(n-s))$.
%Using \ref{inner} we obtain the formula (\ref{NsigmaN}) for $x:=y$.
%By \ref{G0}, this implies the formula (\ref{NsigmaN}) for each $x\in G_0y$
%and thus establishes (\ref{NsigmaN}).

\subsubsection{}\label{xyh}
For each  $s$ as in (\ref{sminmn}) we fix  $x_s$ of rank $s$ with
$$\supp(x_s):=S_s$$
and set  $\DS_s:=\DS_{x_s}$. Notice that  $x_s\in\fg_s$.

\subsubsection{}\label{dependence-on-x}
In general different $x$ (even of the same rank) give rise to different functors $DS_x$. For $\mathfrak{osp}(m|2n)$ any $x$ of rank $1$ induces the same functor $DS_x:\tilde{\cF}(\fg)\to \tilde{\cF}(\fg_x)$ by \cite{CH}, Lemma 7.7 (this is not true for $\mathfrak{gl}(m|n)$); this will later imply that $DS_x(L(\lambda))$ depends only on the rank of $x$. By \cite{CH} $ker(DS_1) = Proj$, the thick ideal of projective objects.

\subsection{The algebra $\fg_x$}\label{subsecgx}
Take $x:=x_s$. Set
$$\Delta^+_x:=\{\alpha\in\Delta^+|\ (\alpha|S_s)=0\}\setminus S_s$$
and denote by $\fg_x$ the algebra generated by $\fg_{\pm\alpha}$ with
$\alpha\in\Delta^+_x$.
Clearly, $\fg_x$ is a subalgebra of $\fg^{\ad x}$. By~\cite{DS},
$\DS_x(\fg)=\fg^{\ad x}/[x,\fg]$  can be identified with $\fg_x$.
One has
$$\fg_x\cong \osp(2(m-s)+t|2(n-s))$$
and $\fh_x:=\fg_x\cap \fh$
is a Cartan subalgebra of $\fg_x$. The triangular decomposition
given by
 $$\Delta^+(\fg_x):=\Delta^+_x$$
is of the same form as the triangular decomposition fixed in \ref{tri}. We denote the corresponding base by $\Sigma^x$ and the Weyl vector by $\rho_x$.

\subsubsection{Examples}
Take $s:=2$.

For  $\osp(11|8)$ ($m=5,n=4, t=1$)  we have $\fg_x\cong\osp(7|4)$ with 
$$\begin{array}{l}
\Sigma=\{\vareps_1-\vareps_2,\vareps_2-\delta_1,\delta_1-\vareps_3,\vareps_3-\delta_2,\delta_2-\vareps_4,\vareps_4-\delta_3,\delta_3-\vareps_5,\vareps_5-\delta_4,\delta_4\}\\
\Sigma^x=\{\vareps_1-\vareps_2,\vareps_2-\delta_1,\delta_1-\vareps_3,\vareps_3-\delta_2,\delta_2\},\\
2\rho=
\delta_4+\delta_3+\delta_2+\delta_1-\vareps_5-\vareps_4-\vareps_3-
\vareps_2+\vareps_1,\\
2\rho_x=\delta_2+\delta_1-\vareps_3-
\vareps_2
+\vareps_1.
   \end{array}$$

For $\osp(12|8)$ we have  $\fg_x\cong\osp(8|4)$ and $\rho=0,\rho_x=0$. In this case
$$\Sigma=\{\vareps_1-\vareps_2,\vareps_2-\delta_1,\delta_1-\vareps_3,\vareps_3-\delta_2,\delta_2-\vareps_4,\vareps_4-\delta_3,\delta_3-\vareps_5,\vareps_5-\delta_4,\delta_4\pm\vareps_6\};$$
for   $t=0$ we have $m=6,n=4$ and
 $$\begin{array}{l}
\Sigma^x=\{\vareps_1-\vareps_2,\vareps_2-\delta_1,\delta_1-\vareps_3,\vareps_3-\delta_2,\delta_2\pm\vareps_4\};
   \end{array}$$
 and  for $t=2$ we have $m=5,n=4$ and
 $$\begin{array}{l}
\Sigma^x=\{\vareps_1-\vareps_2,\vareps_2-\delta_1,\delta_1-\vareps_3,\vareps_3-\delta_2,\delta_2\pm\vareps_6\}.
   \end{array}$$

\subsubsection{}\label{rest}
Recall that $\fh^*$ has a basis $\{\vareps_i\}_{i=1}^{m+\ell}\cup\{\delta_i\}_{i=1}^n$
(where $\ell=0$ for $t=0,1$ and $\ell=1$ for $t=2$); it is easy to see that
 $\fh_x^*$ has a basis $\{\vareps_i\}_{i=1}^{m-s}\cup\{\delta_i\}_{i=1}^{n-s}\cup\{\vareps_{m+\ell}\}$.  For
the restriction map $\gamma\mapsto \gamma|_{\fh_x}$ we have
$$\begin{array}{l}
\vareps_i\mapsto \vareps_i \text{ for }i=1,\ldots,m-s;\ \ \
\delta_i\mapsto \delta_i \text{ for }i=1,\ldots,n-s;\\
\vareps_{m-i}\mapsto 0,\ \  \delta_{n-i}\mapsto 0\text{ for }i=0,\ldots,s-1;
\end{array}$$
and, for $\ell=1$, $\vareps_{m+1}\mapsto \vareps_{m+1}$.
One has $\rho|_{\fh_x}=\rho_x$.

Let $\{\vareps'_i\}_{i=1}^{m+\ell-s}\cup\{\delta'_i\}_{i=1}^{n-s}$ be the
standard basis in $(\fh')^*$, where $\fh'$ is the Cartan subalgebra
of $\osp(2(m-s)+t|2(n-s))$.
The isomorphism $\fg_x\iso\osp(2(m-s)+t|2(n-s))$ gives $\vareps_i\to\vareps_i'$
for $i=1,\ldots, m-s$, $\delta_j\mapsto \delta_j'$ for $j=1,\ldots,n-s$
and $\vareps_{m+1}\mapsto \vareps_{m+1-s}$ if $\ell=1$.

For $t=0,2$ we denote by $\sigma_x$ the analogue of $\sigma$ for $\fg_x$;
note that $\Sigma_x$ is $\sigma_x$-invariant. 
Retain the notation of \ref{d0d1}.

\subsection{} %from august
\begin{lem}{lemsigma}
Take $\fg:=\osp(2(m+\ell)|2n)$.

(i) If $N$ is finite-dimensional, then 
$\DS_x(N^{\sigma})\cong (\DS_x(N))^{\sigma_x}$.

(ii)
Let $x$ be of rank $1$. If $L$ (resp., $L'$) is a simple finite-dimensional
$\fg$ (resp., $\fg_x$)-module, then 
$$[\DS_x(L):L']= [\DS_x(L): (L')^{\sigma_x}].$$
\end{lem}
\begin{proof}
For each $s$ with $0<s\leq \min(m,n)$ take
 $y\in  X_{iso}$ with
$$\supp(y)=\{\vareps_{m+\ell-i}-\delta_{n+1-i}\}_{i=1}^s.$$
Note that  $y$ has rank $s$ and $\sigma(y)=y$.
Moreover, $\sigma$ induces the  involution $\sigma_y$ on the algebra 
$\DS_y(\fg)\cong \osp(2(m+\ell-s)|2(n-s))$.
Using~\ref{inner} we obtain  (i) for $x:=y$. 
By~\ref{G0}, this implies (i) for each $x\in G_0y$ and thus establishes (i) 
for all $x$.

For (ii) notice that for a simple $\osp(2m|2n)$-module $L$ one has
$L^*\cong L$ if $m$ is even and $L^*\cong L^{\sigma}$ if $m$ is odd
(since $-\Id$ lies in the Weyl group $W(\osp(2m|2n))$ if and only if $m$ is even).

If  $m+\ell$ is even, then
$L\cong L^*$, so  $\DS_x(L)\cong \DS_x(L)^*$ and  $(L')^*\cong (L')^{\sigma_x}$;
this gives the required formula.
Consider the case when  $m+\ell$ is odd.  In this case 
$L^*\cong L^{\sigma}$ and $(L')^*\cong L'$. Using (i) we get
$$\begin{array}{r}
[\DS_x(L): (L')^{\sigma_x}]=[\DS_x(L^{\sigma}): L']=[\DS_x(L^*):L']=
[\DS_x(L)^*:L']\\=
[\DS_x(L):(L')^*]=[\DS_x(L):L']\end{array}$$
as required.\end{proof}

\subsection{}
\begin{lem}{lemst}
Let $\lambda\in \Lambda^{(t)}_{m+\ell|n}$ be a stable  weight of atypicality $k$ and let $x:=x_s$.
If $\nu\in \fh^*$  satisfies
$$\nu\leq \lambda;\ \ \ \nu|_{\fh_x}\text{ is dominant };\ \ \core(\lambda)=\core(\nu|_{\fh_x}),$$
then $\nu|_{\fh_x}$ is stable.
\end{lem}
\begin{proof}
Write
$$\lambda+\rho=:\sum_{i=1}^{m+\ell} a_i\vareps_i +\sum_{j=1}^n b_j\delta_j,\ \ \
\nu+\rho=:\sum_{i=1}^{m+\ell} a'_i\vareps_i +\sum_{j=1}^n b'_j\delta_j.$$
Set $\nu':=\nu|_{\fh_x}$.   Since $\rho|_{\fh_x}=\rho_x$ we obtain
$$\nu'+\rho_x=(\nu+\rho)|_{\fh_x}.$$

Since $\nu'$ is dominant and
$\core(\lambda)=\core(\nu')$ we have
$\nu'\in \Lambda^{(t)}_{m+\ell-s|n-s}$. If
 $\ell=1$ we have $t=2$ and thus $a'_{m+1}=0$. Therefore for all $t$ we have
$$\nu'+\rho_x=(\nu+\rho)|_{\fh_x}=\sum_{i=1}^{m-s} a'_i\vareps_i +\sum_{j=1}^{n-s} b'_j\delta_j.$$

If $k=0$, then $s=0$ and $\nu'=\nu$ is stable.
If $\nu'$ is typical, it is stable.
Thus we assume that $\lambda,\nu'$ are atypical. Since $\lambda$ is stable, \ref{stablecore} gives
\begin{equation}\label{corkinet}
\core(\nu')=\core(\lambda)=\sum_{i=1}^{m-k} a_i\vareps_i+\sum_{j=1}^{n-k} b_j\delta_j.
\end{equation}
Moreover, by~\ref{stablecore}, for stability of $\nu'$ it is enough to verify that $a_i=a'_i$ for $i=1,\ldots,m-k$ and
$b_j=b'_j$ for $j=1,\ldots,n-k$.
Let $p$ (resp., $q$) be minimal such that $a_p\not=a'_p$
(resp.,  $b_q\not=b'_q$). By above
it suffices to show that
\begin{equation}\label{kmpnq}
m-p,n-q<k.\end{equation}
Using~(\ref{corkinet}) (and the atypicality of $\nu'$ for $t=0$ case) we obtain
\begin{equation}\label{corki}
a_p\in  \{a'_i\}_{i=p}^{m-s}\ \text{ if } p\leq m-k; \ \ \ \ \ b_q\in  \{b'_j\}_{i=q}^{n-s}\ \text{ if } q\leq n-k.
\end{equation}

The assumption  $\nu\leq \lambda$ gives
\begin{equation}\label{apbqq}
\lambda-\nu=(a_p-a'_p)\vareps_p+(b_q-b'_q)\delta_q+\sum_{i>p} (a_i-a'_i)\vareps_i+\sum_{j>q} (b_j-b'_j)\delta_j\in\mathbb{N}\Sigma.
\end{equation}
Consider the case when $\vareps_p-\delta_q\in\Delta^+$.
Then (\ref{apbqq}) implies
$a_p>a'_p$.  Since  $\nu'$ is dominant
$$a'_i\leq a'_p<a_p\ \ \text{ for } i=p,p+1,\ldots, m-s,$$
so $p\geq m-k$ by (\ref{corki}). Notice that $\delta_{n-j}-\vareps_{m-i}\in\Delta^+$
for $j>i$,  so the assumption $\vareps_p-\delta_q\in\Delta^+$ gives
$n-q\leq m-p$ and thus implies (\ref{kmpnq}). For
the remaining case $\delta_q-\vareps_p\in\Delta^+$ the proof is similar.
 \end{proof}

Taking $s=0$ we obtain the following corollary, which is a reformulation of Lemma 6.2 in~\cite{Skw}.

\subsection{}
\begin{cor}{corst}
Assume that  $\lambda,\nu$ are dominant weights, $\lambda$ is stable and
$$\core(\lambda)=\core(\nu)\ \ \ \ \nu\leq\lambda.$$
Then $\nu$ is stable.
\end{cor}

\subsection{}
\begin{cor}{propDsstable}
If $N\in \cF^g$ is stable, then $\DS_x(N)$ is stable.
\end{cor}
\begin{proof}
Assume that $[\DS_x(N):L_{\fg_x}(\nu')]\not=0$. Then $\nu'$ is dominant and
$$\core(\nu')=g.$$

Let $\ol{v}$ be a  vector in $\DS_s(N)=N^x/xN$ which has weight $\nu'$.
 By~\cite{Gaugusta}, Lem. 2.3 we can choose a preimage $v$ of
$\ol{v}$ in the space
$$\sum_{\mu\in X} N_{\mu}^x,\ \ \text{ where }
X:=\{\nu\in\Omega(N)|\ \nu|_{\fh_x}=\nu', \ \  (\nu|S_s)=0\}.$$
Take $\nu\in X$. Since $\nu\in\Omega(N)$ there exists a stable dominant weight
$\lambda$ (a maximal weight in $\Omega(N)$) such that
$$\nu\leq \lambda,\ \ \ \ \ \core(\lambda)=g.$$
In the light of \ref{rest} the condition
$(\nu|S_s)=0$ implies
$$\core(\nu)=\core(\nu|_{\fh_x})=g.$$
By \Lem{lemst}, $\nu'$ is stable.
\end{proof}

\section{Reduction to principal blocks}\label{Fgp}
In this section we reduce the computation of multiplicities
$[\DS_s(L(\lambda)):L_{\fg_x}(\nu')]$  to the case of principal blocks.

In this section  $g$ stands for a core diagram of type $t$ for 
 $\fg:=\osp(2m+t|2n)$.  Let $\mu\in \Lambda^{(t)}_{m-k|n-k}$
be the typical weight corresponding to $g$.

\subsection{Notation}
We denote by $\cF^{pr}(\fg_k)$ the principal block
for $\fg_k$:  
$$\cF^{pr}(\fg_k)=\left\{\begin{array}{ll}
\cF^{\emptyset}(\osp(2k+t|2k)\ \text{ for }t=0,1\\
\cF^{>}(\osp(2k+2|2k) \text{ for }t=2.
\end{array}\right.$$

We denote by $\cF^{g}_{st}$ the subcategory of $\cF^g(\fg)$ consisting
 of stable modules. This category is zero if and only if
 the zero position of $g$ is non-empty.

 We denote by $\Lambda^+(g;i)$ the set of diagrams
$f$ with the following properties:

$\ \ \ \ \ \  \core(f)=g$ and
all symbols $\times$ lie in the positions $0,\ldots, i$.

We denote by $\cF^{g}_i(\fg)$  the Serre subcategory of
$\cF(\fg)$ generated by the modules $L(\lambda)$
with $\lambda\in \Lambda^+(g;i)\cap \Lambda_{m|n}$
and denote by   $\cF^{pr}_i(\fg_k)$ the corresponding subcategory of $\cF^{pr}(\fg_k)$.
Note that $\cF^g(\fg)$ can be viewed as  a ``limit''
of subcategories $\cF^{g}_i(\fg)$.

\subsection{The functor $\Res$ for $(t;k)\not=(0;0)$}\label{notres}
We assume that  $\cF^g(\fg)$ is a non-principal block and
$\cF^g_{st}(\fg)\not=0$. This means that we exclude the case $t=0,k=0$
(in this case  $\cF^g(\fg)$ is a direct sum of two blocks) and that
$g$ has a non-empty symbols
at a non-zero position. One has
$$\cF^{g}_{st}(\fg)=\cF^{g}_q(\fg),$$
where $q+1$ is the coordinate of the first occupied non-zero
position in $g$.

We retain the notation of \ref{funcnot}. Fix $z\in\fh$ such that $\alpha(z)=0$
for $\alpha\in\Delta(\fg_k)$ and  $\alpha(z)\in\mathbb{R}_{>0}$ for $\alpha\in
\Delta^+\setminus \Delta(\fg_k)$. Then
$$\fg^z=\fg_k+\fh=\fg_k\times\fh'',$$
where $\fh''$ is the centralizer of $\fg_k$. 
Set 
$a:=(\mu-\rho)(z)$
and define the functor $\Res:=\Res_a$ using the construction of~\ref{functdef}
for $\fl:=\fg_k$ and $\mu$ as above.
%We use the construction of
%the functor $\Res$ given in~\ref{funcnot} for the pair $\fg_k\subset \fg$.
%Namely, we fix any $z\in\fh$ given by
%$$\alpha(z)=0\ \text{ for }\alpha\in\Sigma_k;\ \
%\alpha(z)\in\mathbb{R}_{>0}\ \text{ for }\alpha\in\Sigma\setminus\Sigma_k$$
%and for a $\fg$-module $N$ we view
%$$\Res(N):=\{v\in N|\ \  zv=av\},\ \text{ where }\ a:=(\mu-\rho)(z)$$
%as a $\fg_k$-module.

\subsubsection{}
\begin{prop}{corstabeq}
The functor
$$\Res:\ \cF^{g}_{st}(\fg)\iso \cF^{pr}_q(\fg_k)$$
is  an equivalence of the categories and
$$\Res (L(f))=L(\howl(f))$$
for each stable $f\in\Lambda^{(t)}_{m+\ell|n}$ with $\core(f)=g$.
\end{prop}
\begin{proof}
Take
$\fh':=\fg_k\cap\fh$. Then $\fh=\fh'\times\fh''$. Setting
$$E:=\{\vareps_i\}_{i=1}^{m+\ell},\ \ D:=\{\delta_i\}_{i=1}^n,\ \
E':=\{\vareps_{i}\}_{i=m+1-k}^{m+\ell},\ \ D':=\{\delta_i\}_{n+1-k}^n.$$
 we  see that
$(\fh')^*$ is spanned by $E'\cup D'$ and
$\ (\fh'')^*$ is spanned by $(E\setminus E')\cup (D \setminus D')$.

Let
$\upsilon:\fh^*\to (\fh')^*$,  $\upsilon'':\fh^*\to (\fh'')^*$
be the projections given by the  decomposition $\fh=\fh'\oplus \fh''$. Notice that
$\mu\in (\fh'')^*$ and set
$$A:=\{\lambda\in\Lambda^+_{m+\ell|n}|\ \core(\lambda)=g,\ \upsilon''(\lambda)=
\mu-\upsilon''(\rho) \},\ \ \  A':=\upsilon(A).$$
One has
$$\begin{array}{rl}
A & =\{\lambda\in\Lambda^+_{m+\ell|n}|\ \core(\lambda)=g,\ \upsilon''(\lambda)=\mu-\upsilon''(\rho)\}\\
& =\{\lambda\in\Lambda^+_{m+\ell|n}|\ \core(\lambda)=g,\ \upsilon''(\lambda+\rho)=\mu\}.
\end{array}$$
In the light of~\ref{stablecore} we get
$A=\Lambda^+_{m+\ell|n}(g;q)$, so
$\cF^g_{st}(\fg)=\cF^{g}_q(\fg)=\cF(A)$.
By~\Cor{corst}, $\Lambda^+_{m+\ell|n}(g;q)$
satisfies the assumption~(\ref{assmstable}). By~\ref{gequivl}, 
$\Res$ provides an equivalence of categories
$\cF^{g}_{st}(\fg)\iso \cF(A')$.

Take $\lambda\in A$. Since $\lambda$ is stable,~\ref{howlst} implies
that
$$\lambda+\rho=\mu+\howl(\lambda)+\rho_k,$$
where $\rho_k$ is the Weyl vector of $\fg_k$.
Note that $\rho_k=\upsilon(\rho)$.

Since $\mu\in (\fh'')^*$ and  $\howl(\lambda)\in (\fh')^*$, we have
 $\upsilon(\lambda)=\howl(\lambda)$, so
$$A'=\upsilon( \Lambda^+_{m+\ell|n}(g;q)) =\Lambda^+(g_{pr};q),$$
where $g_{pr}$ is empty for $t=0,1$ and $g_{pr}=>$ for $t=2$.
Hence  $\cF(A')=\cF^{pr}_q(\fg_k)$. This completes the proof.
\end{proof}

\subsubsection{}
Retain the notation of \ref{mapdex}.

\begin{cor}{corext}
For each $\lambda,\nu\in\Lambda^{(t)}_{m|n}$ one has 
 $$\pari(\lambda)=\pari(\nu)\ \ \Longrightarrow\ \ \Ext^1(L(\lambda), L(\nu))=0.$$
 \end{cor}
\begin{proof}
For the core-free diagrams the assertion is established in~\cite{Gdex}.
The general case follows from~\Prop{corstabeq} (note that
$\Ext^1(L(\lambda), L(\nu))=0$ if  $\lambda$ is typical).
\end{proof}

\subsubsection{Remark}
Consider the case $k=0$ (and $t\not=0$). In this case
$\cF^{g}_{st}(\fg)=\cF^g(\fg)$ is a typical block containing
$L(\mu)$;
this block is isomorphic to the category of finite-dimensional even vector spaces.

For $t=1$ one has $\fg_k=0$ and $\cF^{pr}_q(\fg_k)=\cF(\fg_k)$
is  the category of finite-dimensional even vector spaces
($\Res(L(\mu))=L(\emptyset)=\mathbb{C}$).

For $t=2$ one has $\fg_k=\mathbb{C}$ and
 $\cF(\fg_k)$ is
the category of finite-dimensional $\fg_k$-module
with the zero action of $\fg_k$ (since $t=2$), i.e.
($\Res(L(\mu))$ is the trivial $\fg_k$-module).

\subsection{$\DS_x$ and $\Res$}\label{DSResx}
Take $x:=x_s$. Since $\DS_s(\cF^g)=0$ for $s>k$ we assume
$$0<s\leq k.$$

By \ref{xyh} one has $x\in \fg_s\subset\fg_k$.

Let $(\fg_x)_{i}$ be the subalgebra of $\fg_x$ constructed for $\fg_x$
in the same way as $\fg_i$ to $\fg$.
Consider the base $\Sigma^x_{k-s}\subset\Delta^+(\fg_x)$.
We define the functors $\Res,\Res_x$ as in \ref{notres}: we take
$$\begin{array}{llll}
z\in\fh\ \ \text{ such that }& \alpha(z)=0\ \text{ for }\alpha\in\Sigma_k, & &  \alpha(z)=1\
\text{ for }\alpha\in\Sigma\setminus\Sigma_k\\
z_x\in\fh_x\ \ \text{ such that }&
\alpha(z_x)=0\ \text{ for }\alpha\in\Sigma^x_{k-s}, & &
\alpha(z_x)=1\ \text{ for }\alpha\in(\Sigma^x\setminus\Sigma^x_{k-s})
\end{array}$$ and set
$a:=(\mu-\rho)(z)$, $a_x:=(\mu-\rho_x)(z_x)$.
We obtain $\ \Res: \cF^{g}_{st}(\fg)\longrightarrow\ \cF_{q}^{pr}(\fg_k)$ and
$$\begin{array}{l}
\Res_x:  \tilde{\cF}^{g}_{st}(\fg_x)
\ {\longrightarrow} \ \tilde{\cF}^{pr}_q((\fg_{x})_{k-s})\end{array}$$
 given by
$$\Res(N):=\{v\in N|\ zv=a\},\ \ \ \ \ \ \
\Res_x(N):=\{v\in N|\ z_xv=a_xv\}.$$

\subsubsection{}
\begin{lem}{lem44}
One has  $\DS_x(\fg_{k})=(\fg_x)_{k-s}$. In addition,
except for the case $t=0, k=s$, we have $z_x=z, a_x=a$.
\end{lem}
\begin{proof}
The partial order $\geq$ gives a total order on the standard basis
of $\fh^*$, i.e.
$$\{\vareps_i\}_{i=1}^{m+\ell}\cup \{\delta_j\}_{j=1}^n=
\{e_i\}_{i=1}^{m+n+\ell},\ \ e_1<e_2<\ldots<e_{m+n+\ell}$$
(where $e_{1}=\delta_n$ for $\osp(2m+1|2n)$ and $e_1=\vareps_{m+\ell}$
for $\osp(2m|2n)$). Let $\{e_i^*\}_{i=1}^{m+n+\ell}$ be the dual basis
of $\fh$.

Denote the base of $\DS_x(\fg_k)$ by $(\Sigma_k)^x$.
For $k=s$ we have $(\Sigma_{k})^x=\emptyset=(\Sigma^x)_0$.

In the $t=0$-case one has
$$\begin{array}{l}
\Sigma=\{e_1+e_2,e_2-e_1,e_3-e_2,\ldots, e_{m+n}-e_{m+n-1}\},\\
\Sigma_k=\{e_1+e_2,e_2-e_1,\ldots, e_{2k}-e_{2k-1}\}\\
\Sigma^x=\{e_{2s+1}+e_{2s+2},e_{2s+2}-e_{2s+1},\ldots,e_{m+n}-e_{m+n-1}\}\\
(\Sigma^x)_{k-s}=\{e_{2s+1}+e_{2s+2},e_{2s+2}-e_{2s+1},\ldots,e_{2k}-e_{2k-1}\}=(\Sigma_k)^x.
\end{array}$$

In the $t=1$-case one has
$$\begin{array}{l}
\Sigma=\{e_1,e_2-e_1,e_3-e-2,\ldots, e_{m+n}-e_{m+n-1}\},\\
\Sigma_k=\{e_1,e_2-e_1,\ldots, e_{2k}-e_{2k-1}\}\\
\Sigma^x=\{e_{2s+1},e_{2s+2}-e_{2s+1},e_{2s+3}-e_{2s+2},\ldots,e_{m+n}-e_{m+n-1}\}\\
(\Sigma^x)_{k-s}=\{e_{2s+1},e_{2s+2}-e_{2s+1},\ldots,e_{2k}-e_{2k-1}\}=(\Sigma_k)^x.
\end{array}$$
In both cases
$$\Sigma\setminus\Sigma_k=\{e_{2k+1}-e_{2k},e_{2k+2}-e_{2k+1},\ldots,
e_{m+n}-e_{m+n-1}\}$$
and
$$z=e^*_{2k+1}+2e^*_{2k+2}+\ldots+(m+n-2k) e^*_{m+n}.$$
One readily sees that $z=z_x$ (except for $t=0$ and $k=s$).

In the $t=2$-case one has
$$\begin{array}{l}
\Sigma=\{e_1+e_2,e_2-e_1,e_3-e_2,\ldots, e_{m+n+1}-e_{m+n}\},\\
\Sigma_k=\{e_1+e_2,e_2-e_1,\ldots, e_{2k}-e_{2k-1}, e_{2k+1}-e_{2k}\}\\
\Sigma^x=\{e_1+e_{2s+2},e_{2s+2}-e_{1},e_{2s+3}-e_{2s+2},
\ldots,e_{m+n+1}-e_{m+n}\}\\
(\Sigma^x)_{k-s}=\{e_{1}+e_{2s+2},e_{2s+2}-e_{2s+1},\ldots,e_{2k+1}-e_{2k}\}=(\Sigma_k)^x.
\end{array}$$
Therefore
$$z=e^*_{2k+2}+2e^*_{2k+3}+\ldots+(m+n-2k) e^*_{m+n+1}=
z_x.$$

By~\ref{rest}, $\rho_x=\rho|_{\fh_x}$; since
$z_x\in\fh_x$  we get
$$a_x-a=(\rho_x-\rho)(z_x)=0.$$
Finally, $(\Sigma^x)_{k-s}=(\Sigma_k)^x$ gives
 $\DS_x(\fg_{k})=(\fg_x)_{k-s}$
as required.
\end{proof}

\subsubsection{}\label{commu}
Assume that $0<s\leq k$ and $s\not=k$ for $t=0$.
Combining \Lem{lem44} and \ref{propDsstable}, \ref{corstabeq}
we obtain the following diagram
\begin{equation}\label{commdiag}
\begin{array}{ccc}
\cF^{g}_{st}(\fg) & \overset{\Res}{\longrightarrow} & \cF^{pr}_q(\fg_k)\\
\ & &\\
\DS_x\downarrow & & \DS_x\downarrow\\
\ & &\\
\tilde{\cF}^{g}_{st}(\fg_x) &  \overset{\Res_x}{\longrightarrow} & \tilde{\cF}^{pr}_q((\fg_{x})_{k-s})
\end{array}\end{equation}
where $\Res,\Res_x$ are equivalence of categories.

Let us show that this diagram is commutative.
Take $N\in\cF^{g;q}(\fg)$.  Since $z=z_x\in\fh_x$,
the spaces $N^x$ and $xN$ are $z$-stable, so
$$\Res_x(\DS_x(N))=(\DS_x(N))_a=(N^x)_a/(xN)_a.$$
On the other hand,
$$\DS_x(\Res(N))=\DS_x(N_a)=(N_a)^x/(xN_a).$$
 Since $[x,z]=0$ one has $(xN)_a=x(N_a)$ and $(N^x)_a=(N_a)^x$.
Hence
$\Res_x(\DS_x(N))=\DS_x(\Res(N))$
as required.

\subsection{The case $t=0$ and  $k=s$}\label{Resremain}

Consider the case $t=0$ and $k=s>0$.
By~\ref{corstabeq}  we have
 $\Res: \cF^{g}_{st}(\fg)\iso\cF^{pr}_q(\fg_k)$.
Note that
\begin{equation}\label{Resnik}
\Res(N):=\{v\in N|\ zv=\mu(z)v\},\end{equation}
 where $z$ is as in \ref{DSResx}.
From  the proof of \Lem{lem44} we see that $z\in\fh_x$.

Note that $k-s=0$, so $(\fg_x)_{k-s}=0$ and $\Fin((\fg_x)_{k-s})$
is the category of finite-dimensional supervector spaces, which we denote by
 $sVect$. Define
$$\Res_x: \Fin(\fg_x)\to \ sVect$$ 
by formula (\ref{Resnik}). Using the arguments of~\ref{commu}  we obtain
the following commutative diagram

$$
\begin{array}{ccc}
\cF^{g}_{st}(\fg) & \overset{\Res}{\longrightarrow} & \cF^{pr}_q(\fg_k)\\
\ & &\\
\DS_x\downarrow & & \DS_x\downarrow\\
\ & &\\
\tilde{\cF}^{g}_{st}(\fg_x) &  \overset{\Res_x}{\longrightarrow} & sVect
\end{array}$$

Observe that  $\tilde{\cF}^{g}_{st}(\fg_x)=\tilde{\cF}^{g}(\fg_x)$
is a semisimple category and
$$\Irr\bigl(\tilde{\cF}^{g}(\fg_x)\bigr)=\{
L',\Pi(L'), (L')^{\sigma_x}, \Pi((L')^{\sigma_x})\},\ \text{ where }
L':=L_{\fg_x}(\mu).$$
The functor $\Res_x$ commutes with $\Pi$; clearly,
$$\Res_x\bigl(L_{\fg_x}(\mu)\bigr)=L(\emptyset)=\mathbb{C}.$$

The eigenvalues of $z$ on $(L')^{\sigma_x}=L_{\fg_x}(\sigma(\mu))$ lie
in the set $\sigma_x(\mu)(z)-\mathbb{N}$.
Recall that   $t=0$
and $\mu$ has the diagram with the sign $+$. Using the notation
of the proof of \Lem{lem44} we have
$$(\mu-\sigma_x(\mu))(z)=2(\mu|e^*_{2k+1})>0,$$
so
$$\Res_x\bigl((L')^{\sigma_x})\bigr)=0.$$

\subsection{Graded multiplicity}\label{gradedmult}
Retain notation of~\ref{d0d1}.
We fix $r$ and denote the graded multiplicity
$[\DS_r(L(f)):L(f')]$ by $\frac{f}{f'}$.

Note that $\frac{f}{f'}=0$ if $\core(f)\not=\core(f')$
or $\atyp f-\atyp f'\not=r$ (where $\atyp$ stands for the atypicality, i.e.
the number of the symbols $\times$ in the diagram).

\subsubsection{Case: $\osp(2n+1|2n)$: switch functor}\label{switch}
By \cite{GS}, Lemma 19
the translation functor $T^{\emptyset}_{\emptyset}$ ("switch functor")
acts  on simple modules $L(\mu)\in\cF^{\emptyset}$
as follows:
$$T^{\emptyset}_{\emptyset}(L(\mu))=L(\mu^{sw}),$$
where the diagram of $\mu^{sw}$ is obtained form the diagram
of $\mu$ by sign change ($\mu^{sw}=\mu$
if the diagram of $\mu$ does not have a sign). Since $\DS$ commutes
with the translation functors, we get
$$[\DS_r(L(\lambda)):L(\nu)]=[\DS_r(L(\lambda^{sw})):L(\nu^{sw})].$$

\subsubsection{}
Let $\sigma(f)$ be the diagram  obtained from 
$f$ by the change of sign (for $\osp(2m|2n)$ this notation was used before).

\subsubsection{}
\begin{cor}{corka}
\begin{enumerate}
\item $\frac{f}{f'}=\frac{\howl(f)}{\howl(f')}$;

\item  $\frac{\sigma(f)}{\sigma(f')}=\frac{f}{f'}$;

\item for $t=0$ and $r=1$ one has $\frac{f}{\sigma(f')}=\frac{f}{f'}$.
\end{enumerate}
\end{cor}
\begin{proof}
Combining~\ref{commu} and~\Cor{corstabi} we obtain (i) 
for all cases except $t=0$ with $\howl(f')=\emptyset$.

For a diagram $f$ denote by $\lambda(f)$ the corresponding dominant weight.
 For $t=2$ one has $\sigma(f)=f$. For $t=1$ one has 
  $\lambda(\sigma(f))=(\lambda(f))^{sw}$ and (ii) follows from~\ref{switch}.

Consider the remaing  case $t=0$. The formulae (ii), (iii) follow from~\Lem{lemsigma}. 
For (i) take $f'$ such that $\howl(f')=\emptyset$.
By \ref{Resremain}, $\frac{f}{f'}=\frac{\howl(f)}{\howl(f')}$ if
the diagram of $f'$ has the sign $+$.
Assume that $f'$ has the sign $-$. Since $\sigma(f')$ has sign $+$ we have 
$$\frac{\sigma(f)}{\sigma(f')}=
\frac{\howl(\sigma(f))}{\howl(\sigma(f'))}.$$
Using (ii) and the formulae
$\howl(\sigma(f))=\sigma(\howl(f))$ 
we get
$$\frac{f}{f'}=\frac{\sigma(f)}{\sigma(f')}=\frac{\howl(\sigma(f))}{\howl(\sigma(f'))}=
\frac{\sigma(\howl(f))}{\sigma(\emptyset)}=\frac{\howl(f)}{\emptyset}
$$
which implies (i). 

\end{proof}

\section{Recursive formulae for $[\DS_x(L(\lambda)):L_{\fg_x}(\nu)]$}\label{sect4}

Using that $DS$ commutes with translation functors, we establish recursive formulas for the multiplicities $[\DS_x(L(\lambda)):L_{\fg_x}(\nu)]$. This will ultimately allow us to reduce the computation of the multiplicities to the atypicality 1 case, i.e. to $\mathfrak{osp}(2|2)$, $\mathfrak{osp}(3|2)$ and $\mathfrak{osp}(4|2)$. For a similar reduction in the $\mathfrak{gl}(m|n)$-case see \cite{HW}.

We fix $r$ and retain the notation of~\ref{gradedmult}.
%
%We always assume that $\lambda\in\Lambda^+_{m+\ell|n}$
%and $\nu\in\Lambda^+_{m+\ell-r|n-r}$;
%for instance, using the notation $\frac{f}{\emptyset}$
%we assume that the diagram $f$ correspond to $\lambda\in\Lambda^+_{r|r}$.
%

\subsection{Translation functors}\label{DSTg}
Recall that $\DS$ commutes with the translation functors. Let $g_0,g_1$ be two core diagrams. Consider the translation functors
$$T^{g_0}_{g_1}: \tilde{\cF}(\fg)^{g_1}\to \tilde{\cF}(\fg)^{g_0},\ \ \
T^{g_0}_{g_1}: \tilde{\cF}(\fg_x)^{g_1}\to \tilde{\cF}(\fg_x)^{g_0}.$$

Since translation functors are exact, they induce morphisms on the Grothendieck ring. For  $N\in\cF(\fg)^{g_1}$ and
 $L'\in \Irr(\fg_x)^{g_0}$ we obtain
\begin{equation}\label{LL'}\begin{array}{ll}
[\DS_r(T^{g_0}_{g_1}(N)):L']&=[T^{g_0}_{g_1}(DS_r(N)):L']\\
&=\displaystyle\sum_{L_1\in\Irr(\fg_x)^{g_1}}
[DS_r(N):L_1][T^{g_0}_{g_1}(L_1):L'].\end{array}\end{equation}

\subsubsection{}\label{trans}
Assume that the number of core symbols in $g_0$ is larger than
 the number of core symbols in $g_1$, i.e.
the atypicality of $g_1$ is larger than the atypicality of $g_0$.
We will use the following results of \cite{GSBGG}, Lemmatta 7, 13, 14:

the image of a simple module $L\in\Irr(\fg_x)^{g_1}$ is either zero or simple;

for $L_1,L_2\in\Irr(\fg_x)^{g_1}$ with $T^{g_0}_{g_1}(L_1)\cong T^{g_0}_{g_1}(L_2)\not=0$
one has $L_1\cong L_2$.

\subsubsection{}\label{ES-trans} These results will be sufficient for us. A complete description of translation functors on irreducible modules and projective covers can be obtained from \cite{ES},\cite{ES2},\cite{ES3}.

\subsubsection{}\label{g0g1}
Take $g_0,g_1$ satisfying the assumption in \ref{trans}. For each 
$L' \in\Irr(\fg_x)^{g_0}$
there exists at most one (up to isomorphism)  $L_1\in \Irr(\fg_x)^{g_1}$ such that $T^{g_0}_{g_1}(L_1)\cong L'$.
Then (\ref{LL'}) gives for $N\in\cF(\fg)^{g_1}$ 
$$[\DS_r(T^{g_0}_{g_1}(N)):L']=[T^{g_0}_{g_1}(DS_r(N)):L']=[DS_r(N):L_1].$$
In particular, 
\begin{equation}\label{DSrN}
[DS_r(N):L_1]\not=0\ \ \Longrightarrow\ \ T^{g_0}_{g_1}(N)\not=0.
\end{equation}

Now take $N:=L(\lambda)$ and $\nu$ such that
$[\DS_r(L(\lambda)):L_{\fg_x}(\nu)]\not=0$
and $T^{g_0}_{g_1}(L_{\fg_x}(\nu))\not=0$. Then, by~\ref{trans},
$$T^{g_0}_{g_1}(L_{\fg_x}(\nu))=L_{\fg_x}(\nu_1)$$
for some $\nu_1$
and $T^{g_0}_{g_1}(L(\lambda))$ is either zero or simple. By above,
$T^{g_0}_{g_1}(L(\lambda))\not=0$, so
$$T^{g_0}_{g_1}(L(\lambda))=L(\lambda_1)$$
for some $\lambda_1$.  Using (\ref{LL'}) and \Cor{corka} we conclude
\begin{equation}\label{DST}
\frac{f(\lambda)}{f(\nu)}=\frac{\howl (f(\lambda_1))}{\howl(f(\nu_1))},
\end{equation}
where $f(\lambda)$ stands for the weight diagram of $\lambda$.

\subsection{Translation functors $T_u$}
\label{pos}
We describe some translation functors via their effect (called \emph{elementary change} in~\cite{GSBGG}, Section 6.3) on core/weight diagrams. Let $a$ be a non-negative integer.
For each diagram $f$ we denote by $\pos_a(f)$
 the subdiagram corresponding to the positions
$a,a+1$.
For a core diagram
$g_1$  with $\pos_a(g_1)=\circ\circ$  we denote by $\phi'_a(g_1)$
the core diagram obtained from $g_1$ by changing
 $\pos_a(g_1)=\circ\circ$ to $\pos_a(g)=><$; for instance,
 $$\phi'_1(<\circ\circ>)==<><>.$$
We denote by $T_a$ the functor which acts as $T^{\phi'_a(g_1)}_{g_1}$
on $\tilde{\cF}^{g_1}$ with $\pos_a(g_1)=\circ\circ$ and by zero on
$\tilde{\cF}^{g'}$ with $\pos_a(g')\not=\circ\circ$. Note that
$T_a$ reduces the atypicality by $1$.

\subsubsection{}\label{posu}
Take $u>0$.
For a diagram $f$ with $\pos_u(f)=\times\circ$ we define two new diagrams:

the diagram
$\phi_u(f)$   obtained from $f$ by changing
$\times\circ $ in the positions $(u,u+1)$ to $><$ (i.e.,
$f$ and $\phi_u(f)$ have
 the same signs and the same symbols in all positions except $u,u+1$
and $\pos_u(\phi_u(f))=><$);

 the diagram $\shr_u(f)$   obtained from $f$ by ``shrinking''
$\times\circ $ at the positions $u,u+1$, i.e.  changing the diagram $f_-\times\circ f_+$ to the diagram
$f_-f_+$;

 For instance, $\phi_1(\times\times\times\circ),\shr_1(\times\times\times\circ)$ are not defined;
 $$\phi_2(\times\times\times\circ)=\times\times><,\ \ \ \ \
\shr_2(\times\times\times\circ)=\times\times.$$

\subsubsection{}
By~\cite{GSBGG}, for $u>0$ one has
$$ T_{u}(L(\mu))=L(\phi_u(\mu))$$
if $\phi_u(\mu)$ is  defined and $T_{u}(L(\mu))=0$ otherwise.

\subsection{}
\begin{cor}{mx}
For $u>0$ one has
$$\frac{f'_u\times\times f'}{f_u\times\circ f}=\frac{f'_u\circ\times f'}{f_u\times\circ f}=\frac{f'_u\circ\circ f'}{f_u\times\circ f}=0$$
(with all possible signs) and
$$\frac{f'_u\times\circ f'}{f_u\times\circ f}=\frac{f'_u f'}{f_u f}$$
where $f_u,f'_u$ stands for the subdiagrams corresponding to the positions
$0,1,\ldots,u-1$.
\end{cor}

\begin{proof} Recall from (\ref{DSrN}) that for irreducible $N = L(\lambda)$ \[ [\DS_r(T^{g_0}_{g_1}(L(\lambda)):T_{g_1}^{g_0}(L_1)] =[DS_r(L(\lambda):L_1].\]
Take $u>0$. Let $\tilde{f}$ be a diagram with 
  $\pos_u(\tilde{f})=\times\circ$. Applying $T_u$ sends $f_u\times\circ f$ to $f_u >< f$. Combining~(\ref{DSrN}) and~\ref{posu}
 we get $\frac{\tilde{f}'}{\tilde{f}}=0$ if $\pos_u(\tilde{f}')\not=\times\circ$;
this gives the first formula.  Combining \ref{posu} and (\ref{DST}) we have
$$\frac{f'_u\times\circ f'}{f_u\times\circ f}=\frac{\howl(f'_u>< f')}{\howl(f_u>< f)}.$$
Using~\Cor{corka} and 
 $\howl(f'_u>< f')=\howl(f'_u f')$ we obtain the second formula.
\end{proof}

\subsection{Reduction to the case $\nu=0$}
Let $f$ be a weight diagram and
$u$ be the coordinate of the rightmost symbol $\times$ in $f$.
Using~\ref{corka} (i) and~\ref{mx} we reduce the computation 
of $\frac{f'}{f}$ to the situation when $f,f'$ are core-free and
$u=0$; in the light of~\ref{corka} (ii) we can assume that
$f$ has sign $-$  for $t=1$ (for $t=0,2$ the diagram
$f$ does not have sign if $u=0$). Notice that in this case
$\lambda(f')=0$, i.e. we reduced the problem to the computation
of $[\DS_r(L(\lambda)): L_{\fg_x}(0)]$ for $\fg=\osp(2n+t|2n)$.

For $\fgl(m|n)$ a similar reasoning reduces the computation
of $[\DS_r(L(\lambda)):L(\nu)]$ to the case when $\nu$ has the empty diagram, i.e.
to the case when $\fg_x=0$ (see also \cite{HW}).
In the $\osp$-case  this is done in Corollaries \ref{prop2m12n},\ref{prop2m2n} below.

\subsection{The case $\osp(2m+1|2n)$}
Consider the case $\fg=\osp(2m+1|2n)$.
We assume, as always, that "the signs
disappear" if the zero position is empty, i.e. $\pm\times^{i-1}f$ stands for
$\circ f$ for $i=1$. 

\subsubsection{}\label{T02m12n}
By \cite{GSBGG},
for the translation functor $T_{0}:=T^{><}_{\emptyset}$ one has
$T_0(L(\mu))\not=0$ if and only if the diagram of $\mu$ has the sign $+$ and
$\pos_0(\mu)=\times^i\circ$ for $i>0$; moreover,
$$T_0(L(+\times^i\circ f)=L(\overset{\ \ \times^{i-1}}{>}<f).$$
Using~(\ref{DSrN})  we obtain
 for $j\geq 1$
$$\frac{\tilde{f}}{+\times^j\circ f}\not=0\ \ \Longrightarrow\ \ \tilde{f}=+\times^p \circ f'\ \text{ for some } p\geq 1.$$

\subsubsection{}
\begin{lem}{prop2m12n}
 For $p,i\geq 1$ one has
$$\begin{array}{cccc}
\frac{+\times^p\circ \circ f}{+\times^i}=\frac{-\times^{p-1}f}{-\times^{i-1}}
 & & &\ \ \ \ 
\frac{+\times^p\circ \times f}{+\times^i}=\frac{+\times^{p}f}{-\times^{i-1}}\end{array}$$
\end{lem}
\begin{proof}
Using~(\ref{DST}) for $T_0$ we obtain for
 $*\in \{\circ, \times\}$
$$\begin{array}{l}
\frac{+\times^p\circ * f}{+\times^i}=\frac{\overset{\ \ \times^{p-1}}{>}< *f}{\overset{\ \ \times^{i-1}}{>}<}.
\end{array}$$
 One has
$$ \howl(\overset{\ \ \times^{p-1}}{>}< *f)=\left\{\begin{array}{ll}
-\times^{p-1}f\ & \text{ if }\ *=\circ \\
+\times^{p}f\ & \text{ if }\ *=\times\end{array}\right.$$
and $\howl(\overset{\ \ \times^{i-1}}{>}<)=-\times^{i-1}$. Now~(\ref{DST}) gives
$$\begin{array}{l}
\frac{+\times^p\circ \circ f}{+\times^i}=\frac{-\times^{p-1}f}{-\times^{i-1}}\ \ \ \ \ \ \ \ \  \frac{+\times^p\circ \times f}{+\times^i}=\frac{+\times^{p}f}{-\times^{i-1}}\end{array}
$$
as required.
\end{proof}

\subsubsection{}
\begin{cor}{cor2m12n}
Take $i\geq 1$. One has
$$\begin{array}{cc}
\frac{-\times^p \underbrace{\circ...\circ}_{2i-1\text{ times}}\times f}{-\times^i}=
\frac{-\times^{p-i+1}f}{\emptyset} & \ \ \ \ \ \ 
\frac{-\times^p \underbrace{\circ...\circ}_{j\text{ times}}f}{-\times^i}=\frac{-\times^{p-i} \underbrace{\circ...\circ}_{j-2i\text{ times}}f}{\emptyset}\ \text{ for }\ j\geq 2i.\end{array}$$
Moreover, $\frac{\tilde{f}}{-\times^i}=0$ if  $\tilde{f}$ is not as above, i.e. 
$\tilde{f}\not=-\times^p \underbrace{\circ...\circ}_{j\text{ times}} f$ for some $p\geq i$ and $j\geq 2i-1$. 
\end{cor}
\begin{proof}
The statement follows by induction from~\Lem{prop2m12n} and~\Cor{corka} (ii).
\end{proof}

\subsection{The case  $\osp(2m|2n)$}\label{2n2n0}
Recall that the simple $OSp(2m|2n)$-modules are in one-to-one
correspondence with the unsigned diagram (see \ref{OSP2m2n}).
For a non-empty diagram $f$ we will use the sign
$\circ f$ for $+\circ f\oplus -\circ f$, i.e. $L(\circ f)$ is a simple
$OSp(2m|2n)$-module which is the direct sum of $\osp(2m|2n)$-modules
$L(+\circ f)$ and $L(-\circ f)$.

For an empty diagram we have
$L(\emptyset)=\mathbb{C}$; notice that
$$\howl(\pm \circ >)=\emptyset;\ \ \ \ \howl(\circ >)\not=\emptyset.$$

\subsubsection{}
For the translation functor $T_{0}:=T^{><}_{\emptyset}$ one has
$T_0(L(\mu))\not=0$ if and only if
$\pos_0(\mu)=\times^i\circ$  for $i>0$ and
$$T_0(L(\times^i\circ f))=L(\overset{\ \ \times^{i-1}}{>}<f).$$
Using \ref{g0g1} we obtain
 for $j\geq 1$
\begin{equation}\label{nunk}
\frac{\tilde{f}}{\times^j\circ f}\not=0\ \ \Longrightarrow\ \ \tilde{f}=\times^p \circ f'\ \text{ for some } p>0.\end{equation}

\subsubsection{}\label{Tecirc}
The translation functor $T^{\circ >}_{>}$ is
given by
$$T^{\circ>}_{>}(L(\overset{\ \ \times^{p}}{>}\circ f)=L(\times^p >f),\ \ \
T^{\circ>}_{>}(L(\overset{\ \ \times^{p}}{>}\times f)=L(\times^{p+1} >f)
 $$
for each $p\geq 0$ (where,  $L(\circ f)=L(+\circ f)\oplus L(-\circ f)$).

\subsubsection{}\label{Tcirc}
The translation functor $T^{>}_{\circ >}$ is  adjoint to the functor $T^{\circ >}_{>}$. 
It is given by
$$T^{>}_{\circ>}(L(\pm \circ >f)=L(>\circ f),\ \ \ \ T^{>}_{\circ>}(L(\times^p >f))=
L(\overset{\ \ \times^{p}}{>}\circ f)\oplus L(\overset{\ \ \times^{p-1}}{>}\times f)
$$
for each $p>0$ (there is a misprint in \cite{GSBGG}).

\subsubsection{}
\begin{lem}{prop2m2n}
Take $i,p\geq 1$.

(i) $\frac{\times^p\circ f}{\times^{i}}=\frac{\overset{\ \ \times^{p-1}}{>}f}{\overset{\ \ \times^{i-1}}{>}}\ \ \  \ \ \ \ \ \ \ \frac{\overset{\ \ \times^{p}}{>}\circ f}{\overset{\ \ \times^{i}}{>}}=\frac{\times^p f}{\times^{i}}.$

(ii) $\frac{\overset{\ \ \times^{p-1}}{>}\times f}{\overset{\ \ \times^{i}}{>}}=0$.

(iii) $\frac{\overset{\ \ \times^{p}}{>}\circ f}{>}=\frac{\overset{\ \ \times^{p-1}}{>}\times f}{>}=
\frac{\times^p f}{\emptyset}$.
\end{lem}

\begin{proof}
Using (\ref{DST}) for $T_0$ we obtain
\begin{equation}\label{eqnik}\begin{array}{l}
\frac{\times^p\circ f}{\times^{i}}=\frac{\overset{\ \ \times^{p-1}}{>}< f}{\overset{\ \ \times^{i-1}}{>}<}=\frac{\overset{\ \ \times^{p-1}}{>}f}{\overset{\ \ \times^{i-1}}{>}}\end{array}\end{equation}
which establishes the first formula.
Using \ref{DSTg} for $T^{\circ>}_{>}$ we obtain
$$\begin{array}{l}
\frac{\times^p f}{\times^i}=\frac{\times^p>f}{\times^i>}=
[\DS_r(T^{\circ >}_{>}(L(\overset{\ \ \times^{p}}{>}\circ f)):L(\times^i>)]\\
=\sum [\DS_r(L(\overset{\ \ \times^{p}}{>}\circ f)):L_1]
[T^{\circ >}_{>}(L_1):L(\times^i >)]\\=
[\DS_r(L(\overset{\ \ \times^{p}}{>}\circ f)):L(\overset{\ \ \times^{i}}{>})]+
[\DS_r(L(\overset{\ \ \times^{p}}{>}\circ f)):L(\overset{\ \ \times^{i-1}}{>}\times)]\\
=\frac{\overset{\ \ \times^{p}}{>}\circ f}{\overset{\ \ \times^{i}}{>}}+
\frac{\overset{\ \ \times^{p}}{>}\circ f}{\overset{\ \ \times^{i-1}}{>}\times}.
\end{array}$$
By \Cor{mx} the second summand in the last formula is zero; this
implies the second formula.
Similarly,  $T^{\circ >}_{>}(L(\overset{\ \ \times^{p-1}}{>}\times f))=
L(\times^p> f)$ implies
$$\begin{array}{l}\frac{\times^p f}{\times^i}=\frac{\overset{\ \ \times^{p-1}}{>}\times f}{
\overset{\ \ \times^{i}}{>}}+\frac{\overset{\ \ \times^{p-1}}{>}\times f}
{\overset{\ \ \times^{i-1}}{>}\times}.\end{array}$$
By (\ref{nunk}) if $f=\times f'$, then $\frac{\times^p f}{\times^i}=0$, so the both summands in the right-hand side are equal to $0$; in particular, $\ \frac{\overset{\ \ \times^{p-1}}{>}\times \times f'}{
\overset{\ \ \times^{i}}{>}}=0$.

If $f=\circ f'$ we have
$$\begin{array}{l}\frac{\times^p \circ f'}{\times^i}=\frac{\overset{\ \ \times^{p-1}}{>}\times \circ f'}{
\overset{\ \ \times^{i}}{>}}+\frac{\overset{\ \ \times^{p-1}}{>}\times \circ f'}
{\overset{\ \ \times^{i-1}}{>}\times}=
\frac{\overset{\ \ \times^{p-1}}{>}\times \circ f'}{
\overset{\ \ \times^{i}}{>}}+\frac{\overset{\ \ \times^{p-1}}{>}f'}
{\overset{\ \ \times^{i-1}}{>}}.
\end{array}$$
Using (i) we conclude
$$\begin{array}{l}\frac{\overset{\ \ \times^{p-1}}{>}\times \circ f'}{
\overset{\ \ \times^{i}}{>}}=0.\end{array}$$
This establishes (ii).

Using~\ref{DSTg} for $T^{\circ>}_{>}$ we obtain
$$\begin{array}{l}
\frac{\times^p f}{\emptyset}=\frac{\times^p>f}{+\circ >}=
[\DS_r(T^{\circ>}_{>}(L(\overset{\ \ \times^{p}}{>}\circ f)):L(+\circ >)]\\
=\sum [\DS_r(L(\overset{\ \ \times^{p}}{>}\circ f)):L_1]
[T^{\circ>}_{>}(L_1):L(+\circ >)]\\=
[\DS_r(L(\overset{\ \ \times^{p}}{>}\circ f)):L(>)]
=\frac{\overset{\ \ \times^{p}}{>}\circ f}{>}
\end{array}$$
and, similarly,
$$\begin{array}{l}
\frac{\times^p f}{\emptyset}=\frac{\overset{\ \ \times^{p-1}}{>}\times f}{>}
\end{array}.$$
This establishes (iii).
\end{proof}

\subsubsection{}
\begin{cor}{cor2m2n}
Take $i\geq 1$. One has $\frac{\times^p\times f}{\times^{i}}=0$ and
$\frac{\times^p\circ f}{\times^{i}}=\frac{\overset{\ \ \times^{p-1}}{>}f}
{\overset{\ \ \times^{i-1}}{>}}$. Moreover,
$$\begin{array}{l}
\ \ \ \ \ \frac{\overset{\ \ \times^{p}}{>} \underbrace{\circ...\circ}_{j\text{ times}}f'}{\overset{\ \ \times^{i}}{>}}=\frac{\overset{\ \ \times^{p-i}}{>} \underbrace{\circ...\circ}_{j-2i\text{ times}} f'}{>}\end{array}$$
and  $\frac{\overset{\ \ \times^{p}}{>}f}{\overset{\ \ \times^{i}}{>}}\not=0$ 
implies that
$p\geq i$ and
the diagram $f$ is as above (i.e., 
$f=\underbrace{\circ...\circ}_{j\text{ times}} f'$
for some $j\geq 2i$).
\end{cor}

\section{Computation of $\DS_1(L)$ in terms of arc diagrams}\label{Sectlast}

\subsection{Arc diagrams}\label{arcs}
We assign to each core-free diagram  an
 {\em arc diagram} $Arc(f)$ ; if $f$ is contains
 at most one $\times$ at the zero position ("$\fgl$-type"), the corresponding
 arc  diagram coincides with the usual arc diagram introduced in \cite{GSBGG}. It differs from the arc- or cup diagrams of \cite{ES}. Advantages of our weight and arc diagrams are that they can immediately be read off from the weight $\lambda$ and describe the effect of $DS$ very nicely. On the other they do not connect to Khovanov algebras of type $D$ and therefore to the Kazhdan-Lusztig combinatoric of parabolic category $\mathcal{O}$.

An {\em arc diagram} is the following data: a diagram $f$, where
the symbols $\times$ at the zero position are drawn vertically and
and a collection of non-intersecting arcs. 
 Each arc connects  one symbol $\times$ (the left end) with
 one or two empty symbols according to the following rules:
 
if the symbol $\times$ has non-zero coordinate, 
the arc connects
 this symbol with one empty symbol;

for $\ell=0$  (i.e., $t=0,1$)  the lowest symbol $\times$ in the zero position,
is connected by an arc with one empty symbol and the
 other symbols $\times$ in the zero position
 are connected with two empty symbols; 

for $\ell=1$ (i.e., $t=2$) all symbols $\times$ in the zero position
 are connected with two empty symbols. 

An empty position in $f$ is called {\em free} in the arc diagram
 if this position is not an end of an arc.

We say that an arc is supported by a symbol $\times$ if
this symbol is the left end of the arc; if the arc
is supported by a  symbol $\times$ with the coordinate $a$
we denote this arc by  $arc(a;b)$ (resp., $arc(0;b_1,b_2)$) where $b$
(resp., $b_1<b_2$) is the coordinate of the right end (resp., right ends)
of the arc.

We remark that we can similarly define an arc diagram for any weight diagram by just fixing and ignoring the core symbols.

\subsubsection{Partial order}
We consider a partial order on the set of arcs by saying that one arc is smaller than
another one if the first one is "below" the second one, that is

$arc(a;b)>arc(a';b')$ if and only if $a<a'<b'<b$;

$arc(0;b_1,b_2)>arc(a';b')$  if and only if $b'<b_2$;

$arc(0;b_1,b_2)>arc(0;b'_1,b'_2)$  if and only if $b_2'<b_2$.

Since the arcs do not intersect, one has
$$\begin{array}{l}
arc(a;b)>arc(a';b')\ \ \Longleftrightarrow\ \ a<a'<b\\
arc(0;b_1,b_2)>arc(a';b')\ \ \Longleftrightarrow\ \ a'<b_2
\end{array}$$
and
any two distinct arcs of the form $arc(0;b_1,b_2),arc(0;b'_1,b'_2)$ are comparable: either 
$arc(0;b_1,b_2)>arc(0;b'_1,b'_2)$  or
$arc(0;b_1,b_2)<arc(0;b'_1,b'_2)$.

\subsubsection{Definition}
We assign to a core-free diagram $f$ the
arc diagram $Arc(f)$
with the following properties:

each symbol $\times$ is the left end of exactly one arc;

there are no free positions under the arcs.

\subsubsection{Example}\label{exaMitka}

%\begin{figure}
%\begin{center}
%\includegraphics[width=0.3\textwidth]{ospdiag.pdf}
%\caption{Arc diagram for $\times^2\circ \times \circ\circ\times\times$}
%\end{center}
%\end{figure}

The weight diagram below
does not have a sign
for $t=0$ ($\fg=\osp(10|10)$) or have one of the signs $\pm$
for $t=1$ ($\fg=\osp(11|10)$).

\begin{center}
%\medskip
 
 \scalebox{0.7}{
\begin{tikzpicture}
 %\draw (-1,0) -- (7,0);
\foreach \x in {} %vee
     \draw[very thick] (\x-.1, .1) -- (\x,-0.1) -- (\x +.1, .1);
\foreach \x in {} %wedge
     \draw[very thick] (\x-.1, -.1) -- (\x,0.1) -- (\x +.1, -.1);
\foreach \x in {0,2,5,6,11} %cross
     \draw[very thick] (\x-.1, .1) -- (\x +.1, -.1) (\x-.1, -.1) -- (\x +.1, .1);
%\foreach \x in {1,3,4,7,8,9,10,12}  \draw[semithick] \circ; %circle
     %\draw[very thick]  node at (0,0) [fill=white,draw,circle,inner sep=0pt,minimum size=6pt]{};
     \draw[very thick]  node at (1,0) [fill=white,draw,circle,inner sep=0pt,minimum size=6pt]{};
     \draw[very thick]  node at (3,0) [fill=white,draw,circle,inner sep=0pt,minimum size=6pt]{};
     \draw[very thick]  node at (4,0) [fill=white,draw,circle,inner sep=0pt,minimum size=6pt]{};
     \draw[very thick]  node at (7,0) [fill=white,draw,circle,inner sep=0pt,minimum size=6pt]{};
     \draw[very thick]  node at (8,0) [fill=white,draw,circle,inner sep=0pt,minimum size=6pt]{};
     \draw[very thick]  node at (9,0) [fill=white,draw,circle,inner sep=0pt,minimum size=6pt]{};
     \draw[very thick]  node at (10,0) [fill=white,draw,circle,inner sep=0pt,minimum size=6pt]{};
     \draw[very thick]  node at (12,0) [fill=white,draw,circle,inner sep=0pt,minimum size=6pt]{};
\foreach \x in {0} %cross
     \draw[very thick] (\x-.1, +0.8) -- (\x +.1, +0.6) (\x-.1, +0.6) -- (\x +.1, +0.8);

\draw (0,-0.5) node {0};
\draw (1,-0.5) node {1};
\draw (2,-0.5) node {2};
\draw (3,-0.5) node {3};
\draw (4,-0.5) node {4};
\draw (5,-0.5) node {5};
\draw (6,-0.5) node {6};
\draw (7,-0.5) node {7};
\draw (8,-0.5) node {8};
\draw (9,-0.5) node {9};
\draw (10,-0.5) node {10};
\draw (11,-0.5) node {11};
\draw (12,-0.5) node {12};

% node[pos=(0, -0,5)]{0};

%%caps,cups
\draw[very thick] [-,black,out=90, in=90](0,0.2) to (1,0.2);
\draw[very thick] [-,black,out=90, in=90](2,0.2) to (3,0.2);
\draw[very thick] [-,black,out=90, in=90](11,0.2) to (12,0.2);
\draw[very thick] [-,black,out=90, in=90](6,0.2) to (7,0.2);
\draw[very thick] [-,black,out=90, in=90](5,0.2) to (8,0.2);
\draw[very thick] [-,black,out=90, in=90](0,+0.9) to (4,0.2);
\draw[very thick] [-,black,out=90, in=90](0,+0.9) to (9,0.2);

%\foreach \x in {} \draw + at (-1,0);

\end{tikzpicture} }
\medskip

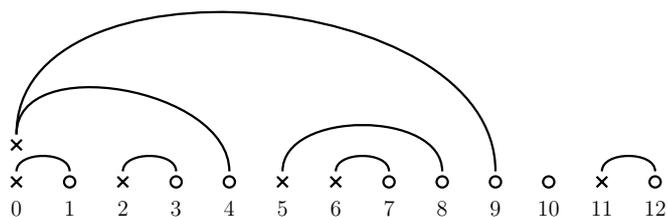
\text{Arc diagram for $\times^2\circ \times \circ\circ\times\times \circ \circ \circ \circ\times \circ$}
\end{center}

\subsubsection{Examples}

In the following diagrams we ignore signs. Note that in case the arc connects $\times$ with two empty symbols, this still counts as one arc (also for Theorem \ref{thmDS1osp}).

For $\osp(5|4)$ the arc diagram of $+\times\circ \circ \times$ is given by 

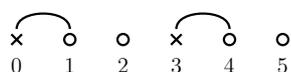
\begin{center}
%\medskip
 
 \scalebox{0.7}{
\begin{tikzpicture}
 %\draw (-1,0) -- (7,0);
\foreach \x in {} %vee
     \draw (\x-.1, .2) -- (\x,0) -- (\x +.1, .2);
\foreach \x in {} %wedge
     \draw (\x-.1, -.2) -- (\x,0) -- (\x +.1, -.2);
\foreach \x in {0,3} %cross
     \draw[very thick] (\x-.1, .1) -- (\x +.1, -.1) (\x-.1, -.1) -- (\x +.1, .1);
%\foreach \x in {} %circle
     %\draw  node at (1,0) [fill=white,draw,circle,inner sep=0pt,minimum size=6pt]{};
     %\draw  node at (2,0) [fill=white,draw,circle,inner sep=0pt,minimum size=6pt]{};
     %\draw  node at (4,0) [fill=white,draw,circle,inner sep=0pt,minimum size=6pt]{};

\draw (0,-0.5) node {0};
\draw (1,-0.5) node {1};
\draw (2,-0.5) node {2};
\draw (3,-0.5) node {3};
\draw (4,-0.5) node {4};
\draw (5,-0.5) node {5};
%\draw (6,-0.5) node {6};

     \draw[very thick]  node at (1,0) [fill=white,draw,circle,inner sep=0pt,minimum size=6pt]{};
     \draw[very thick]  node at (2,0) [fill=white,draw,circle,inner sep=0pt,minimum size=6pt]{};
     \draw[very thick]  node at (4,0) [fill=white,draw,circle,inner sep=0pt,minimum size=6pt]{};
     \draw[very thick]  node at (5,0) [fill=white,draw,circle,inner sep=0pt,minimum size=6pt]{};
     %\draw[very thick]  node at (6,0) [fill=white,draw,circle,inner sep=0pt,minimum size=6pt]{};

%%caps,cups
\draw[very thick] [-,black,out=90, in=90](0,0.2) to (1,0.2);
\draw[very thick] [-,black,out=90, in=90](3,0.2) to (4,0.2);
%\draw [-,black,out=90, in=90](-2,0) to (3,0);
%\draw [-,black,out=90, in=90](-3,0) to (4,0);

%\foreach \x in {} \draw + at (-1,0);

\end{tikzpicture} }
%\medskip
\end{center}

The arc diagrams of $-\times^2\underbrace{\circ...\circ}_{j\text{ times}} \times$ are the following:

\begin{center}
%\medskip
 
 \scalebox{0.7}{
\begin{tikzpicture}
 %\draw (-1,0) -- (7,0);
\foreach \x in {} %vee
     \draw[very thick] (\x-.1, .1) -- (\x,-0.1) -- (\x +.1, .1);
\foreach \x in {} %wedge
     \draw[very thick] (\x-.1, -.1) -- (\x,0.1) -- (\x +.1, -.1);
\foreach \x in {0,1} %cross
     \draw[very thick] (\x-.1, .1) -- (\x +.1, -.1) (\x-.1, -.1) -- (\x +.1, .1);
%\foreach \x in {1,3,4,7,8,9,10,12}  \draw[semithick] \circ; %circle
     %\draw[very thick]  node at (0,0) [fill=white,draw,circle,inner sep=0pt,minimum size=6pt]{};
     \draw[very thick]  node at (2,0) [fill=white,draw,circle,inner sep=0pt,minimum size=6pt]{};
     \draw[very thick]  node at (3,0) [fill=white,draw,circle,inner sep=0pt,minimum size=6pt]{};
     \draw[very thick]  node at (4,0) [fill=white,draw,circle,inner sep=0pt,minimum size=6pt]{};
     \draw[very thick]  node at (5,0) [fill=white,draw,circle,inner sep=0pt,minimum size=6pt]{};
     %\draw[very thick]  node at (8,0) [fill=white,draw,circle,inner sep=0pt,minimum size=6pt]{};
     %\draw[very thick]  node at (9,0) [fill=white,draw,circle,inner sep=0pt,minimum size=6pt]{};
     %\draw[very thick]  node at (10,0) [fill=white,draw,circle,inner sep=0pt,minimum size=6pt]{};
     %\draw[very thick]  node at (12,0) [fill=white,draw,circle,inner sep=0pt,minimum size=6pt]{};
\foreach \x in {0} %cross
     \draw[very thick] (\x-.1, +0.8) -- (\x +.1, +0.6) (\x-.1, +0.6) -- (\x +.1, +0.8);

\draw (0,-0.5) node {0};
\draw (1,-0.5) node {1};
\draw (2,-0.5) node {2};
\draw (3,-0.5) node {3};
\draw (4,-0.5) node {4};
\draw (5,-0.5) node {5};
%\draw (6,-0.5) node {6};
%\draw (7,-0.5) node {7};
%\draw (8,-0.5) node {8};
%\draw (9,-0.5) node {9};
%\draw (10,-0.5) node {10};
%\draw (11,-0.5) node {11};
%\draw (12,-0.5) node {12};

% node[pos=(0, -0,5)]{0};

%%caps,cups
\draw[very thick] [-,black,out=90, in=90](0,0.2) to (3,0.2);
\draw[very thick] [-,black,out=90, in=90](1,0.2) to (2,0.2);
%\draw[very thick] [-,black,out=90, in=90](11,0.2) to (12,0.2);
%\draw[very thick] [-,black,out=90, in=90](6,0.2) to (7,0.2);
%\draw[very thick] [-,black,out=90, in=90](5,0.2) to (8,0.2);
\draw[very thick] [-,black,out=90, in=90](0,+0.9) to (4,0.2);
\draw[very thick] [-,black,out=90, in=90](0,+0.9) to (5,0.2);

%\foreach \x in {} \draw + at (-1,0);

\end{tikzpicture} }
\smallskip

\text{The arc diagram of $-\times^2 \times \circ$}
\end{center}

\begin{center}
%\medskip
 
 \scalebox{0.7}{
\begin{tikzpicture}
 %\draw (-1,0) -- (7,0);
\foreach \x in {} %vee
     \draw[very thick] (\x-.1, .1) -- (\x,-0.1) -- (\x +.1, .1);
\foreach \x in {} %wedge
     \draw[very thick] (\x-.1, -.1) -- (\x,0.1) -- (\x +.1, -.1);
\foreach \x in {0,2} %cross
     \draw[very thick] (\x-.1, .1) -- (\x +.1, -.1) (\x-.1, -.1) -- (\x +.1, .1);
%\foreach \x in {1,3,4,7,8,9,10,12}  \draw[semithick] \circ; %circle
     %\draw[very thick]  node at (0,0) [fill=white,draw,circle,inner sep=0pt,minimum size=6pt]{};
     \draw[very thick]  node at (1,0) [fill=white,draw,circle,inner sep=0pt,minimum size=6pt]{};
     \draw[very thick]  node at (3,0) [fill=white,draw,circle,inner sep=0pt,minimum size=6pt]{};
     \draw[very thick]  node at (4,0) [fill=white,draw,circle,inner sep=0pt,minimum size=6pt]{};
     \draw[very thick]  node at (5,0) [fill=white,draw,circle,inner sep=0pt,minimum size=6pt]{};
     %\draw[very thick]  node at (8,0) [fill=white,draw,circle,inner sep=0pt,minimum size=6pt]{};
     %\draw[very thick]  node at (9,0) [fill=white,draw,circle,inner sep=0pt,minimum size=6pt]{};
     %\draw[very thick]  node at (10,0) [fill=white,draw,circle,inner sep=0pt,minimum size=6pt]{};
     %\draw[very thick]  node at (12,0) [fill=white,draw,circle,inner sep=0pt,minimum size=6pt]{};
\foreach \x in {0} %cross
     \draw[very thick] (\x-.1, +0.8) -- (\x +.1, +0.6) (\x-.1, +0.6) -- (\x +.1, +0.8);

\draw (0,-0.5) node {0};
\draw (1,-0.5) node {1};
\draw (2,-0.5) node {2};
\draw (3,-0.5) node {3};
\draw (4,-0.5) node {4};
\draw (5,-0.5) node {5};
%\draw (6,-0.5) node {6};
%\draw (7,-0.5) node {7};
%\draw (8,-0.5) node {8};
%\draw (9,-0.5) node {9};
%\draw (10,-0.5) node {10};
%\draw (11,-0.5) node {11};
%\draw (12,-0.5) node {12};

% node[pos=(0, -0,5)]{0};

%%caps,cups
\draw[very thick] [-,black,out=90, in=90](0,0.2) to (1,0.2);
\draw[very thick] [-,black,out=90, in=90](2,0.2) to (3,0.2);
%\draw[very thick] [-,black,out=90, in=90](11,0.2) to (12,0.2);
%\draw[very thick] [-,black,out=90, in=90](6,0.2) to (7,0.2);
%\draw[very thick] [-,black,out=90, in=90](5,0.2) to (8,0.2);
\draw[very thick] [-,black,out=90, in=90](0,+0.9) to (4,0.2);
\draw[very thick] [-,black,out=90, in=90](0,+0.9) to (5,0.2);

%\foreach \x in {} \draw + at (-1,0);

\end{tikzpicture} }
\smallskip

\text{The arc diagram of $-\times^2 \circ \times \circ$}
\end{center}

\begin{center}
%\medskip
 
 \scalebox{0.7}{
\begin{tikzpicture}
 %\draw (-1,0) -- (7,0);
\foreach \x in {} %vee
     \draw[very thick] (\x-.1, .1) -- (\x,-0.1) -- (\x +.1, .1);
\foreach \x in {} %wedge
     \draw[very thick] (\x-.1, -.1) -- (\x,0.1) -- (\x +.1, -.1);
\foreach \x in {0,3} %cross
     \draw[very thick] (\x-.1, .1) -- (\x +.1, -.1) (\x-.1, -.1) -- (\x +.1, .1);
%\foreach \x in {1,3,4,7,8,9,10,12}  \draw[semithick] \circ; %circle
     %\draw[very thick]  node at (0,0) [fill=white,draw,circle,inner sep=0pt,minimum size=6pt]{};
     \draw[very thick]  node at (2,0) [fill=white,draw,circle,inner sep=0pt,minimum size=6pt]{};
     \draw[very thick]  node at (1,0) [fill=white,draw,circle,inner sep=0pt,minimum size=6pt]{};
     \draw[very thick]  node at (4,0) [fill=white,draw,circle,inner sep=0pt,minimum size=6pt]{};
     \draw[very thick]  node at (5,0) [fill=white,draw,circle,inner sep=0pt,minimum size=6pt]{};
     %\draw[very thick]  node at (8,0) [fill=white,draw,circle,inner sep=0pt,minimum size=6pt]{};
     %\draw[very thick]  node at (9,0) [fill=white,draw,circle,inner sep=0pt,minimum size=6pt]{};
     %\draw[very thick]  node at (10,0) [fill=white,draw,circle,inner sep=0pt,minimum size=6pt]{};
     %\draw[very thick]  node at (12,0) [fill=white,draw,circle,inner sep=0pt,minimum size=6pt]{};
\foreach \x in {0} %cross
     \draw[very thick] (\x-.1, +0.8) -- (\x +.1, +0.6) (\x-.1, +0.6) -- (\x +.1, +0.8);

\draw (0,-0.5) node {0};
\draw (1,-0.5) node {1};
\draw (2,-0.5) node {2};
\draw (3,-0.5) node {3};
\draw (4,-0.5) node {4};
\draw (5,-0.5) node {5};
%\draw (6,-0.5) node {6};
%\draw (7,-0.5) node {7};
%\draw (8,-0.5) node {8};
%\draw (9,-0.5) node {9};
%\draw (10,-0.5) node {10};
%\draw (11,-0.5) node {11};
%\draw (12,-0.5) node {12};

% node[pos=(0, -0,5)]{0};

%%caps,cups
\draw[very thick] [-,black,out=90, in=90](0,0.2) to (1,0.2);
\draw[very thick] [-,black,out=90, in=90](3,0.2) to (4,0.2);
%\draw[very thick] [-,black,out=90, in=90](11,0.2) to (12,0.2);
%\draw[very thick] [-,black,out=90, in=90](6,0.2) to (7,0.2);
%\draw[very thick] [-,black,out=90, in=90](5,0.2) to (8,0.2);
\draw[very thick] [-,black,out=90, in=90](0,+0.9) to (2,0.2);
\draw[very thick] [-,black,out=90, in=90](0,+0.9) to (5,0.2);

%\foreach \x in {} \draw + at (-1,0);

\end{tikzpicture} }
\medskip

\text{The arc diagram of $-\times^2\circ \circ \times \circ$}
\end{center}

\begin{center}
%\medskip
 
 \scalebox{0.7}{
\begin{tikzpicture}
 %\draw (-1,0) -- (7,0);
\foreach \x in {} %vee
     \draw[very thick] (\x-.1, .1) -- (\x,-0.1) -- (\x +.1, .1);
\foreach \x in {} %wedge
     \draw[very thick] (\x-.1, -.1) -- (\x,0.1) -- (\x +.1, -.1);
\foreach \x in {0,4} %cross
     \draw[very thick] (\x-.1, .1) -- (\x +.1, -.1) (\x-.1, -.1) -- (\x +.1, .1);
%\foreach \x in {1,3,4,7,8,9,10,12}  \draw[semithick] \circ; %circle
     %\draw[very thick]  node at (0,0) [fill=white,draw,circle,inner sep=0pt,minimum size=6pt]{};
     \draw[very thick]  node at (2,0) [fill=white,draw,circle,inner sep=0pt,minimum size=6pt]{};
     \draw[very thick]  node at (3,0) [fill=white,draw,circle,inner sep=0pt,minimum size=6pt]{};
     \draw[very thick]  node at (1,0) [fill=white,draw,circle,inner sep=0pt,minimum size=6pt]{};
     \draw[very thick]  node at (5,0) [fill=white,draw,circle,inner sep=0pt,minimum size=6pt]{};
     %\draw[very thick]  node at (8,0) [fill=white,draw,circle,inner sep=0pt,minimum size=6pt]{};
     %\draw[very thick]  node at (9,0) [fill=white,draw,circle,inner sep=0pt,minimum size=6pt]{};
     %\draw[very thick]  node at (10,0) [fill=white,draw,circle,inner sep=0pt,minimum size=6pt]{};
     %\draw[very thick]  node at (12,0) [fill=white,draw,circle,inner sep=0pt,minimum size=6pt]{};
\foreach \x in {0} %cross
     \draw[very thick] (\x-.1, +0.8) -- (\x +.1, +0.6) (\x-.1, +0.6) -- (\x +.1, +0.8);

\draw (0,-0.5) node {0};
\draw (1,-0.5) node {1};
\draw (2,-0.5) node {2};
\draw (3,-0.5) node {3};
\draw (4,-0.5) node {4};
\draw (5,-0.5) node {5};
%\draw (6,-0.5) node {6};
%\draw (7,-0.5) node {7};
%\draw (8,-0.5) node {8};
%\draw (9,-0.5) node {9};
%\draw (10,-0.5) node {10};
%\draw (11,-0.5) node {11};
%\draw (12,-0.5) node {12};

% node[pos=(0, -0,5)]{0};

%%caps,cups
\draw[very thick] [-,black,out=90, in=90](0,0.2) to (1,0.2);
\draw[very thick] [-,black,out=90, in=90](4,0.2) to (5,0.2);
%\draw[very thick] [-,black,out=90, in=90](11,0.2) to (12,0.2);
%\draw[very thick] [-,black,out=90, in=90](6,0.2) to (7,0.2);
%\draw[very thick] [-,black,out=90, in=90](5,0.2) to (8,0.2);
\draw[very thick] [-,black,out=90, in=90](0,+0.9) to (2,0.2);
\draw[very thick] [-,black,out=90, in=90](0,+0.9) to (3,0.2);

%\foreach \x in {} \draw + at (-1,0);

\end{tikzpicture} }
\smallskip

\text{The arc diagram of $-\times^2 \circ \circ\ \circ \times \circ$}
\end{center}

For $\mathfrak{osp}(6|4)$ the arc diagram of $>\times\circ \circ \times$ is

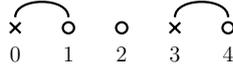
\begin{center}
%\medskip
 
 \scalebox{0.7}{
\begin{tikzpicture}
 %\draw (-1,0) -- (7,0);
\foreach \x in {} %vee
     \draw[very thick] (\x-.1, .1) -- (\x,-0.1) -- (\x +.1, .1);
\foreach \x in {} %wedge
     \draw[very thick] (\x-.1, -.1) -- (\x,0.1) -- (\x +.1, -.1);
\foreach \x in {0,3} %cross
     \draw[very thick] (\x-.1, .1) -- (\x +.1, -.1) (\x-.1, -.1) -- (\x +.1, .1);
%\foreach \x in {1,3,4,7,8,9,10,12}  \draw[semithick] \circ; %circle
     %\draw[very thick]  node at (0,0) [fill=white,draw,circle,inner sep=0pt,minimum size=6pt]{};
     \draw[very thick]  node at (2,0) [fill=white,draw,circle,inner sep=0pt,minimum size=6pt]{};
     \draw[very thick]  node at (1,0) [fill=white,draw,circle,inner sep=0pt,minimum size=6pt]{};
     \draw[very thick]  node at (4,0) [fill=white,draw,circle,inner sep=0pt,minimum size=6pt]{};
     %\draw[very thick]  node at (5,0) [fill=white,draw,circle,inner sep=0pt,minimum size=6pt]{};
     %\draw[very thick]  node at (8,0) [fill=white,draw,circle,inner sep=0pt,minimum size=6pt]{};
     %\draw[very thick]  node at (9,0) [fill=white,draw,circle,inner sep=0pt,minimum size=6pt]{};
     %\draw[very thick]  node at (10,0) [fill=white,draw,circle,inner sep=0pt,minimum size=6pt]{};
     %\draw[very thick]  node at (12,0) [fill=white,draw,circle,inner sep=0pt,minimum size=6pt]{};
%\foreach \x in {0} %cross
 %    \draw[very thick] (\x-.1, +0.8) -- (\x +.1, +0.6) (\x-.1, +0.6) -- (\x +.1, +0.8);

\draw (0,-0.5) node {0};
\draw (1,-0.5) node {1};
\draw (2,-0.5) node {2};
\draw (3,-0.5) node {3};
\draw (4,-0.5) node {4};
%\draw (5,-0.5) node {5};
%\draw (6,-0.5) node {6};
%\draw (7,-0.5) node {7};
%\draw (8,-0.5) node {8};
%\draw (9,-0.5) node {9};
%\draw (10,-0.5) node {10};
%\draw (11,-0.5) node {11};
%\draw (12,-0.5) node {12};

% node[pos=(0, -0,5)]{0};

%%caps,cups
\draw[very thick] [-,black,out=90, in=90](0,0.2) to (1,0.2);
\draw[very thick] [-,black,out=90, in=90](3,0.2) to (4,0.2);
%\draw[very thick] [-,black,out=90, in=90](11,0.2) to (12,0.2);
%\draw[very thick] [-,black,out=90, in=90](6,0.2) to (7,0.2);
%\draw[very thick] [-,black,out=90, in=90](5,0.2) to (8,0.2);
%\draw[very thick] [-,black,out=90, in=90](0,+0.9) to (4,0.2);
%\draw[very thick] [-,black,out=90, in=90](0,+0.9) to (5,0.2);

%\foreach \x in {} \draw + at (-1,0);

\end{tikzpicture} }

\end{center}

and the arc diagram of $+\times\circ \circ \times$ is

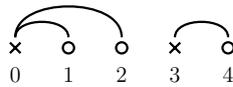
\begin{center}
%\medskip
 
 \scalebox{0.7}{
\begin{tikzpicture}
 %\draw (-1,0) -- (7,0);
\foreach \x in {} %vee
     \draw[very thick] (\x-.1, .1) -- (\x,-0.1) -- (\x +.1, .1);
\foreach \x in {} %wedge
     \draw[very thick] (\x-.1, -.1) -- (\x,0.1) -- (\x +.1, -.1);
\foreach \x in {0,3} %cross
     \draw[very thick] (\x-.1, .1) -- (\x +.1, -.1) (\x-.1, -.1) -- (\x +.1, .1);
%\foreach \x in {1,3,4,7,8,9,10,12}  \draw[semithick] \circ; %circle
     %\draw[very thick]  node at (0,0) [fill=white,draw,circle,inner sep=0pt,minimum size=6pt]{};
     \draw[very thick]  node at (2,0) [fill=white,draw,circle,inner sep=0pt,minimum size=6pt]{};
     \draw[very thick]  node at (1,0) [fill=white,draw,circle,inner sep=0pt,minimum size=6pt]{};
     \draw[very thick]  node at (4,0) [fill=white,draw,circle,inner sep=0pt,minimum size=6pt]{};
     %\draw[very thick]  node at (5,0) [fill=white,draw,circle,inner sep=0pt,minimum size=6pt]{};
     %\draw[very thick]  node at (8,0) [fill=white,draw,circle,inner sep=0pt,minimum size=6pt]{};
     %\draw[very thick]  node at (9,0) [fill=white,draw,circle,inner sep=0pt,minimum size=6pt]{};
     %\draw[very thick]  node at (10,0) [fill=white,draw,circle,inner sep=0pt,minimum size=6pt]{};
     %\draw[very thick]  node at (12,0) [fill=white,draw,circle,inner sep=0pt,minimum size=6pt]{};
%\foreach \x in {0} %cross
 %    \draw[very thick] (\x-.1, +0.8) -- (\x +.1, +0.6) (\x-.1, +0.6) -- (\x +.1, +0.8);

\draw (0,-0.5) node {0};
\draw (1,-0.5) node {1};
\draw (2,-0.5) node {2};
\draw (3,-0.5) node {3};
\draw (4,-0.5) node {4};
%\draw (5,-0.5) node {5};
%\draw (6,-0.5) node {6};
%\draw (7,-0.5) node {7};
%\draw (8,-0.5) node {8};
%\draw (9,-0.5) node {9};
%\draw (10,-0.5) node {10};
%\draw (11,-0.5) node {11};
%\draw (12,-0.5) node {12};

% node[pos=(0, -0,5)]{0};

%%caps,cups
\draw[very thick] [-,black,out=90, in=90](0,0.2) to (1,0.2);
\draw[very thick] [-,black,out=90, in=90](0,0.2) to (2,0.2);
\draw[very thick] [-,black,out=90, in=90](3,0.2) to (4,0.2);
%\draw[very thick] [-,black,out=90, in=90](6,0.2) to (7,0.2);
%\draw[very thick] [-,black,out=90, in=90](5,0.2) to (8,0.2);
%\draw[very thick] [-,black,out=90, in=90](0,+0.9) to (4,0.2);
%\draw[very thick] [-,black,out=90, in=90](0,+0.9) to (5,0.2);

%\foreach \x in {} \draw + at (-1,0);

\end{tikzpicture} }

\end{center}

\subsubsection{Description}
The arc diagram $Arc(f)$ is unique. For instance, for $\ell=0$
the arc diagram $Arc(f)$ can be constructed as follows:

first we consider the diagram $f'$ obtained from $f$ by
removing all but one symbol $\times$ from the zero position
($f'=f$ if the zero position contains $\circ$ or $\times$);

for $f'$ we construct the arc diagram in the usual way (we
connect each symbol $\times$ with the next free empty symbol
starting from the rightmost symbol $\times$ and going to the left);

if the zero in $f$ is occupied by $\times^p$ with $p>1$,
we connect each of the remaining $p-1$ symbol $\times$ with next
two free positions starting from the lowest symbol $\times$ and going up.

\subsubsection{Maximal arcs}\label{maxarc}
Note that an arc $arc(0;b)$ (resp., $(0;b_1,b_2)$)
is maximal if and only if it is supported by the top symbol $\times$
in the zero position.  Note that  $Arc(\times^{p}f)$ is obtained from $Arc(\times^{p+1} f)$
by "removal" of the maximal arc supported by the top symbol $\times$
in the zero position. Notice that
for $\osp(2n+1|2n)$, when we
erase a symbol $\times$ from a signed diagram we either obtain
a signed diagram with the same sign or an unsigned diagram  if the resulting diagram
does not have $\times$ in the zero position.

Consider a diagram of the form $f_1\times f_2$ where the symbol
$\times$ occupies a  position $a>0$.
This symbol $\times$ supports a maximal arc $arc(a;b)$ if and only if
$Arc(f_1\circ f_2)$ is obtained from
$Arc(f_1\times f_2)$ by "removal" of the arc $arc(a;0)$, i.e.
$$Arc(f_1\circ f_2)=Arc(f_1)\circ Arc(f_2)$$
(if we remove an arc which is not maximal, the resulting diagram is not an arc diagram).

An arc $arc(a;b)$ with  $a>0$ is maximal in $Arc(f)$ if and only if for each $0<c<a$ one has $\ell_f(c;a)\geq 0$, where
 $\ell_f(c;a)$ is the number of the symbols $\times$
minus the number of empty positions strictly between the positions $c$ and $a$.

For instance, in~\Exa{exaMitka}, the maximal arcs are $arc(11; 12)$
and $arc(0;4,9)$.

\subsection{}
\begin{thm}{thmDS1osp}
(i) The module $L(\nu)$ is a subquotient of
$\DS_1(L(\lambda))$ if and only if
$Arc(\howl(\nu))$ is obtained from $Arc(\howl(\lambda))$ by removing
a maximal arc and, in addition,
in the $\osp(2m+1|2n)$-case,
if $\nu$ has a sign, then the signs of $\lambda$ and $\nu$ are equal.

(ii) Let $e$ be the number of free positions in $Arc(\howl(\lambda))$ which are to the left of the maximal arc.

For $\osp(2m+1|2n)$
one has
\begin{equation}\label{eee}
[\DS_1(L(\lambda)):L(\nu)]=\left\{\begin{array}{ll}
(1|0) & \text{ if }e=0;\\
(2|0) & \text{ if $e$ is even and }e\not=0;\\
(0|2)  & \text{ if $e$ is odd}.\end{array}
\right.\end{equation}

The same formula holds for
the case $t=2$. For  $t=0$
\begin{equation}\label{eee-2}
[\DS_1(L(\lambda)):L(\nu)]=\left\{\begin{array}{ll}
(1|0) & \text{ if $e$ is even};\\
(0|1)  & \text{ if $e$ is odd}.\end{array}
\right.\end{equation}

\end{thm}
%
%\subsubsection{Remark}
%For $\osp(2m|2n)$-case the formula (\ref{eee}) gives
%$$[\DS_1(L(\lambda)):L(\nu)]=[\DS_1(L(\lambda)):L(\nu^{\sigma})].$$
%
%
%

\subsubsection{Examples}
Take $r=1$ and use the notations of Section \ref{sect4}.

For $\osp(5|4)$ we have
$$\frac{+\times\circ \circ \times}{\circ\circ \circ\times}=(1|0)\ \ \ \ \ \ \ \
\frac{+\times\circ \circ \times}{+\times}=(0|2)$$
and $\frac{+\times\circ \circ \times}{f'}=0$ in other cases. The multiplicity is
$$\frac{-\times^2\underbrace{\circ...\circ}_{j\text{ times}} \times}{-\times^2}=\left\{\begin{array}{ll}
0 & \text{ for } j<3,\\
(1|0) & \text{ for } j=3,\\
(2|0) & \text{ if $j>3$ is odd},\\
(0|2) & \text{ if $j>3$ is even}.\\
\end{array}\right.$$

For $\osp(4|4)$ we have
$$\frac{+\circ\times\circ\times}{+\circ\times}=\frac{+\circ\times\circ\times}{-\circ\times}=(0|1).$$

For $\osp(6|4)$ we have
$$\frac{>\times\circ \circ \times}{>\circ\circ \circ\times}=(1|0)\ \ \ \ \ \ \ \
\frac{+\times\circ \circ \times}{+\times}=(0|2)$$
and $\frac{+\times\circ \circ \times}{f'}=0$ in other cases.

\subsection{Low rank cases} \label{low-rank} The proof will be a reduction to the atypicality 1 case, so consider the cases $\osp(2+t|2)$ for $t=0,1,2$. In these cases $DS(L(\lambda))$ was computed in \cite{Gdex} (note that $x$ is necessarily of rank 1 and we omit it from the notation). For $t=0,1$ $\mathfrak{g}_x = 0$, so $\mathfrak{g}_x$-modules are supervector spaces. For $t=2$ $\mathfrak{g}_x = \mathbb{C}$. For the principal block $\cF^>(\mathfrak{g})$ $DS_x(\cF^>(\mathfrak{g}))$ is the category of finite-dimensional supervector spaces with trivial action of $\mathfrak{g}_x$. The simple modules in the principal block are of the form $L(\lambda_j)$, $j \in \mathbb{Z}$, for $t=0$ and $j \in \mathbb{N}$ for $t=1,2$.

For $t=0$ \[ DS(L(\lambda_j)) \cong \Pi^j (\mathbb{C}) \ \ \forall j \in \mathbb{Z} \] where $\lambda_j := j \varepsilon_1 + |j|\delta_1$ for $j \in \mathbb{Z}$.

For $t=1$ we put $\lambda_0 := 0$ and $\lambda_j:= j\varepsilon_1 + (j-1)\delta_1$ for $j \geq 1$. For $t=2$ we put $\lambda_j := j\varepsilon_1 + j \delta_1$ for $j \geq 0$. Then one has the uniform rule \cite{Gdex} \begin{align*} DS (L(\lambda_0) = DS(L(\lambda_1)) & \cong \mathbb{C} \text{ and } \\ DS(L(\lambda_j)) & \cong \Pi^{j-1} (\mathbb{C})^{\oplus 2} \text{ for } j\geq 2.\end{align*}

\subsection{Proof}
Set $f:=\howl(\nu)$. Denote by
$u$  the coordinate of the rightmost symbol $\times$ in $f$.

If $u>0$, then $Arc(f)$ has a minimal arc $arc(u;u+1)$.
If $\frac{f'}{f}\not=0$, then, by \ref{mx}, $f'$ has a minimal arc $arc(u;u+1)$ and
$$\frac{f'}{f}=\frac{\ol{f}'}{\ol{f}},$$
where $\ol{f}'$ (resp., $\ol{f}$) are obtained from $f'$ (resp., $f$)
by ``shrinking'' the minimal arc $arc(u;u+1)$.
Clearly, $Arc(f)$ can be obtained from $Arc(f')$
by deleting a maximal arc if and only if
 $Arc(\ol{f})$ can be obtained from $Arc(\ol{f}')$ by deleting a maximal arc. Each maximal arc
in $Arc(\ol{f})$ corresponds to a maximal arc in $Arc(f)$ and
the number of free positions to the left of each maximal arc is the same.
This reduces the statement to the case $u=0$. Corollaries \ref{cor2m12n},\ref{cor2m2n}  reduce the statement to the
cases when $f=\emptyset$ or $f=>$. Since $r=1$, the case $f=\emptyset$
correspond to $\osp(2|2),\osp(3|2)$ and the case $f=>$
correspond to the case $\osp(4|2)$. For
these cases the formula was checked in \cite{Gdex}(see \ref{low-rank}).
\qed

We will use the following lemma.
\subsubsection{}
\begin{lem}{lemtau}
Take $\tau$ as in~\ref{t12}. For $r=1$ one has
$\frac{f}{f'}=\frac{\tau(f)}{\tau(f')}$.
\end{lem}
\begin{proof}
One has to check the following cases
$$\begin{array}{lcl}
f=\overset{\ \ \times^{p}}{>}\circ f_1\times f_2 & & f'=\overset{\ \ \times^{p}}{>}\circ f_1\circ f_2\\
f=\overset{\ \ \times^{p}}{>}\times f_1\times f_2 & &  f'=\overset{\ \ \times^{p}}{>}\times f_1\circ f_2\\
f=\overset{\ \ \times^{p}}{>}\times f_1 & & f'=\overset{\ \ \times^{p}}{>}\circ f_1\\
f=\overset{\ \ \times^{p+1}}{>}\circ f_1 & &  f'=\overset{\ \ \times^{p}}{>}\circ f_1\\
f=\overset{\ \ \times^{p+1}}{>}\times f_1 & & f'=\overset{\ \ \times^{p}}{>}\times f_1.
\end{array}$$
This can be easily  done using the properties
of maximal arcs discussed in \ref{maxarc}.
\end{proof}

\section{Semisimplicity of $\DS_x(L)$} \label{sec:Semi}
Retain the notation of \ref{mapdex}. Denote by $\cF_+(\fg)$ the Serre subcategory 
of $\tilde{\cF}(\fg)$ generated
by $L\in\Irr\bigr(\tilde{\cF}(\fg)\bigr)$ with $\pari(L)=1$.
By~\Cor{corext}, $\cF_+(\fg)$ is semisimple. 

In~\Thm{thmssDS}
we will show that $\DS_x(\cF_+(\fg))\subset\cF_+(\fg_x)$. As a result,
$\DS_x(L)$ is semisimple
for each simple finite-dimensional module $L$.

\subsection{Supercharacters} Recall that besides the usual character we have the supercharacter $$\sch L(\lambda)=(-1)^{p(\lambda)}\pi (\ch L(\lambda)),$$
where $\pi:\mathbb{Z}[\Lambda_{m|n}]\to \mathbb{Z}[\Lambda_{m|n}]$ is the linear involution
given by
$$\pi(e^{\mu}):=(-1)^{p(\mu)}e^{\mu}.$$ 
The supercharacter ring of $\Fin$ is the image of the map
$\sch: \Fin\to \mathbb{Z}[\fh^*]$; we denote this ring by $\mathcal{J}(\fg)$.

%The algebra
%$\fh_x:=\fh^{\ad x}/([x,\fg]\cap\fh)$ is a Cartan subalgebra of $\fg_x$.
 %By \cite{HR}, Lemma 4 for  each $N\in\tilde{\cF}(\mathfrak{g})$ we have
%$$\sch \DS_x(N)=\ds_x(\sch N),$$
 %where $\ds_x:\mathcal{J}(\fg)\to \mathcal{J}(\fg_x)$ is a ring homomorphism
%given by $f\mapsto f|_{\fh_x}$.

\subsection{}
We will use the following lemma.

\begin{lem}{lemF+}
Let $M\in\Fin(\fg)$ and $N\in\cF_+(\fg)$ be such that $\sch M=\sch N$ in the supercharacter ring. If  $\dim M\leq \dim N$, then
$M\cong N$.
\end{lem}

\begin{proof}
For each $\nu\in\Lambda^{(t)}_{m+\ell|n}$ we set 
$$(d_0(\nu)| d_1(\nu)):=[N:L(\nu)],\ \ (d'_0(\nu)| d'_1(\nu)):=[M:L(\nu)].$$
Since $N\in \cF_+(\fg)$  we have
$d_0(\nu)d_1(\nu)=0$.

Since $\{\sch L(\nu)|\ \nu\in\Lambda^{(t)}_{m+\ell|n}\}$ are linearly independent, 
the equality $\sch M=\sch N$
implies
$$(d'_0(\nu)| d'_1(\nu))=(d_0(\nu)+j(\nu)| d_1(\nu)+j(\nu))\ \ \text{ 
for some }j(\nu)\in\mathbb{Z}.$$
 Combining $d_0(\nu)d_1(\nu)=0$
with $d_0'(\nu),d_1'(\nu)\geq 0$, we obtain $j(\nu)\geq 0$
for each $\nu$.  Using $\dim M\leq \dim N$ we get
$j(\nu)=0$ for each $\nu$, that is
\begin{equation}\label{MNLnu}
\forall \nu\ \ \ [M:L(\nu)]=[N:L(\nu)].\end{equation}
Hence $M\in \cF_+(\fg)$. Since $\cF_+(\fg)$ is semisimple, $M$ and $N$
are completely reducible. Thus (\ref{MNLnu}) gives $M\cong N$.
\end{proof}

\subsection{}
\begin{thm}{thmssDS}
\begin{enumerate}
\item For each $x$ one has  $\DS_x\bigl(\cF_+(\fg)\bigr)=\cF_+(\fg_x)$.

\item
For each  $L\in\Irr(\tilde{\cF}(\fg))$ one has  $\DS_{r+1}(L)\cong \DS_{1}(\DS_{r}(L))$.
\end{enumerate}
\end{thm}
\begin{proof}
It is enough to consider the case $x:=x_r$.

For (i) we have to verify that for   each  $L(\lambda)\in\Irr(\cF_+(\fg))$,
and $L_{\fg_x}(\nu)\in \Irr\bigr({\cF}(\fg_x)\bigr)$ the graded multiplicity 
$$(d_0|d_1):= [\DS_x(L(\lambda)):L_{\fg_x}(\nu)]$$
satisfies
\begin{equation}\label{eqd0d1}
\left\{
\begin{array}{ll}
d_1=0\ & \text{ if }\ \pari(\nu)=1,\\
d_0=0 \ & \text{ if }\ \pari(\nu)=-1.
\end{array}
\right.\end{equation}
 Recall that 
$\pari(\lambda)=\pari(\howl(\lambda))$. Using~\Cor{corka} we reduce
(i) to the case when $\lambda,\nu$ are core-free. In this case
$\fg=\osp(2n+t|2n)=\fg_n$ (for some $n$) and $\fg_x=\fg_{n-r}$.

We proceed by induction on $r$.  Note that (ii) for $r=0$ is  tautological. 

Consider the case $r=1$ for $t=0,1$.
For  a core-free diagram $f$ 
denote by $||f||$  the sum of the coordinates of the symbols
$\times$ in $f$. One has $\pari(\lambda(f))=(-1)^{||f||}$, so~(\ref{eqd0d1})
follows from~\Thm{thmDS1osp}.
The case $r=1$ for $t=2$ follows from $t=1$ and~\Lem{lemtau}.
 This establishes (i) for $r=1$.

Now fix any $t$ and take
$r\geq 2$. By induction,
$\DS_{r-1}(L(\lambda))\in \cF_+(\fg_{n-r+1})$.
Using (i) for $r=1$ we get
$$N:=\DS_{1}(\DS_{r-1}(L(\lambda))\in  \cF_+(\fg_x).$$
By~\cite{HR}, 
$\sch N=\sch\DS_{r}(L(\lambda))$;
by~\cite{Gaugusta}, Lem. 2.4.1,  $\dim \DS_{r}(L(\lambda))\leq \dim N$.
Using~\Lem{lemF+} we obtain $N\cong \DS_{r}(L(\lambda))$
as required.
\end{proof}

\subsection{}
\begin{cor}{cortau}
Take $\tau:\Lambda^{(2)}_{m+1|n}\to\Lambda^{(1)}_{m|n}$ as in~\ref{t12}. One has
$$[\DS_r(L(\lambda)):L(\fg_x(\nu))]=[\DS_r(L(\tau(\lambda))):L(\fg_x(\tau(\nu)))].$$
\end{cor}
\begin{proof}
The case $r=1$ was treated in~\Lem{lemtau}. The general case follows from~\Thm{thmssDS} (ii). 
\end{proof}

\subsubsection{} We say that a module $M$ is \emph{pure} if for any subquotient $L$ of $M$, $\Pi(L)$ is not a subquotient of $M$. Theorem \ref{thmssDS} implies immediately the following assertion.

\subsubsection{}
\begin{cor}{pure} $DS_x(L(\lambda))$ is pure for any $x$ and any $\lambda$.
\end{cor}

%\Thorsten{Need to say a little bit about $OSp$}

\subsubsection{}
\begin{cor}{corDSsemidecomp} For irreducible $L(\lambda)$ \[ \DS_1(L(\lambda)) \cong \bigoplus_i m_i \Pi^{n_i} L(\lambda_i) \] where the arc diagram of $L(\lambda_i)$ is obtained by removing the $i$-th maximal arc and the associated $\times$ from the arc diagram of $\lambda$. The multiplicity $m_i$ is 1 or 2 according to the rules of Theorem \ref{thmDS1osp} and $n_i =1 \ mod \ 2$ if and only if the parities of $\lambda$ and $\lambda_i$ differ.
\end{cor}

Since $\DS_{r+1}(L)\cong \DS_{1}(\DS_{r}(L))$ we can calculate any $\DS_{r+1}(L)$ by repeated application of $DS_1$.

\subsubsection{Remark} In the $\mathfrak{gl}(m|n)$-case $\DS(L(\lambda))$ is even multiplicity free \cite{HW}.

\subsection{The $OSp$-case}
Using~\Cor{corDSsemidecomp} and~\Lem{lemsigma} it is easy to
describe the effect of $DS_1$ on irreducible $OSp(M|N)$-modules $L_{OSp}(\lambda)$, see below.

Let $\tilde{\cF}'(M|N)$ denote the category of algebraic representations of $OSp(M|N)$, then we have a commutative diagram \[ \xymatrix{ \tilde{\cF}'(M|N) \ar[d]^{DS_1} \ar[r]^{Res} & \tilde{\cF}(M|N) \ar[d]^{DS_1} \\ \tilde{\cF}'(M-2|N-2) \ar[r]^{Res} & \tilde{\cF}(M-2|N-2).} \] 

Formula (\ref{eee}) holds for $OSp(2m+1|2n)$-modules if the signs of $L_{OSp}(\lambda,\pm)$ and $L_{OSp}(\nu,\pm)$ are equal; otherwise the multiplicity is zero. 

Consider $OSp(2m|2n)$-case.  
Combining~\Rem{OSP2m2n}, \Cor{corDSsemidecomp} and~\Lem{lemsigma} we conclude that
$\DS_x(L(\lambda))$ has a structure of $OSp(2m|2n)$-module. Let us show that the
multiplicity $d_{OSp}(\lambda;\nu):=[\DS_1(L_{OSp}(\lambda)):L_{OSp}(\nu)]$ is given by the formula~(\ref{eee-2})
in the case $\lambda^{\sigma}=\lambda$, $t=0$ and by the formula~(\ref{eee}) otherwise.

Indeed,  combining~\Rem{OSP2m2n} and~\Lem{lemsigma} we get
$$d_{OSp}(\lambda;\nu)=[\DS_1(L_{OSp}(\lambda)):L(\nu)]=\left\{\begin{array}{ll} [\DS_1(L(\lambda)):L(\nu)]
 & \text{ if $\lambda^{\sigma}=\lambda$}\\
2[\DS_1(L(\lambda)):L(\nu)] & \text{ otherwise}.\end{array}\right.$$
Thus $d_{OSp}(\lambda;\nu)$ is given
by  the formula~(\ref{eee-2})  (resp.,~(\ref{eee}))  for the case
when $\lambda^{\sigma}=\lambda$ and $t=0$ (resp., $t=2$).
For the remaining case $t=0$ and $\lambda^{\sigma}\not=\lambda$ we get
$$d_{OSp}(\lambda;\nu)=\left\{\begin{array}{ll} (2|0) & \text{ if $e$ is even}\\
(0|2) & \text{ if $e$ is odd},\end{array}\right.
$$
where $e$  as in~\Thm{thmDS1osp}.
Notice that $e\not=0$, since the condition $\lambda^{\sigma}\not=\lambda$
implies that the zero position of the diagram of $\lambda$ is empty.
Hence $d_{OSp}(\lambda;\nu)$ is given
by  the formula~(\ref{eee}).

\subsection{Question} For a simple module $L$ in $\tilde{\cF}(\osp(2m|2n))$, 
 \Lem{lemsigma} and~\Thm{thmssDS} imply that $\DS_x(L)$ is an $OSp(2m|2n)$-module for any $x$. Is this still true for an arbitrary module in $\tilde{\cF}(2m|2n)$?
 
A related question can be asked about 
finite-dimensional modules over $\fg=D(2|1,a), F(4)$ (for these cases $\DS_x(\fg)$
admits an involution $\sigma$ and $\DS_x(L)$ is $\sigma$-invariant
for each simple module $L$ in $\tilde{\cF}(\fg)$, see~\cite{Gdex}).

\subsection{Superdimensions} Similarly to \cite{HW} in the $\mathfrak{gl}(m|n)$-case this allows us now to compute the superdimension of any irreducible $L(\lambda)$. Let $\lambda$ be maximal atypical and $x$ of rank equal to the atypicality. Then $\mathfrak{g}_x$ is either an orthogonal or symplectic Lie algebra or $\osp(1|2r)$ for some $r$. In each case the superdimension of an irreducible module is known. Since $DS$ is a symmetric monoidal functor it preserves the superdimension. So $sdim(L(\lambda)) = sdim(DS(L(\lambda)))$. So in order to compute $\sdim(L(\lambda))$ it suffices to compute the multiplicity $m(\lambda)$ of the isotypic representation $DS(L(\lambda))$ of $\mathfrak{g}_x$, but this multiplicity is computed exactly by Theorem \ref{thmDS1osp} since $DS_r (L(\lambda)) = DS_1(\ldots (DS_1(L(\lambda))))$. Alternatively the multiplicity can also be computed via~\cite{GH}.

\subsubsection{Example} 
Consider a weight $\lambda$ with the diagram consisting of
$m$ maximal arcs (with no other arcs) and $i$ symbols $\times$ in the zero position
(since all arcs are maximal $i\leq 1$). Using the results of~\cite{HW} for $\fgl$-case we get
$$\sdim L_{\fg}(\lambda)=\left\{\begin{array}{ll}
m!  & \text{ for }\fg=\fgl(m|m)\\
2^{m-1} m!  & \text{ for }\fg=\fosp(2m|2m)\\
2^{m-i} m!  & \text{ for }\fg=\fosp(2m+1|2m), \fosp(2m+2|2m).\\
\end{array}
\right.$$

\subsubsection{Example} 
For an
$\osp(2m+1|2n)$-diagram with the empty zero position the multiplicity is $2^l |F|!/F!$ for $l = min(m|n)$ where $F$ is the forest associated to the arc diagram of $\lambda$ (see \cite{HW}). For $t=0$ the multiplicity $m(\lambda)$ is given by a forest factorial exactly as for $\mathfrak{gl}(n|n)$, but one has to take into account
that in this case
a removal of one arc may produce  two arc diagrams, which differ by the sign
(if the  resulting arc diagram has  $\times$ in the zero position).

\appendix

\section{The language of Gruson-Serganova and Ehrig-Stroppel} Weight and arc diagrams also appear in the work of Ehrig and Stroppel \cite{ES}, but our conventions differ considerably. They attach to a weight $\lambda$ a weight diagram $\pla^{\owedge}$ and a cup diagram $\underline{\pla^{\owedge}}$. We shortly explain the differences, mostly following \cite{ES}, Section 6.3. 

\subsection{The associated tailless weight} Gruson-Serganova \cite{GSBGG}, Section 8, associate to every weight $\lambda$ a tailless weight $\bar{\lambda}$. 

Assume that $M = 2m+1$ and consider the core-free weight diagram $f_{\lambda}$. Remove the tail of $f_{\lambda}$. Note that in case the diagram has a sign $+$, one symbol $\times$ remains at the zero position. Form the arc diagram of the new weight diagram (with removed tail) and number the free vertices $\circ$ from left to right by $1,2,3,\ldots$. Label the vertices at positions $1,3,\ldots$ then by an $\times$ until $2l-1$ where $l$ is the number of removed tail symbols $\times$ (these $\times$ are called coloured in \cite{GSBGG}). This gives the weight diagram $f_{\bar{\lambda}}$. To assign an arc diagram, first form the arc diagram of the weight diagram with removed tail and then connect the coloured $\times$' in the usual way with free positions. Finally place a dot on each arc which starts at one of the edges $1,3,\ldots,2l-1$. For $M = 2m$ even, proceed in the same way as if $f_{\lambda}$ was a weight diagram for $M = 2m+1$ with a $+$-sign. 

\subsection{Example} If $f_{\lambda} = + \times^3 \times \circ \circ \ldots$, then removing the tail gives the diagram $\times \times \circ \circ \ldots$. The weight diagram $f_{\bar{\lambda}}$ is given by $+ \times \times \circ \circ \times \circ \times \circ \circ \ldots$. Its dotted arc diagram is given by 

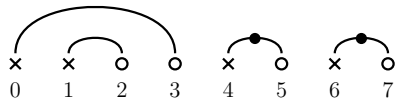
\begin{center}
%\medskip
 
 \scalebox{0.7}{
\begin{tikzpicture}
 %\draw (-1,0) -- (7,0);
\foreach \x in {} %vee
     \draw[very thick] (\x-.1, .1) -- (\x,-0.1) -- (\x +.1, .1);
\foreach \x in {} %wedge
     \draw[very thick] (\x-.1, -.1) -- (\x,0.1) -- (\x +.1, -.1);
\foreach \x in {0,1,4,6} %cross
     \draw[very thick] (\x-.1, .1) -- (\x +.1, -.1) (\x-.1, -.1) -- (\x +.1, .1);
%\foreach \x in {1,3,4,7,8,9,10,12}  \draw[semithick] \circ; %circle
     %\draw[very thick]  node at (0,0) [fill=white,draw,circle,inner sep=0pt,minimum size=6pt]{};
     \draw[very thick]  node at (2,0) [fill=white,draw,circle,inner sep=0pt,minimum size=6pt]{};
     \draw[very thick]  node at (3,0) [fill=white,draw,circle,inner sep=0pt,minimum size=6pt]{};
     \draw[very thick]  node at (7,0) [fill=white,draw,circle,inner sep=0pt,minimum size=6pt]{};
     \draw[very thick]  node at (5,0) [fill=white,draw,circle,inner sep=0pt,minimum size=6pt]{};
     %\draw[very thick]  node at (8,0) [fill=white,draw,circle,inner sep=0pt,minimum size=6pt]{};
     %\draw[very thick]  node at (9,0) [fill=white,draw,circle,inner sep=0pt,minimum size=6pt]{};
     %\draw[very thick]  node at (10,0) [fill=white,draw,circle,inner sep=0pt,minimum size=6pt]{};
     %\draw[very thick]  node at (12,0) [fill=white,draw,circle,inner sep=0pt,minimum size=6pt]{};
%\foreach \x in {0} %cross
 %    \draw[very thick] (\x-.1, +0.8) -- (\x +.1, +0.6) (\x-.1, +0.6) -- (\x +.1, +0.8);

\draw (0,-0.5) node {0};
\draw (1,-0.5) node {1};
\draw (2,-0.5) node {2};
\draw (3,-0.5) node {3};
\draw (4,-0.5) node {4};
\draw (5,-0.5) node {5};
\draw (6,-0.5) node {6};
\draw (7,-0.5) node {7};
%\draw (8,-0.5) node {8};
%\draw (9,-0.5) node {9};
%\draw (10,-0.5) node {10};
%\draw (11,-0.5) node {11};
%\draw (12,-0.5) node {12};

% node[pos=(0, -0,5)]{0};

\fill (4.5,0.48) circle(3pt);
\fill (6.5,0.48) circle(3pt);

%%caps,cups
\draw[very thick] [-,black,out=90, in=90](0,0.2) to (3,0.2);
\draw[very thick] [-,black,out=90, in=90](1,0.2) to (2,0.2);
\draw[very thick] [-,black,out=90, in=90](4,0.2) to (5,0.2);
\draw[very thick] [-,black,out=90, in=90](6,0.2) to (7,0.2);
%\draw[very thick] [-,black,out=90, in=90](5,0.2) to (8,0.2);
%\draw[very thick] [-,black,out=90, in=90](0,+0.9) to (4,0.2);
%\draw[very thick] [-,black,out=90, in=90](0,+0.9) to (5,0.2);

%\foreach \x in {} \draw + at (-1,0);

\end{tikzpicture} }
\end{center}

The arc diagrams of \cite{ES} do not have labelings $\times$ etc. and contain moreover infinite rays (undotted and dotted) at all other free positions. Ignoring this, the arc diagram of $\bar{\lambda}$ agrees with the cup diagram of \cite{ES}: By \cite{ES}, Proposition 6.1, the dotted arc diagram of $\bar{\lambda}$ is the mirror image of the cup diagram $\underline{\pla^{\owedge}}$. Note that the labelings for weight diagrams differ also from the ones used in \cite{ES}, for instance our $\times$ corresponds to the $\vee$ in loc. cit.

\subsubsection{Example} Consider again the weight diagram $f_{\lambda} =  + \times^3 \times \circ \circ \ldots$. Its associated arc diagram is

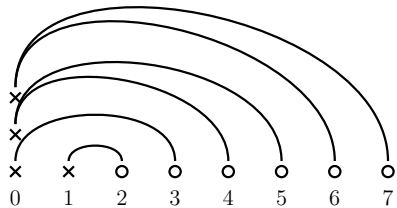
\begin{center}
%\medskip
 
 \scalebox{0.7}{
\begin{tikzpicture}
 %\draw (-1,0) -- (7,0);
\foreach \x in {} %vee
     \draw[very thick] (\x-.1, .1) -- (\x,-0.1) -- (\x +.1, .1);
\foreach \x in {} %wedge
     \draw[very thick] (\x-.1, -.1) -- (\x,0.1) -- (\x +.1, -.1);
\foreach \x in {0,1} %cross
     \draw[very thick] (\x-.1, .1) -- (\x +.1, -.1) (\x-.1, -.1) -- (\x +.1, .1);
%\foreach \x in {1,3,4,7,8,9,10,12}  \draw[semithick] \circ; %circle
     %\draw[very thick]  node at (0,0) [fill=white,draw,circle,inner sep=0pt,minimum size=6pt]{};
     \draw[very thick]  node at (2,0) [fill=white,draw,circle,inner sep=0pt,minimum size=6pt]{};
     \draw[very thick]  node at (3,0) [fill=white,draw,circle,inner sep=0pt,minimum size=6pt]{};
     \draw[very thick]  node at (4,0) [fill=white,draw,circle,inner sep=0pt,minimum size=6pt]{};
     \draw[very thick]  node at (5,0) [fill=white,draw,circle,inner sep=0pt,minimum size=6pt]{};
     \draw[very thick]  node at (6,0) [fill=white,draw,circle,inner sep=0pt,minimum size=6pt]{};
     \draw[very thick]  node at (7,0) [fill=white,draw,circle,inner sep=0pt,minimum size=6pt]{};
     %\draw[very thick]  node at (10,0) [fill=white,draw,circle,inner sep=0pt,minimum size=6pt]{};
     %\draw[very thick]  node at (12,0) [fill=white,draw,circle,inner sep=0pt,minimum size=6pt]{};
\foreach \x in {0} %cross
     \draw[very thick] (\x-.1, +0.8) -- (\x +.1, +0.6) (\x-.1, +0.6) -- (\x +.1, +0.8);
\foreach \x in {0} %cross
     \draw[very thick] (\x-.1, +1.5) -- (\x +.1, +1.3) (\x-.1, +1.3) -- (\x +.1, +1.5);

\draw (0,-0.5) node {0};
\draw (1,-0.5) node {1};
\draw (2,-0.5) node {2};
\draw (3,-0.5) node {3};
\draw (4,-0.5) node {4};
\draw (5,-0.5) node {5};
\draw (6,-0.5) node {6};
7\draw (7,-0.5) node {7};
%\draw (8,-0.5) node {8};
%\draw (9,-0.5) node {9};
%\draw (10,-0.5) node {10};
%\draw (11,-0.5) node {11};
%\draw (12,-0.5) node {12};

% node[pos=(0, -0,5)]{0};

%%caps,cups
\draw[very thick] [-,black,out=90, in=90](0,0.2) to (3,0.2);
\draw[very thick] [-,black,out=90, in=90](1,0.2) to (2,0.2);
%\draw[very thick] [-,black,out=90, in=90](11,0.2) to (12,0.2);
%\draw[very thick] [-,black,out=90, in=90](6,0.2) to (7,0.2);
%\draw[very thick] [-,black,out=90, in=90](5,0.2) to (8,0.2);
\draw[very thick] [-,black,out=90, in=90](0,+0.9) to (4,0.2);
\draw[very thick] [-,black,out=90, in=90](0,+0.9) to (5,0.2);
\draw[very thick] [-,black,out=90, in=90](0,+1.6) to (6,0.2);
\draw[very thick] [-,black,out=90, in=90](0,+1.6) to (7,0.2);

%\foreach \x in {} \draw + at (-1,0);

\end{tikzpicture} }

\end{center}

By Theorem \ref{thmDS1osp} $DS_1(L(\lambda)) \cong L(\lambda_1)$ with $f_{\lambda_1} = \times^2 \times \circ \circ \ldots$. The dotted arc diagram of $\lambda_1$ is given by 

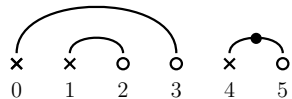
\begin{center}
%\medskip
 
 \scalebox{0.7}{
\begin{tikzpicture}
 %\draw (-1,0) -- (7,0);
\foreach \x in {} %vee
     \draw[very thick] (\x-.1, .1) -- (\x,-0.1) -- (\x +.1, .1);
\foreach \x in {} %wedge
     \draw[very thick] (\x-.1, -.1) -- (\x,0.1) -- (\x +.1, -.1);
\foreach \x in {0,1,4,} %cross
     \draw[very thick] (\x-.1, .1) -- (\x +.1, -.1) (\x-.1, -.1) -- (\x +.1, .1);
%\foreach \x in {1,3,4,7,8,9,10,12}  \draw[semithick] \circ; %circle
     %\draw[very thick]  node at (0,0) [fill=white,draw,circle,inner sep=0pt,minimum size=6pt]{};
     \draw[very thick]  node at (2,0) [fill=white,draw,circle,inner sep=0pt,minimum size=6pt]{};
     \draw[very thick]  node at (3,0) [fill=white,draw,circle,inner sep=0pt,minimum size=6pt]{};
     %\draw[very thick]  node at (4,0) [fill=white,draw,circle,inner sep=0pt,minimum size=6pt]{};
     \draw[very thick]  node at (5,0) [fill=white,draw,circle,inner sep=0pt,minimum size=6pt]{};
     %\draw[very thick]  node at (8,0) [fill=white,draw,circle,inner sep=0pt,minimum size=6pt]{};
     %\draw[very thick]  node at (9,0) [fill=white,draw,circle,inner sep=0pt,minimum size=6pt]{};
     %\draw[very thick]  node at (10,0) [fill=white,draw,circle,inner sep=0pt,minimum size=6pt]{};
     %\draw[very thick]  node at (12,0) [fill=white,draw,circle,inner sep=0pt,minimum size=6pt]{};
%\foreach \x in {0} %cross
 %    \draw[very thick] (\x-.1, +0.8) -- (\x +.1, +0.6) (\x-.1, +0.6) -- (\x +.1, +0.8);

\draw (0,-0.5) node {0};
\draw (1,-0.5) node {1};
\draw (2,-0.5) node {2};
\draw (3,-0.5) node {3};
\draw (4,-0.5) node {4};
\draw (5,-0.5) node {5};
%\draw (6,-0.5) node {6};
%\draw (7,-0.5) node {7};
%\draw (8,-0.5) node {8};
%\draw (9,-0.5) node {9};
%\draw (10,-0.5) node {10};
%\draw (11,-0.5) node {11};
%\draw (12,-0.5) node {12};

% node[pos=(0, -0,5)]{0};

\fill (4.5,0.48) circle(3pt);
%\fill (6.5,0.48) circle(3pt);

%%caps,cups
\draw[very thick] [-,black,out=90, in=90](0,0.2) to (3,0.2);
\draw[very thick] [-,black,out=90, in=90](1,0.2) to (2,0.2);
\draw[very thick] [-,black,out=90, in=90](4,0.2) to (5,0.2);
%\draw[very thick] [-,black,out=90, in=90](6,0.2) to (7,0.2);
%\draw[very thick] [-,black,out=90, in=90](5,0.2) to (8,0.2);
%\draw[very thick] [-,black,out=90, in=90](0,+0.9) to (4,0.2);
%\draw[very thick] [-,black,out=90, in=90](0,+0.9) to (5,0.2);

%\foreach \x in {} \draw + at (-1,0);

\end{tikzpicture} }
\end{center}

i.e. by removing the rightmost dotted cup.

\section{Functors $\Res$ and $\ol{Ind}$}
Let $\fg$ be a  finite-dimensional Lie superalgebra, $\fh\subset\fg$ be a subalgebra
and $z\in\fh$ be an element with the following properties:

(H1) $\fh$ is the centralizer of $\fh_{\ol{0}}$ in $\fg$;

(H2) $\fh_{\ol{0}}$ acts diagonally on $\fg$;

(H3) all eigenvalues of $\ad h$ are real and  $\fg^{\ad h}=\fh$.

In this case we have a usual triangular decomposition
$$\begin{array}{l}
\fg=\fh\oplus (\oplus_{\alpha\in\Delta(\fg)}\fg_{\alpha})\ \ 
 \text{ where }\\
 \Delta(\fg)\subset\fh^*_{\ol{0}},\ \ \Delta(\fg)=\Delta^+(\fg)\coprod \Delta^-(\fg)\\ 
 \fg_{\alpha}:=\{g\in\fg|\ [h,g]=\alpha(h)g\ \text{ for all }h\in\fh_{\ol{0}}\}\\
\Delta^+(\fg):=\{\alpha\in\Delta(\fg)|\  \alpha(h)>0\},\ \ \ \Delta^-(\fg):=\{\alpha\in\Delta(\fg)|\ \alpha(h)<0\}.
\end{array}$$
 We set
$\fn^{\pm}:=\oplus_{\alpha\in\Delta^{\pm}}\fg_{\alpha},\ \ \ \fb:=\fn\oplus\fn^+$.
We consider the  partial order on $\fh^*_{\ol{0}}$ given by
$$\lambda>\nu\ \ \text{  if } \nu-\lambda\in\mathbb{N}\Delta^-.$$

\subsection{}\label{funcnot}
We fix  $z\in\fh^*_{\ol{0}}$ satisfying
\begin{equation}\label{condz}
\alpha(z)\in\mathbb{R}_{\geq 0}\ \text{ for all }\alpha\in\Delta^+\ \text{ and }
\alpha(z)\in\mathbb{R}_{\leq 0}\ \text{ for all }\alpha\in\Delta^-.
\end{equation}
and introduce
$$\ft:=\fg^z,\ \ \ \fm:= \displaystyle\bigoplus_{\alpha\in\Delta: \alpha(z)>0}\fg_{\alpha},\ \ \ 
\fp=\fg^z+\fb=\ft\rtimes\fm.$$
The triples $(\fp(z),\fh,h)$ and $(\ft(z),\fh,h)$ satisfy (H1)--(H3) and
$$\begin{array}{l}
\Delta^+(\fp)=\Delta^+(\fg),\ \ \ \ 
\ \Delta^+(\ft)=\{\alpha\in\Delta^{+}(\fg)|\ \alpha(z)=0\}\\
\Delta^-(\fp)=\Delta^-(\ft)=\{\alpha\in\Delta^{-}(\fg)|\ \alpha(z)=0\}
\end{array}$$

\subsection{Functors $\Res_a$ and  $\ol{\Ind}$}\label{functdef}

We denote by $\CO$ the full category of 
finitely generated modules with a diagonal action of $\fh_{\ol{0}}$
and locally nilpotent action of $\fn$.
It is easy to see that,
up to a parity change, the simple $\fh$-modules are 
parametrized by
$\lambda\in\fh^*_{\ol{0}}$; we denote by $C_{\lambda}$ a simple
$\fh$-module, where $\fh_{\ol{0}}$ acts by $\lambda$.
We view $C(\lambda)$ as a $\fb$-module with the zero action of $\fn$
and set
$$M(\lambda):=\Ind^{\fg}_{\fb} C(\lambda);$$  
 this module
has a  unique simple quotients which we denote by 
$L(\lambda)$. 
(The module $M(\lambda)$ is a Verma module
if $\fh^*_{\ol{1}}=0$).

For each $a\in\mathbb{C}$ we define
a functor $\Res_a: \CO(\fg)\to\CO(\ft)$ by
$$\Res_a(N):=\{v\in N|\ \  zv=a v\}\ \text{ for } N\in\CO(\fg).$$

We assume that 

{\em each module  in $\CO$ admits a unique maximal finite-dimensional
quotient}.

This holds, for example,  if $\fg_{\ol{0}}$ is reductive,
since in this case all modules in $\CO$ have finite lengths. 

View $V\in\CO(\ft)$
as a $\fp$-module with the trivial action of $\fm$ and
consider the induced $\fg$-module
$$\Ind(V)=\cU(\fg)\otimes_{\cU(\fp)} V.$$
We denote by $\ol{\Ind}(V)$ the maximal finite-dimensional quotient of
$\Ind(V)$.

For  a semisimple  $\fh_{\ol{0}}$-module $N$ we
 denote by $N_{\nu}$ the weight space of the weight $\nu$ and by
 $\Omega(N)$ the set of weights of $N$.
 
\subsection{Assumption}
We fix a block $\cB$ in the category $\CO$ and set
$$\begin{array}{l}
P^+(\cB)_0:=\{\lambda \in\fh^*_{\ol{0}}|\ \dim L(\lambda)<\infty,\ \ 
L(\lambda)\in\cB\ \text{ or }\ \Pi(L(\lambda))\in\cB\}.
\end{array}$$

We fix $\fh''\subset  \fh_{\ol{0}}$ such that
$$z\in\fh'',\ \ \ [\fh'',\ft]=0$$

and $\mu\in (\fh'')^*$ satisfying 

\begin{equation}\label{assmstable}
\forall \lambda,\nu\in P^+(\cB)\ \ \lambda\geq \nu,\ \ 
\lambda|_{\fh''}=\mu\ \Longrightarrow\ \ \nu|_{\fh''}=\mu.
\end{equation}

We set
$$A:=\{\lambda\in  P^+(\cB)|\ \lambda|_{\fh''}=\mu\}.$$
By above, if $\lambda\in A$ and $\nu\in P^+(\cB)$ are such that $\lambda\geq \nu$,
then $\nu\in A$.

Let $\cF_{\fg}(A)$ (resp., $\cF_{\ft}(A)$) be the Serre subcategory of $\CO(\fg))$ (resp.,
of $\CO(\ft)$)
generated by the simple modules
$\{L(\lambda),\Pi(L(\lambda))\}_{\lambda\in A}$ 
(resp., by $\{L_{\ft}(\lambda),\Pi(L_{\ft}(\lambda))\}_{\lambda\in A}$).

\subsection{}
\begin{prop}{propind}
Set $a:=\mu(z)$ and $\Res:=\Res_a$.

(i) For $N\in \cF_{\fg}(A)$ one has
$\Res(N)=N^{\fm}$. For  $\lambda\in A$ one has
$\Res(L(\lambda))=L_{\ft}(\lambda)$.

(ii) The restrictions of $\Res$ and $\ol{\Ind}$ give an equivalence of categories
 $\cF_{\fg}(A)$ and $\cF_{\ft}(A)$.
\end{prop}
\begin{proof}
For each $\fh$-module $M$ we denote by $\Spec(M)$ the set of
$z$-eigenvalues. Recall that
$$\Spec(\fm)\subset\mathbb{R}_{>0},\ \ \Spec(\fn^-)\subset\mathbb{R}_{\leq 0}.$$
Identifying $V$ with $1\otimes V\subset\Ind(V)$ we obtain
\begin{equation}\label{resind}
\Res(\Ind(V))=V\ \text{ for } V\in\CO(\ft).
\end{equation}

Take $N\in \cF_{\fg}(A)$.
For $\lambda\in A$ one has
$\Spec(L(\lambda))\subset (\mu(z)-\mathbb{R}_{\geq 0})=:U$; this gives
$\Spec(N)\subset U$ and implies
$$\Res(N)\subset N^{\fm}.$$

Set $\fn':=\fn\cap \ft$. Let $v\in N^{\fm}$ be a non-zero vector of weight $\nu$.
The subspace $\cU(\fn)v$ contains a singular weight 
vector $v'$; let $\lambda'$ be the weight of $v'$.
Since $N\in \cF_{\fg}(A)$ we have $\lambda'\in A$, that is $\lambda'(z)=\mu(z)$. 
Since $\fm v=0$ we have $\cU(\fn)v=\cU(\fn')v$, so $v'\in \cU(\fn')v$,
that is $\nu(z)=\lambda'(z)$. Therefore $v\in \Res(N)$. 
Hence $\Res(N)=N^{\fm}$.

Take $\lambda\in A$ and denote by $v_{\lambda}$ a highest weight vector in $L(\lambda)$. By above, for each non-zero vector $v$ the space
$\cU(\fn')v$  contains $v_{\lambda}$. Hence
$L(\lambda)^{\fm}$ is simple, so $L(\lambda)^{\fm}=L_{\ft}(\lambda)$.
 This establishes (i).

Clearly, $\Res$ is exact.
By (i), $\Res$ maps simple modules in $\cF_{\fg}(A)$ to simple modules
in $\cF_{\ft}(A)$. Therefore
the restriction of $\Res$ gives an exact functor
$$\Res:\cF_{\fg}(A)\to \cF_{\ft}(A).$$

Take a module $V\in\cF_{\ft}(A)$. 
Note that for $\nu\in A$ the module $\Ind L_{\ft}(\nu)$
lies in $\cB$; by~(\ref{assmstable})
each simple finite-dimensional subquotient
of $\Ind L_{\ft}(\nu)$ is of the form $L(\nu')$, where $\nu'\in A$.
Since $\Ind$ is exact, each simple subquotient of $\Ind(V)$ is $L(\nu')$
or $\Pi(L(\nu'))$ for some $\nu'\in A$, so
  $\ol{\Ind}(V)\in\cF_{\fg}(A)$. Moreover, using the exactness
of $\Res$ and~(\ref{resind}) we get
$$[\Ind(V):L({\nu})]=[V:L_{\ft}({\nu})],$$
that is
\begin{equation}\label{multresi}
[\ol{\Ind}(V):L(\nu)]\leq [V:L_{\ft}(\nu)].\end{equation}

The module $\Ind(L_{\ft}(\lambda))$
is a quotient of  $M(\lambda)$, so 
$L(\lambda)$ is a quotient of $\Ind(L_{\ft}(\lambda))$.
Using~(\ref{multresi}) we get
$$\ol{\Ind}(L_{\ft}(\lambda))=L(\lambda).$$ 

For each $V\in\cF_{\ft}(A)$ and $N\in\cF_{\fg}(A)$ we have
$$\Hom_{\fg}(\ol{\Ind}(V), N)=\Hom_{\fg}(\Ind(V), N)=\Hom_{\fp}(V,N)=
\Hom_{\ft}(V,N^{\fm}).$$
Using (i) we conclude that
 $\ol{\Ind}:\cF_{\ft}(A)\to\cF_{\fg}(A)$ is a left adjoint to $\Res:\cF_{\fg}(A)\to\cF_{\ft}(A)$; by above,
these functors map simple modules to simple modules and $\Res$ is exact. 

Take any $N\in\cF_{\fg}(A)$ and set $V:=\Res(N)$. Let $\phi\in \Hom_{\fg}(\ol{\Ind}(V), N)$ be the preimage
of the identity map
$$V\iso \Res(N).$$
The image of $\phi$ is the submodule of $N$
generated by $\Res(N)$; since $\Res$ is exact and $\Res(M)\not=0$
for each $M\in\cF(A)$, $\phi$ is surjective. Moreover,
for each $\nu\in A$ one has
$$[N:L(\nu)]=[V:L_{\ft}(\nu)].$$
Combining with~(\ref{multresi}) we conclude that $\phi$ is bijective and
for each $\nu\in A$ one has
$$[\ol{\Ind}(V):L(\nu)]=[V:L_{\ft}\nu)]$$
which gives
\begin{equation}\label{multres}
[\Res(\ol{Ind}(V)):L_{\ft}(\nu)]=[V:L_{\ft}(\nu)].\end{equation}
Take $V\in\cF_{\ft}(A)$. Identifying $V$ and $1\otimes V\subset N$
we have $V=\Res(\Ind V)$. Since $\Res$ is an exact functor on $\CO$
this gives the natural surjective map $V\to \Res(\ol{\Ind}(V))$.
By~(\ref{multres}) this map is bijective.
Hence $\Res$ and $\ol{\Ind}$ provide an equivalence of categories
$\cF_{\fg}(A)$ and $\cF_{\ft}(A)$.
\end{proof}

%\subsubsection{Case $\Sigma'=\emptyset$}\label{RemRes} ???
%Consider the special case of the above construction
%when $\fh''=\fh$ and $\upsilon''=Id$. Therefore
%$A=\{\mu\}$ and each module in $\cF(A)$ is of the form
%$L(\mu)^{\oplus s}$ for $s\geq 0$. For $N\in\cF(A)$ the map
%$\Res$ is given by  $\Res(N)=N_{\mu}$
%and is an equivalence of categories between $\cF(A)$ and
%the category of even finite-dimensional vector spaces.
%
%Since $A'=0$, the category $\cF'(A')$
%is the category of even finite-dimensional vector spaces.
%Hence $\Res$ is an equivalence of categories.

\subsection{Remark}\label{gequivl}
Consider the case $\fh=\fh_{\ol{0}}$ and
$\ft=\fl\times\fh''$. Set $\fh':=\fl\cap\fh$ and
$$A':=\{\lambda|_{\fh'}|\ \lambda\in A\}.$$
Since $\fh=\fh'\oplus \fh''$ the map
$\lambda\mapsto \lambda|_{\fh'}$ gives a bijection between $A$ and $A'$.
Using~\Prop{propind} we obtain
the equivalence of the category $\cF_{\fg}(A)$ and 
the Serre category $\cF_{\fl}(A')$ which is generated by 
the simple $\fl$-modules $\{L_{\fl}(\lambda'),\Pi(L_{\fl}(\lambda'))\}_{\lambda'\in A'}$.

%
%
%This construction provides an equivalence between the subcategory of stable modules with a given central character and a certain Serre subcategory of the corresponding principal block. For $\osp(M|N)$-case we give details in Section~\ref{notres} and show that this equivalence
%commutes with $\DS$-functors. 
%

\end{document}